\documentclass[10pt,12pt]{article}%
\usepackage{amsmath}
\usepackage{amsfonts}
\usepackage{amssymb}
\usepackage{amsthm}
\usepackage{hyperref}
\usepackage{graphicx}%
\providecommand{\U}[1]{\protect\rule{.1in}{.1in}}
\newtheorem{thm}{Theorem}[section]
\newtheorem{lemm}{Lemma}[section]
\newtheorem*{appen}{Lemma A}
\newtheorem{step}{Step}[subsection]
\newtheorem{parti}{Part}[section]
\begin{document}

\title{Carleman commutator approach in\\logarithmic convexity for parabolic equations}
\author{Kim Dang Phung\thanks{ Universit\'{e} d'Orl\'{e}ans, Laboratoire MAPMO, CNRS
UMR 7349, F\'{e}d\'{e}ration Denis Poisson, FR CNRS 2964, B\^{a}timent de
Math\'{e}matiques, B.P. 6759, 45067 Orl\'{e}ans Cedex 2, France. E-mail
address: kim\_dang\_phung@yahoo.fr. }}
\date{}
\maketitle

\begin{abstract}
In this paper we investigate on a new strategy combining the logarithmic
convexity (or frequency function) and the Carleman commutator to obtain an
observation estimate at one time for the heat equation in a bounded domain. We
also consider the heat equation with an inverse square potential. Moreover,
spectral inequality for the associated eigenvalue problem is derived.

\end{abstract}

\bigskip

\section{Introduction and main results}

\bigskip

\bigskip

When we mention the logarithmic convexity method for the heat equation in a
bounded domain $\Omega\subset\mathbb{R}^{n}$ :
\[
\left\{
\begin{array}
[c]{ll}%
{\partial}_{t}u-\Delta u=0\text{ ,} & \quad\text{in}~\Omega\times\left(
0,T\right)  \text{ ,}\\
u=0\text{ ,} & \quad\text{on}~\partial\Omega\times\left(  0,T\right)  \text{
,}\\
u\left(  \cdot,0\right)  =u_{0}\in L^{2}(\Omega)\left\backslash \left\{
0\right\}  \right.  \text{ , } &
\end{array}
\right.
\]
we have in mind that $t\mapsto$ln$\left\Vert u\left(  \cdot,t\right)
\right\Vert _{L^{2}\left(  \Omega\right)  }^{2}$ is a convex function by
evaluating the sign of the derivative of $t\mapsto\displaystyle\frac
{\int_{\Omega}\left\vert \nabla u\left(  x,t\right)  \right\vert ^{2}dx}%
{\int_{\Omega}\left\vert u\left(  x,t\right)  \right\vert ^{2}dx}$ (see
\cite{AN}, \cite[p.11]{Pa}, \cite[p.43]{I}, \cite{Ve}). As a consequence, the
following well-known estimate holds. For any $0\leq t\leq T$,
\[
\left\Vert e^{t\Delta}u_{0}\right\Vert _{L^{2}\left(  \Omega\right)  }%
\leq\left\Vert e^{T\Delta}u_{0}\right\Vert _{L^{2}\left(  \Omega\right)
}^{t/T}\left\Vert u_{0}\right\Vert _{L^{2}\left(  \Omega\right)  }%
^{1-t/T}\text{ .}%
\]

\bigskip

In a series of articles (see \cite{PW1}, \cite{PW2}, \cite{PWZ}, \cite{BP} for
parabolic equations) inspired by \cite{Po} and \cite{EFV}, we were interested
on the function $t\mapsto\int_{\Omega}\left\vert u\left(  x,t\right)
\right\vert ^{2}e^{\Phi\left(  x,t\right)  }dx$ and its frequency function
$t\mapsto\displaystyle\frac{\int_{\Omega}\left\vert \nabla u\left(
x,t\right)  \right\vert ^{2}e^{\Phi\left(  x,t\right)  }dx}{\int_{\Omega
}\left\vert u\left(  x,t\right)  \right\vert ^{2}e^{\Phi\left(  x,t\right)
}dx}$ when $e^{\Phi\left(  x,t\right)  }=\frac{1}{\left(  T-t+\hbar\right)
^{n/2}}e^{-\frac{\left\vert x-x_{0}\right\vert ^{2}}{4\left(  T-t+\hbar
\right)  }}$ with $x_{0}\in\Omega$, $\hbar>0$. It provides us with an
observation estimate at one point in time: For any $T>0$ and any $\omega$
nonempty open subset of $\Omega$,
\[
\left\Vert e^{T\Delta}u_{0}\right\Vert _{L^{2}\left(  \Omega\right)  }%
\leq\left(  ce^{\frac{K}{T}}\left\Vert e^{T\Delta}u_{0}\right\Vert
_{L^{2}\left(  \omega\right)  }\right)  ^{\beta}\left\Vert u_{0}\right\Vert
_{L^{2}\left(  \Omega\right)  }^{1-\beta}\text{ .}%
\]
Here $c,K>0$ and $\beta\in\left(  0,1\right)  $. From the above observation at
one time, many applications were derived as bang-bang control \cite{PW2} and
impulse control \cite{PWX}, fast stabilization \cite{PWX} or local backward
reconstruction \cite{Vo}. In particular, we can also deduce the observability
estimate for parabolic equations on a positive measurable set in time
\cite{PW2}. Recall that observability for parabolic equations have a long
history now from the works of \cite{LR} and \cite{FI} based on Carleman
inequalities. Furthermore, it was remarked in \cite{AEWZ} that the observation
estimate at one point in time is equivalent to the Lebeau-Robbiano spectral
inequality on the sum of eigenfunctions of the Dirichlet Laplacian. Recall
that the Lebeau-Robbiano spectral inequality, originally derived from Carleman
inequalities for elliptic equations (see \cite{JL}, \cite{LRL}, \cite{Lu}),
was used in different contexts as in thermoelasticity (see \cite{LZ},
\cite{BN}), for the Stokes operator \cite{CL}, in transmission problem and
coupled systems (see \cite{Le}, \cite{LLR}), for the Bilaplacian (see
\cite{Ga}, \cite{EMZ}, \cite{LRR3}), in Kolmogorov equation (see \cite{LRM},
\cite{Z}). We also refer to \cite{M}.

\bigskip

In this paper, we study the equation solved by $f\left(  x,t\right)  =u\left(
x,t\right)  e^{\frac{1}{2}\Phi\left(  x,t\right)  }$ for a larger set of
weight functions $\Phi\left(  x,t\right)  $ and establish a kind of convexity
property for $t\mapsto$ln$\left\Vert f\left(  \cdot,t\right)  \right\Vert
_{L^{2}\left(  \Omega\right)  }^{2}$. By such approach we make appear the
Carleman commutator. The link between logarithmic convexity (or frequency
function) and Carleman inequality has already appeared in \cite{EKPV1} (see
also \cite{EKPV2}, \cite{EKPV3}).

\bigskip

Choosing suitable weight functions $\Phi$ (not necessary linked to the heat
kernel) we obtain the following new results:

\bigskip

\begin{thm}
\label{theorem1.1} Suppose that $\Omega\subset\mathbb{R}^{n}$ is a convex
domain or a star-shaped domain with respect to $x_{0}\in\Omega$ such that
$\left\{  x;\left\vert x-x_{0}\right\vert <r\right\}  \Subset\Omega$ for some
$r\in\left(  0,1\right)  $. Then for any $u_{0}\in L^{2}\left(  \Omega\right)
$,\thinspace$T>0$, $\left(  a_{j}\right)  _{j\geq1}\in\mathbb{R}$, $\lambda
>0$, $\varepsilon\in\left(  0,1\right)  $, one has%
\[
\left\Vert e^{T\Delta}u_{0}\right\Vert _{L^{2}\left(  \Omega\right)  }\leq
K_{\varepsilon}e^{\frac{K_{\varepsilon}}{r^{\varepsilon}}\frac{1}{T}}\int
_{0}^{T}\left\Vert e^{t\Delta}u_{0}\right\Vert _{L^{2}\left(  \left\vert
x-x_{0}\right\vert <r\right)  }dt\text{ }%
\]
and%
\[
\sum\limits_{\lambda_{j}\leq\lambda}\left\vert a_{j}\right\vert ^{2}%
\leq4e^{\frac{K_{\varepsilon}}{r^{\varepsilon}}\sqrt{\lambda}}{\int
\nolimits_{\left\vert x-x_{0}\right\vert <r}}\left\vert \sum\limits_{\lambda
_{j}\leq\lambda}a_{j}e_{j}\left(  x\right)  \right\vert ^{2}dx
\]
where $K_{\varepsilon}>0$ is a constant only depending on $\left(
\varepsilon,\normalfont{\text{max}}\left\{  \left\vert x-x_{0}\right\vert
;x\in\overline{\Omega}\right\}  \right)  $. Here $\left(  \lambda_{j}%
,e_{j}\right)  $ denotes the eigenbasis of the Laplace operator with Dirichlet
boundary condition.
\end{thm}

\bigskip

Theorem \ref{theorem1.1} thus states both the observability for the heat
equation and the spectral inequality for the Dirichlet Laplacian in a simple
geometry. One can see how fast the constant cost blows up when the observation
region $\omega$ becomes smaller. Notice that the constant $K_{\varepsilon}$
does not depend on the dimension $n$ (see \cite[Theorem 4.2]{BP}).

\bigskip

\begin{thm}
\label{theorem1.2} Let $n\geq3$ and consider a $C^{2}$ bounded domain
$\Omega\subset\mathbb{R}^{n}$ such that $0\in\Omega$. Let $\omega\subset
\Omega$ be a nonempty open set. Suppose that
\[
\mu\leq\left\vert
\begin{array}
[c]{ll}%
\frac{7}{2\cdot3^{3}}\text{ ,} & \text{if }n=3\text{ ,}\\
\frac{1}{4}\left(  n-1\right)  \left(  n-3\right)  \text{ ,} & \text{if }%
n\geq4\text{ .}%
\end{array}
\right.
\]
Then, there exist constants $c>0$, $K>0$ such that for any $\left(
a_{j}\right)  _{j\geq1}\in\mathbb{R}$ and any $\lambda>0$, we have
\[
\sqrt{\sum\limits_{\lambda_{j}\leq\lambda}\left\vert a_{j}\right\vert ^{2}%
}\leq ce^{K\sqrt{\lambda}}{\int\nolimits_{\omega}}\left\vert \sum
\limits_{\lambda_{j}\leq\lambda}a_{j}e_{j}\left(  x\right)  \right\vert dx
\]
where $\left(  \lambda_{j},e_{j}\right)  $ denotes the eigenbasis of the
Schr\"{o}dinger operator $-\Delta-\frac{\mu}{\left\vert x\right\vert ^{2}}$
with Dirichlet boundary condition
\[
\left\{
\begin{array}
[c]{rll}%
-\Delta e_{j}-\frac{\mu}{\left\vert x\right\vert ^{2}}e_{j} & =\lambda
_{j}e_{j}\text{ ,} & \text{ in }\Omega\text{ ,}\\
e_{j} & =0\text{ ,} & \text{ on }\partial\Omega\text{ .}%
\end{array}
\right.
\]

\end{thm}

\bigskip

Theorem \ref{theorem1.2} gives a spectral inequality for the Schr\"{o}dinger
operator $-\Delta-\frac{\mu}{\left\vert x\right\vert ^{2}}$ under a quite
strong assumption on $\mu<\mu^{\ast}$ where the critical coefficient is
$\mu^{\ast}=\frac{1}{4}\left(  n-2\right)  ^{2}$. Our first motivation was to
be able to choose $0\notin\overline{\omega}$ by performing localization with
annulus. We believe that a similar analysis can be handle with more suitable
weight function $\Phi$ than those considered here and may considerably improve
the results presented here.

\bigskip

We have organized our paper as follows. Section 2 is the important part of
this article. We present the strategy to get the observation at one point by
studying the equation solved by $f=ue^{\Phi/2}$ for a larger set of weight
functions $\Phi$ adapting the energy estimates style of computations in
\cite{BT} (see also \cite[Section 4]{BP}). The Carleman commutator appears
naturally here. Section 3 is devoted to check different possibilities for the
weight function $\Phi$, and in particular for the localization with annulus.
In Section 4, we prove Theorem \ref{theorem1.1}. The proof of Theorem
\ref{theorem1.2} is given in Section 5. In Appendix, we recall the useful link
between the observation at one point and the spectral inequality.

\bigskip

I am happy to dedicate this paper to my friend and colleague Jiongmin Yong on
the occasion of his 60th birthday. I am also grateful for his book \cite{LY}
in where I often found the answer on my questions.

\bigskip

\section{The strategy of logarithmic convexity with the Carleman commutator}

\bigskip

We present an approach to get the observation estimate at one point in time
for a model heat equation in a bounded domain $\Omega\subset\mathbb{R}^{n}$
with Dirichlet boundary condition. We shall present this strategy
step-by-step. Two different geometric cases are discussed: When $\Omega$ is
convex or star-shaped, we can used a global weight function; For the more
general $C^{2}$ domain $\Omega$, we will use localized weight functions
exploiting a covering argument and propagation of interpolation inequalities
along a chain of balls (also called propagation of smallness).

\bigskip

\subsection{Convex domain}

\bigskip

Throughout this subsection, we assume that $\Omega\subset\mathbb{R}^{n}$ is a
convex domain or a star-shaped domain with respect to $x_{0}\in\Omega$. Let
$\left\langle \cdot,\cdot\right\rangle $ denote the usual scalar product in
$L^{2}\left(  \Omega\right)  $ and let $\left\Vert \cdot\right\Vert $ be its
corresponding norm. Here, recall that $u\left(  x,t\right)  =e^{t\Delta}%
u_{0}\left(  x\right)  \in C\left(  \left[  0,T\right]  ;L^{2}\left(
\Omega\right)  \right)  \cap C\left(  \left(  0,T\right]  ;H_{0}^{1}\left(
\Omega\right)  \right)  $ and we aim to check that
\[
\left\Vert u\left(  \cdot,T\right)  \right\Vert \leq\left(  ce^{\frac{K}{T}%
}\left\Vert u\left(  \cdot,T\right)  \right\Vert _{L^{2}\left(  \omega\right)
}\right)  ^{\beta}\left\Vert u\left(  \cdot,0\right)  \right\Vert ^{1-\beta
}\text{ .}%
\]
The strategy to establish the above observation at one time is as follows. We
decompose the proof into six steps.

\begin{step}
\label{step2.1.1} Symmetric part and antisymmetric part.
\end{step}

Let $\Phi$ be a sufficiently smooth function of $\left(  x,t\right)
\in\mathbb{R}^{n}\times\mathbb{R}_{t}$ and define
\[
f\left(  x,t\right)  =u\left(  x,t\right)  e^{\Phi\left(  x,t\right)
/2}\text{ .}%
\]
We look for the equation solved by $f$ by computing $e^{\Phi\left(
x,t\right)  /2}\left(  \partial_{t}-\Delta\right)  \left(  e^{-\Phi\left(
x,t\right)  /2}f\left(  x,t\right)  \right)  $. We find that
\[
\partial_{t}f-\Delta f-\frac{1}{2}f\left(  \partial_{t}\Phi+\frac{1}%
{2}\left\vert \nabla\Phi\right\vert ^{2}\right)  +\nabla\Phi\cdot\nabla
f+\frac{1}{2}\Delta\Phi f=0\quad\text{in}~\Omega\times\left(  0,T\right)
\text{ ,}%
\]
and furthermore, $f_{\left\vert \partial\Omega\right.  }=0$. Introduce
\[
\left\{
\begin{array}
[c]{ll}%
\mathcal{A}f=-\nabla\Phi\cdot\nabla f-\frac{1}{2}\Delta\Phi f\text{ ,} & \\
\mathcal{S}f=\Delta f+\eta f\text{ where }\eta=\frac{1}{2}\left(  \partial
_{t}\Phi+\frac{1}{2}\left\vert \nabla\Phi\right\vert ^{2}\right)  \text{ .} &
\end{array}
\right.
\]
We can check that
\[
\left\{
\begin{array}
[c]{ll}%
\left\langle \mathcal{A}f,g\right\rangle =-\left\langle \mathcal{A}%
g,f\right\rangle \text{ ,} & \\
\left\langle \mathcal{S}f,g\right\rangle =\left\langle \mathcal{S}%
g,f\right\rangle \text{ for any }g\in H_{0}^{1}\left(  \Omega\right)  \text{
.} &
\end{array}
\right.
\]
Furthermore, we have
\[
\partial_{t}f-\mathcal{S}f-\mathcal{A}f=0\text{ .}%
\]

\begin{step}
\label{step2.1.2} Energy estimates.
\end{step}

\nopagebreak Multiplying by $f$ the above equation, integrating over $\Omega$,
we obtain that%
\[
\frac{1}{2}\frac{d}{dt}\left\Vert f\right\Vert ^{2}+\left\langle
-\mathcal{S}f,f\right\rangle =0\text{ .}%
\]
Introduce the frequency function $t\mapsto\mathbf{N}\left(  t\right)  $
defined by%
\[
\mathbf{N}=\frac{\left\langle -\mathcal{S}f,f\right\rangle }{\left\Vert
f\right\Vert ^{2}}\text{ .}%
\]
Thus,
\[
\frac{1}{2}\frac{d}{dt}\left\Vert f\right\Vert ^{2}+\mathbf{N}\left\Vert
f\right\Vert ^{2}=0\text{ .}%
\]
Now, we compute the derivative of $\mathbf{N}$ and claim that:%
\[
\frac{d}{dt}\mathbf{N}\leq\frac{1}{\left\Vert f\right\Vert ^{2}}\left\langle
-\left(  \mathcal{S}^{\prime}+\left[  \mathcal{S},\mathcal{A}\right]  \right)
f,f\right\rangle -\frac{1}{\left\Vert f\right\Vert ^{2}}\int_{\partial\Omega
}\partial_{\nu}f\mathcal{A}fd\sigma\text{ .}%
\]
Indeed,
\[%
\begin{array}
[c]{ll}
& \displaystyle\frac{d}{dt}\mathbf{N}\\
& =\displaystyle\frac{1}{\left\Vert f\right\Vert ^{4}}\left(
\displaystyle\frac{d}{dt}\left\langle -\mathcal{S}f,f\right\rangle \left\Vert
f\right\Vert ^{2}+\left\langle \mathcal{S}f,f\right\rangle \displaystyle\frac
{d}{dt}\left\Vert f\right\Vert ^{2}\right) \\
& =\displaystyle\frac{1}{\left\Vert f\right\Vert ^{2}}\left[  \left\langle
-\mathcal{S}^{\prime}f,f\right\rangle -2\left\langle \mathcal{S}f,f^{\prime
}\right\rangle \right]  +\displaystyle\frac{2}{\left\Vert f\right\Vert ^{4}%
}\left\langle \mathcal{S}f,f\right\rangle ^{2}\\
& =\displaystyle\frac{1}{\left\Vert f\right\Vert ^{2}}\left[  \left\langle
-\mathcal{S}^{\prime}f,f\right\rangle -2\left\langle \mathcal{S}%
f,\mathcal{A}f\right\rangle \right]  +\displaystyle\frac{2}{\left\Vert
f\right\Vert ^{4}}\left[  -\left\Vert \mathcal{S}f\right\Vert ^{2}\left\Vert
f\right\Vert ^{2}+\left\langle \mathcal{S}f,f\right\rangle ^{2}\right] \\
& =\displaystyle\frac{1}{\left\Vert f\right\Vert ^{2}}\left[  \left\langle
-\left(  \mathcal{S}^{\prime}+\left[  \mathcal{S},\mathcal{A}\right]  \right)
f,f\right\rangle -\displaystyle\int_{\partial\Omega}\partial_{\nu}%
f\mathcal{A}fd\sigma\right]  +\displaystyle\frac{2}{\left\Vert f\right\Vert
^{4}}\left[  -\left\Vert \mathcal{S}f\right\Vert ^{2}\left\Vert f\right\Vert
^{2}+\left\langle \mathcal{S}f,f\right\rangle ^{2}\right] \\
& \leq\displaystyle\frac{1}{\left\Vert f\right\Vert ^{2}}\left[  \left\langle
-\left(  \mathcal{S}^{\prime}+\left[  \mathcal{S},\mathcal{A}\right]  \right)
f,f\right\rangle -\displaystyle\int_{\partial\Omega}\partial_{\nu}%
f\mathcal{A}fd\sigma\right]  \text{ .}%
\end{array}
\]
In the third line, we used $\frac{1}{2}\frac{d}{dt}\left\Vert f\right\Vert
^{2}+\left\langle -\mathcal{S}f,f\right\rangle =0$; In the fourth line,
multiplying the equation of $f$ by $\mathcal{S}f$, and integrating over
$\Omega$, give $\left\langle \mathcal{S}f,f^{\prime}\right\rangle =\left\Vert
\mathcal{S}f\right\Vert ^{2}+\left\langle \mathcal{S}f,\mathcal{A}%
f\right\rangle $; In the fifth line, $\partial_{\nu}$ denotes the normal
derivative to the boundary, and we used
\[%
\begin{array}
[c]{ll}%
2\left\langle \mathcal{S}f,\mathcal{A}f\right\rangle  & =\left\langle
\mathcal{SA}f,f\right\rangle -\left\langle \mathcal{AS}f,f\right\rangle
+\displaystyle\int_{\partial\Omega}\partial_{\nu}f\mathcal{A}fd\sigma\\
& :=\left\langle \left[  \mathcal{S},\mathcal{A}\right]  f,f\right\rangle
+\displaystyle\int_{\partial\Omega}\partial_{\nu}f\mathcal{A}fd\sigma\text{ ;}%
\end{array}
\]
In the sixth line, we used Cauchy-Schwarz inequality.

Here we have followed the energy estimates style of computations in \cite{BT}
(see also \cite[p.535]{Ph}) The interested reader may wish here to compare
with \cite[Theorem 3]{EKPV1}.

\begin{step}
\label{step2.1.3} Assumption on Carleman commutator.
\end{step}

Assume that $\displaystyle\int_{\partial\Omega}\partial_{\nu}f\mathcal{A}%
fd\sigma\geq0$ on $\left(  0,T\right)  $ by convexity or star-shaped property
of $\Omega$, and suppose that
\[
\left\langle -\left(  \mathcal{S}^{\prime}+\left[  \mathcal{S},\mathcal{A}%
\right]  \right)  f,f\right\rangle \leq\frac{1}{\Upsilon}\left\langle
-\mathcal{S}f,f\right\rangle \text{ }%
\]
on $\left(  0,T\right)  $ where $\Upsilon\left(  t\right)  =T-t+\hbar$ and
$\hbar>0$. Therefore the following differential inequalities hold.
\[
\left\{
\begin{array}
[c]{ll}%
\displaystyle\frac{1}{2}\frac{d}{dt}\left\Vert f\left(  \cdot,t\right)
\right\Vert ^{2}+\mathbf{N}\left(  t\right)  \left\Vert f\left(
\cdot,t\right)  \right\Vert ^{2}=0\text{ ,} & \\
\displaystyle\frac{d}{dt}\mathbf{N}\left(  t\right)  \leq\displaystyle\frac
{1}{\Upsilon\left(  t\right)  }\mathbf{N}\left(  t\right)  \text{ .} &
\end{array}
\right.
\]
By solving such system of differential inequalities, we obtain (see
\cite[p.655]{BP}): For any $0<t_{1}<t_{2}<t_{3}\leq T$,%
\[
\left(  \left\Vert f\left(  \cdot,t_{2}\right)  \right\Vert ^{2}\right)
^{1+M}\leq\left(  \left\Vert f\left(  \cdot,t_{1}\right)  \right\Vert
^{2}\right)  ^{M}\left\Vert f\left(  \cdot,t_{3}\right)  \right\Vert ^{2}%
\]
where%
\[
M=\frac{-\text{ln}\left(  T-t_{3}+\hbar\right)  +\text{ln}\left(
T-t_{2}+\hbar\right)  }{-\text{ln}\left(  T-t_{2}+\hbar\right)  +\text{ln}%
\left(  T-t_{1}+\hbar\right)  }\text{ .}%
\]
In other words, we have
\[
\left(  \displaystyle\int_{\Omega}\left\vert u\left(  x,t_{2}\right)
\right\vert ^{2}e^{\Phi\left(  x,t_{2}\right)  }dx\right)  ^{1+M}\leq\left(
\displaystyle\int_{\Omega}\left\vert u\left(  x,t_{1}\right)  \right\vert
^{2}e^{\Phi\left(  x,t_{1}\right)  }dx\right)  ^{M}\displaystyle\int_{\Omega
}\left\vert u\left(  x,t_{3}\right)  \right\vert ^{2}e^{\Phi\left(
x,t_{3}\right)  }dx\text{ .}%
\]
\textit{ }

\begin{step}
\label{step2.1.4}
\end{step}

Let $\omega$ be a nonempty open subset of $\Omega$. We take off the weight
function $\Phi$ from the integrals:%
\[%
\begin{array}
[c]{ll}%
\left(  \displaystyle\int_{\Omega}\left\vert u\left(  x,t_{2}\right)
\right\vert ^{2}dx\right)  ^{1+M} & \leq\text{exp}\left[  -\left(  1+M\right)
\underset{x\in\overline{\Omega}}{\text{min}}\Phi\left(  x,t_{2}\right)
+M\underset{x\in\overline{\Omega}}{\text{max}}\Phi\left(  x,t_{1}\right)
\right] \\
& \quad\times\left(  \displaystyle\int_{\Omega}\left\vert u\left(
x,t_{1}\right)  \right\vert ^{2}dx\right)  ^{M}\displaystyle\int_{\Omega
}\left\vert u\left(  x,t_{3}\right)  \right\vert ^{2}e^{\Phi\left(
x,t_{3}\right)  }dx
\end{array}
\]
and%
\[%
\begin{array}
[c]{ll}%
\displaystyle\int_{\Omega}\left\vert u\left(  x,t_{3}\right)  \right\vert
^{2}e^{\Phi\left(  x,t_{3}\right)  }dx & =\displaystyle\int_{\omega}\left\vert
u\left(  x,t_{3}\right)  \right\vert ^{2}e^{\Phi\left(  x,t_{3}\right)
}dx+\displaystyle\int_{\left.  \Omega\right\backslash \omega}\left\vert
u\left(  x,t_{3}\right)  \right\vert ^{2}e^{\Phi\left(  x,t_{3}\right)  }dx\\
& \leq\text{exp}\left[  \underset{x\in\overline{\omega}}{\text{max}}%
\Phi\left(  x,t_{3}\right)  \right]  \displaystyle\int_{\omega}\left\vert
u\left(  x,t_{3}\right)  \right\vert ^{2}dx\\
& \quad+\text{exp}\left[  \underset{x\in\overline{\left.  \Omega
\right\backslash \omega}}{\text{max}}\Phi\left(  x,t_{3}\right)  \right]
\displaystyle\int_{\Omega}\left\vert u\left(  x,t_{3}\right)  \right\vert
^{2}dx\text{ .}%
\end{array}
\]
Therefore, we obtain that
\[%
\begin{array}
[c]{ll}
& \quad\left(  \displaystyle\int_{\Omega}\left\vert u\left(  x,t_{2}\right)
\right\vert ^{2}dx\right)  ^{1+M}\\
& \leq\text{exp}\left[  -\left(  1+M\right)  \underset{x\in\overline{\Omega}%
}{\text{min}}\Phi\left(  x,t_{2}\right)  +M\underset{x\in\overline{\Omega}%
}{\text{max}}\Phi\left(  x,t_{1}\right)  +\underset{x\in\overline{\omega}%
}{\text{max}}\Phi\left(  x,t_{3}\right)  \right] \\
& \quad\times\left(  \displaystyle\int_{\Omega}\left\vert u\left(
x,t_{1}\right)  \right\vert ^{2}dx\right)  ^{M}\displaystyle\int_{\omega
}\left\vert u\left(  x,t_{3}\right)  \right\vert ^{2}dx\\
& \quad+\text{exp}\left[  -\left(  1+M\right)  \underset{x\in\overline{\Omega
}}{\text{min}}\Phi\left(  x,t_{2}\right)  +M\underset{x\in\overline{\Omega}%
}{\text{max}}\Phi\left(  x,t_{1}\right)  +\underset{x\in\overline{\left.
\Omega\right\backslash \omega}}{\text{max}}\Phi\left(  x,t_{3}\right)  \right]
\\
& \quad\times\left(  \displaystyle\int_{\Omega}\left\vert u\left(
x,t_{1}\right)  \right\vert ^{2}dx\right)  ^{M}\displaystyle\int_{\Omega
}\left\vert u\left(  x,t_{3}\right)  \right\vert ^{2}dx\text{ .}%
\end{array}
\]
Using the fact that $\left\Vert u\left(  \cdot,T\right)  \right\Vert
\leq\left\Vert u\left(  \cdot,t\right)  \right\Vert \leq\left\Vert u\left(
\cdot,0\right)  \right\Vert $ $\forall0<t<T$, the above inequality becomes%
\[%
\begin{array}
[c]{ll}%
\left(  \left\Vert u\left(  \cdot,T\right)  \right\Vert ^{2}\right)  ^{1+M} &
\leq\text{exp}\left[  -\left(  1+M\right)  \underset{x\in\overline{\Omega}%
}{\text{min}}\Phi\left(  x,t_{2}\right)  +M\underset{x\in\overline{\Omega}%
}{\text{max}}\Phi\left(  x,t_{1}\right)  +\underset{x\in\overline{\omega}%
}{\text{max}}\Phi\left(  x,t_{3}\right)  \right] \\
& \quad\times\left(  \left\Vert u\left(  \cdot,0\right)  \right\Vert
^{2}\right)  ^{M}\displaystyle\int_{\omega}\left\vert u\left(  x,t_{3}\right)
\right\vert ^{2}dx\\
& \quad+\text{exp}\left[  -\left(  1+M\right)  \underset{x\in\overline{\Omega
}}{\text{min}}\Phi\left(  x,t_{2}\right)  +M\underset{x\in\overline{\Omega}%
}{\text{max}}\Phi\left(  x,t_{1}\right)  +\underset{x\in\overline{\left.
\Omega\right\backslash \omega}}{\text{max}}\Phi\left(  x,t_{3}\right)  \right]
\\
& \quad\times\left(  \left\Vert u\left(  \cdot,0\right)  \right\Vert
^{2}\right)  ^{1+M}\text{ .}%
\end{array}
\]

\begin{step}
\label{step2.1.5} Special weight function.
\end{step}

Assume that $\Phi\left(  x,t\right)  =\displaystyle\frac{\varphi\left(
x\right)  }{T-t+\hbar}$, we get that%
\[%
\begin{array}
[c]{ll}%
\left\Vert u\left(  \cdot,T\right)  \right\Vert ^{1+M} & \leq\text{exp}%
\frac{1}{2}\left[  -\frac{1+M}{T-t_{2}+\hbar}\underset{x\in\overline{\Omega}%
}{\text{min}}\varphi\left(  x\right)  +\frac{M}{T-t_{1}+\hbar}\underset
{x\in\overline{\Omega}}{\text{max}}\varphi\left(  x\right)  +\frac{1}%
{T-t_{3}+\hbar}\underset{x\in\overline{\omega}}{\text{max}}\varphi\left(
x\right)  \right] \\
& \quad\times\left\Vert u\left(  \cdot,0\right)  \right\Vert ^{M}\left\Vert
u\left(  \cdot,t_{3}\right)  \right\Vert _{L^{2}\left(  \omega\right)  }\\
& \quad+\text{exp}\frac{1}{2}\left[  -\frac{1+M}{T-t_{2}+\hbar}\underset
{x\in\overline{\Omega}}{\text{min}}\varphi\left(  x\right)  +\frac{M}%
{T-t_{1}+\hbar}\underset{x\in\overline{\Omega}}{\text{max}}\varphi\left(
x\right)  +\frac{1}{T-t_{3}+\hbar}\underset{x\in\overline{\left.
\Omega\right\backslash \omega}}{\text{max}}\varphi\left(  x\right)  \right] \\
& \quad\times\left\Vert u\left(  \cdot,0\right)  \right\Vert ^{1+M}\text{ .}%
\end{array}
\]
Choose $t_{3}=T$, $t_{2}=T-\ell\hbar$, $t_{1}=T-2\ell\hbar$ with $0<2\ell
\hbar<T$ and $\ell>1$, and denote
\[
M_{\ell}=\frac{\text{ln}\left(  \ell+1\right)  }{\text{ln}\left(  \frac
{2\ell+1}{\ell+1}\right)  }\text{ .}%
\]
Therefore, we have%
\[%
\begin{array}
[c]{ll}%
\left\Vert u\left(  \cdot,T\right)  \right\Vert ^{1+M_{\ell}} & \leq
\text{exp}\frac{1}{2\hbar}\left[  -\frac{1+M_{\ell}}{1+\ell}\underset
{x\in\overline{\Omega}}{\text{min}}\varphi\left(  x\right)  +\frac{M_{\ell}%
}{1+2\ell}\underset{x\in\overline{\Omega}}{\text{max}}\varphi\left(  x\right)
+\underset{x\in\overline{\omega}}{\text{max}}\varphi\left(  x\right)  \right]
\\
& \quad\times\left\Vert u\left(  \cdot,0\right)  \right\Vert ^{M_{\ell}%
}\left\Vert u\left(  \cdot,T\right)  \right\Vert _{L^{2}\left(  \omega\right)
}\\
& \quad+\text{exp}\frac{1}{2\hbar}\left[  -\frac{1+M_{\ell}}{1+\ell}%
\underset{x\in\overline{\Omega}}{\text{min}}\varphi\left(  x\right)
+\frac{M_{\ell}}{1+2\ell}\underset{x\in\overline{\Omega}}{\text{max}}%
\varphi\left(  x\right)  +\underset{x\in\overline{\left.  \Omega
\right\backslash \omega}}{\text{max}}\varphi\left(  x\right)  \right] \\
& \quad\times\left\Vert u\left(  \cdot,0\right)  \right\Vert ^{1+M_{\ell}%
}\text{ .}%
\end{array}
\]

\begin{step}
\label{step2.1.6} Assumption on weight function.
\end{step}

We construct $\varphi\left(  x\right)  $ and choose $\ell>1$ sufficiently
large in order that
\[
\left[  -\frac{1+M_{\ell}}{1+\ell}\underset{x\in\overline{\Omega}}{\text{min}%
}\varphi\left(  x\right)  +\frac{M_{\ell}}{1+2\ell}\underset{x\in
\overline{\Omega}}{\text{max}}\varphi\left(  x\right)  +\underset
{x\in\overline{\left.  \Omega\right\backslash \omega}}{\text{max}}%
\varphi\left(  x\right)  \right]  <0\text{ .}%
\]
Consequently, there are $C_{1}>0$ and $C_{2}>0$ such that for any $\hbar>0$
with $0<2\ell\hbar<T$,
\[
\left\Vert u\left(  \cdot,T\right)  \right\Vert ^{1+M_{\ell}}\leq
e^{C_{1}\frac{1}{\hbar}}\left\Vert u\left(  \cdot,0\right)  \right\Vert
^{M_{\ell}}\left\Vert u\left(  \cdot,T\right)  \right\Vert _{L^{2}\left(
\omega\right)  }+e^{-C_{2}\frac{1}{\hbar}}\left\Vert u\left(  \cdot,0\right)
\right\Vert ^{1+M_{\ell}}\text{ .}%
\]
Notice that $\left\Vert u\left(  \cdot,T\right)  \right\Vert \leq\left\Vert
u\left(  \cdot,0\right)  \right\Vert $ and for any $2\ell\hbar\geq T$, $1\leq
e^{C_{2}\frac{2\ell}{T}}e^{-C_{2}\frac{1}{\hbar}}$. We deduce that for any
$\hbar>0$,%
\[
\left\Vert u\left(  \cdot,T\right)  \right\Vert ^{1+M_{\ell}}\leq
e^{C_{1}\frac{1}{\hbar}}\left\Vert u\left(  \cdot,0\right)  \right\Vert
^{M_{\ell}}\left\Vert u\left(  \cdot,T\right)  \right\Vert _{L^{2}\left(
\omega\right)  }+e^{C_{2}\frac{2\ell}{T}}e^{-C_{2}\frac{1}{\hbar}}\left\Vert
u\left(  \cdot,0\right)  \right\Vert ^{1+M_{\ell}}\text{ .}%
\]
Finally, we choose $\hbar>0$ such that
\[
e^{C_{2}\frac{2\ell}{T}}e^{-C_{2}\frac{1}{\hbar}}\left\Vert u\left(
\cdot,0\right)  \right\Vert ^{1+M_{\ell}}=\frac{1}{2}\left\Vert u\left(
\cdot,T\right)  \right\Vert ^{1+M_{\ell}}\text{ ,}%
\]
that is,
\[
e^{C_{2}\frac{1}{\hbar}}:=2e^{C_{2}\frac{2\ell}{T}}\left(  \frac{\left\Vert
u\left(  \cdot,0\right)  \right\Vert }{\left\Vert u\left(  \cdot,T\right)
\right\Vert }\right)  ^{1+M_{\ell}}%
\]
in order that%
\[
\left\Vert u\left(  \cdot,T\right)  \right\Vert ^{1+M_{\ell}}\leq2\left(
2e^{C_{2}\frac{2\ell}{T}}\left(  \frac{\left\Vert u\left(  \cdot,0\right)
\right\Vert }{\left\Vert u\left(  \cdot,T\right)  \right\Vert }\right)
^{1+M_{\ell}}\right)  ^{\frac{C_{1}}{C_{2}}}\left\Vert u\left(  \cdot
,0\right)  \right\Vert ^{M_{\ell}}\left\Vert u\left(  \cdot,T\right)
\right\Vert _{L^{2}\left(  \omega\right)  }\text{ ,}%
\]
that is,
\[
\left\Vert u\left(  \cdot,T\right)  \right\Vert \leq2^{1+\frac{C_{1}}{C_{2}}%
}e^{C_{1}\frac{2\ell}{T}}\left(  \frac{\left\Vert u\left(  \cdot,0\right)
\right\Vert }{\left\Vert u\left(  \cdot,T\right)  \right\Vert }\right)
^{M_{\ell}+\left(  1+M_{\ell}\right)  \frac{C_{1}}{C_{2}}}\left\Vert u\left(
\cdot,T\right)  \right\Vert _{L^{2}\left(  \omega\right)  }\text{ .}%
\]
This ends to the desired inequality.

\bigskip

\bigskip

\subsection{$C^{2}$ bounded domain}

\bigskip

For $C^{2}$ bounded domain $\Omega$, we will use localized weight functions
exploiting a covering argument and propagation of smallness.

\bigskip

Let $0<r<R$, $x_{0}\in\Omega$ and $\delta\in\left(  0,1\right]  $. Denote
$R_{0}:=\left(  1+2\delta\right)  R$ and $B_{x_{0},r}:=\left\{  x;\left\vert
x-x_{0}\right\vert <r\right\}  $. Assume that $B_{x_{0},r}\Subset\Omega$ and
$\Omega\cap B_{x_{0},R_{0}}$ is star-shaped with respect to $x_{0}$. Let
$\left\langle \cdot,\cdot\right\rangle _{0}$ denote the usual scalar product
in $L^{2}\left(  \Omega\cap B_{x_{0},R_{0}}\right)  $ and let $\left\Vert
\cdot\right\Vert _{0}$ be its corresponding norm.

\bigskip

It suffices to prove the following result to get the desired observation
inequality at one point in time for the heat equation with Dirichlet boundary
condition in a $C^{2}$ bounded domain $\Omega$ (see \cite[Lemma 4 and Lemma 5
at p.493]{PWZ}).

\bigskip

\begin{lemm}
\label{lemma2.1} There is $\omega_{0}$ a nonempty open subset of $B_{x_{0},r}$
and constants $c$, $K>0$ and $\beta\in\left(  0,1\right)  $ such that for any
$T>0$ and $u_{0}\in L^{2}\left(  \Omega\right)  $,%
\[
\left\Vert e^{T\Delta}u_{0}\right\Vert _{L^{2}\left(  \Omega\cap B_{x_{0}%
,R}\right)  }\leq\left(  ce^{\frac{K}{T}}\left\Vert e^{T\Delta}u_{0}%
\right\Vert _{L^{2}\left(  \omega_{0}\right)  }\right)  ^{\beta}\left\Vert
u_{0}\right\Vert ^{1-\beta}\text{ .}%
\]

\end{lemm}

\bigskip

The strategy to establish the above Lemma \ref{lemma2.1} is as follows. It
will be divided into seven steps.

\begin{step}
\label{step2.2.1}Localization, symmetric and antisymmetric parts.
\end{step}

Let $\chi\in C_{0}^{\infty}\left(  B_{x_{0},R_{0}}\right)  $, $0\leq\chi\leq
1$, $\chi=1$ on $\left\{  x;\left\vert x-x_{0}\right\vert \leq\left(
1+3\delta/2\right)  R\right\}  $. Introduce $z=\chi u$. It solves
\[
\partial_{t}z-\Delta z=g:=-2\nabla\chi\cdot\nabla u-\Delta\chi u\text{ ,}%
\]
and furthermore, $z_{\left\vert \partial\left(  \Omega\cap B_{x_{0},R_{0}%
}\right)  \right.  }=0$. Let $\Phi$ be a sufficiently smooth function of
$\left(  x,t\right)  \in\mathbb{R}^{n}\times\mathbb{R}_{t}$ depending on
$x_{0}$. Set
\[
f\left(  x,t\right)  =z\left(  x,t\right)  e^{\Phi\left(  x,t\right)
/2}\text{ .}%
\]
We look for the equation solved by $f$ by computing $e^{\Phi\left(
x,t\right)  /2}\left(  \partial_{t}-\Delta\right)  \left(  e^{-\Phi\left(
x,t\right)  /2}f\left(  x,t\right)  \right)  $. It gives
\[
\partial_{t}f-\Delta f-\frac{1}{2}f\left(  \partial_{t}\Phi+\frac{1}%
{2}\left\vert \nabla\Phi\right\vert ^{2}\right)  +\nabla\Phi\cdot\nabla
f+\frac{1}{2}\Delta\Phi f=e^{\Phi/2}g\quad\text{in}~\left(  \Omega\cap
B_{x_{0},R_{0}}\right)  \times\left(  0,T\right)  \text{ ,}%
\]
and furthermore, $f_{\left\vert \partial\left(  \Omega\cap B_{x_{0},R_{0}%
}\right)  \right.  }=0$. Introduce
\[
\left\{
\begin{array}
[c]{ll}%
\mathcal{A}f=-\nabla\Phi\cdot\nabla f-\frac{1}{2}\Delta\Phi f\text{ ,} & \\
\mathcal{S}f=\Delta f+\eta f\text{ where }\eta=\frac{1}{2}\left(  \partial
_{t}\Phi+\frac{1}{2}\left\vert \nabla\Phi\right\vert ^{2}\right)  \text{ .} &
\end{array}
\right.
\]
It holds
\[
\left\{
\begin{array}
[c]{ll}%
\left\langle \mathcal{A}f,v\right\rangle _{0}=-\left\langle \mathcal{A}%
v,f\right\rangle _{0}\text{ ,} & \\
\left\langle \mathcal{S}f,v\right\rangle _{0}=\left\langle \mathcal{S}%
v,f\right\rangle _{0}\text{ for any }v\in H_{0}^{1}\left(  \Omega\cap
B_{x_{0},R_{0}}\right)  \text{ .} &
\end{array}
\right.
\]
Furthermore, one has
\[
\partial_{t}f-\mathcal{S}f-\mathcal{A}f=e^{\Phi/2}g\text{ .}%
\]

\begin{step}
\label{step2.2.2} Energy estimates.
\end{step}

Multiplying by $f$ the above equation and integrating over $\Omega\cap
B_{x_{0},R_{0}}$, we find that%
\[
\frac{1}{2}\frac{d}{dt}\left\Vert f\right\Vert _{0}^{2}+\left\langle
-\mathcal{S}f,f\right\rangle _{0}=\left\langle f,e^{\Phi/2}g\right\rangle
_{0}\text{ .}%
\]
Introduce the frequency function $t\mapsto\mathbf{N}\left(  t\right)  $
defined by%
\[
\mathbf{N}=\frac{\left\langle -\mathcal{S}f,f\right\rangle _{0}}{\left\Vert
f\right\Vert _{0}^{2}}\text{ .}%
\]
Thus,
\[
\frac{1}{2}\frac{d}{dt}\left\Vert f\right\Vert _{0}^{2}+\mathbf{N}\left\Vert
f\right\Vert _{0}^{2}=\left\langle e^{\Phi/2}g,f\right\rangle _{0}\text{ .}%
\]
Now, we compute the derivative of $\mathbf{N}$ and claim that:%
\[
\frac{d}{dt}\mathbf{N}\leq\frac{1}{\left\Vert f\right\Vert _{0}^{2}%
}\left\langle -\left(  \mathcal{S}^{\prime}+\left[  \mathcal{S},\mathcal{A}%
\right]  \right)  f,f\right\rangle _{0}-\frac{1}{\left\Vert f\right\Vert
_{0}^{2}}\int_{\partial\left(  \Omega\cap B_{x_{0},R_{0}}\right)  }%
\partial_{\nu}f\mathcal{A}fd\sigma+\frac{1}{2\left\Vert f\right\Vert _{0}^{2}%
}\left\Vert e^{\Phi/2}g\right\Vert _{0}^{2}\text{ .}%
\]
Indeed,%
\[%
\begin{array}
[c]{ll}
& \displaystyle\frac{d}{dt}\mathbf{N}\\
& =\displaystyle\frac{1}{\left\Vert f\right\Vert _{0}^{4}}\left(
-\displaystyle\frac{d}{dt}\left\langle \mathcal{S}f,f\right\rangle
_{0}\left\Vert f\right\Vert _{0}^{2}+\left\langle \mathcal{S}f,f\right\rangle
_{0}\displaystyle\frac{d}{dt}\left\Vert f\right\Vert _{0}^{2}\right) \\
& =\displaystyle\frac{-1}{\left\Vert f\right\Vert _{0}^{2}}\left[
\left\langle \mathcal{S}^{\prime}f,f\right\rangle _{0}+2\left\langle
\mathcal{S}f,f^{\prime}\right\rangle _{0}\right]  +\displaystyle\frac
{2}{\left\Vert f\right\Vert _{0}^{4}}\left[  \left\langle \mathcal{S}%
f,f\right\rangle _{0}^{2}+\left\langle \mathcal{S}f,f\right\rangle
_{0}\left\langle f,e^{\Phi/2}g\right\rangle _{0}\right] \\
& =\displaystyle\frac{-1}{\left\Vert f\right\Vert _{0}^{2}}\left[
\left\langle \mathcal{S}^{\prime}f,f\right\rangle _{0}+2\left\langle
\mathcal{S}f,f^{\prime}\right\rangle _{0}\right]  +\displaystyle\frac
{2}{\left\Vert f\right\Vert _{0}^{4}}\left[  \left\vert \left\langle
\mathcal{S}f,f\right\rangle _{0}+\frac{1}{2}\left\langle f,e^{\Phi
/2}g\right\rangle _{0}\right\vert ^{2}-\left\vert \frac{1}{2}\left\langle
f,e^{\Phi/2}g\right\rangle _{0}\right\vert ^{2}\right] \\
& =\displaystyle\frac{-1}{\left\Vert f\right\Vert _{0}^{2}}\left[
\left\langle \mathcal{S}^{\prime}f,f\right\rangle _{0}+2\left\langle
\mathcal{S}f,\mathcal{A}f\right\rangle _{0}\right]  +\displaystyle\frac
{-2}{\left\Vert f\right\Vert _{0}^{2}}\left[  \left\Vert \mathcal{S}%
f\right\Vert _{0}^{2}+\left\langle \mathcal{S}f,e^{\Phi/2}g\right\rangle
_{0}\right] \\
& \quad+\displaystyle\frac{2}{\left\Vert f\right\Vert _{0}^{4}}\left[
\left\vert \left\langle \mathcal{S}f+\frac{1}{2}e^{\Phi/2}g,f\right\rangle
_{0}\right\vert ^{2}-\left\vert \frac{1}{2}\left\langle f,e^{\Phi
/2}g\right\rangle _{0}\right\vert ^{2}\right] \\
& \leq\displaystyle\frac{-1}{\left\Vert f\right\Vert _{0}^{2}}\left[
\left\langle \mathcal{S}^{\prime}f,f\right\rangle _{0}+2\left\langle
\mathcal{S}f,\mathcal{A}f\right\rangle _{0}\right]  +\displaystyle\frac
{-2}{\left\Vert f\right\Vert _{0}^{2}}\left[  \left\Vert \mathcal{S}%
f\right\Vert _{0}^{2}+\left\langle \mathcal{S}f,e^{\Phi/2}g\right\rangle
_{0}\right] \\
& \quad+\displaystyle\frac{2}{\left\Vert f\right\Vert _{0}^{4}}\left\Vert
\mathcal{S}f+\frac{1}{2}e^{\Phi/2}g\right\Vert _{0}^{2}\left\Vert f\right\Vert
_{0}^{2}\\
& =\displaystyle\frac{-1}{\left\Vert f\right\Vert _{0}^{2}}\left[
\left\langle \mathcal{S}^{\prime}f,f\right\rangle _{0}+2\left\langle
\mathcal{S}f,\mathcal{A}f\right\rangle _{0}\right]  +\displaystyle\frac
{2}{\left\Vert f\right\Vert _{0}^{2}}\left\Vert \frac{1}{2}e^{\Phi
/2}g\right\Vert _{0}^{2}\text{ .}%
\end{array}
\]
In the third line, we used $\frac{1}{2}\frac{d}{dt}\left\Vert f\right\Vert
_{0}^{2}-\left\langle \mathcal{S}f,f\right\rangle _{0}=\left\langle
f,e^{\Phi/2}g\right\rangle _{0}$; In the fifth line, multiplying the equation
of $f$ by $\mathcal{S}f$, and integrating over $\Omega\cap B_{x_{0},R_{0}}$,
give
\[%
\begin{array}
[c]{ll}%
\left\langle \mathcal{S}f,f^{\prime}\right\rangle _{0} & =\left\langle
\mathcal{S}f,\left(  \mathcal{S}f+\mathcal{A}f+e^{\Phi/2}g\right)
\right\rangle _{0}\\
& =\left\Vert \mathcal{S}f\right\Vert _{0}^{2}+\left\langle \mathcal{S}%
f,\mathcal{A}f\right\rangle _{0}+\left\langle \mathcal{S}f,e^{\Phi
/2}g\right\rangle _{0}\text{ ;}%
\end{array}
\]
In the sixth line, we used Cauchy-Schwarz inequality. Finally, recall that
\[
2\left\langle \mathcal{S}f,\mathcal{A}f\right\rangle _{0}:=\left\langle
\left[  \mathcal{S},\mathcal{A}\right]  f,f\right\rangle _{0}%
+\displaystyle\int_{\partial\left(  \Omega\cap B_{x_{0},R_{0}}\right)
}\partial_{\nu}f\mathcal{A}fd\sigma\text{ .}%
\]

\begin{step}
\label{step2.2.3} Assumption on Carleman commutator.
\end{step}

Assume that $\displaystyle\int_{\partial\left(  \Omega\cap B_{x_{0},R_{0}%
}\right)  }\partial_{\nu}f\mathcal{A}fd\sigma\geq0$ on $\left(  0,T\right)  $
by the star-shaped property of $\Omega\cap B_{x_{0},R_{0}}$, and suppose that
\[
\left\langle -\left(  \mathcal{S}^{\prime}+\left[  \mathcal{S},\mathcal{A}%
\right]  \right)  f,f\right\rangle _{0}\leq\frac{1}{\Upsilon}\left\langle
-\mathcal{S}f,f\right\rangle _{0}%
\]
on $\left(  0,T\right)  $ where $\Upsilon\left(  t\right)  =T-t+\hbar$ and
$\hbar>0$. Therefore, the following differential inequalities hold.
\[
\left\{
\begin{array}
[c]{ll}%
\displaystyle\left\vert \frac{1}{2}\frac{d}{dt}\left\Vert f\left(
\cdot,t\right)  \right\Vert _{0}^{2}+\mathbf{N}\left(  t\right)  \left\Vert
f\left(  \cdot,t\right)  \right\Vert _{0}^{2}\right\vert \leq\left\Vert
e^{\Phi/2}g\left(  \cdot,t\right)  \right\Vert _{0}\left\Vert f\left(
\cdot,t\right)  \right\Vert _{0}\text{ ,} & \\
\displaystyle\frac{d}{dt}\mathbf{N}\left(  t\right)  \leq\displaystyle\frac
{1}{\Upsilon\left(  t\right)  }\mathbf{N}\left(  t\right)  +\displaystyle\frac
{\left\Vert e^{\Phi/2}g\left(  \cdot,t\right)  \right\Vert _{0}^{2}%
}{\left\Vert f\left(  \cdot,t\right)  \right\Vert _{0}^{2}}\text{ .} &
\end{array}
\right.
\]
By solving such system of differential inequalities, we have: For any
$0<t_{1}<t_{2}<t_{3}\leq T$,%
\[
\left(  \left\Vert f\left(  \cdot,t_{2}\right)  \right\Vert _{0}^{2}\right)
^{1+M}\leq\left(  \left\Vert f\left(  \cdot,t_{1}\right)  \right\Vert _{0}%
^{2}\right)  ^{M}\left\Vert f\left(  \cdot,t_{3}\right)  \right\Vert _{0}%
^{2}e^{2D}%
\]
where%
\[
M=\frac{\displaystyle\int_{t_{2}}^{t_{3}}\frac{1}{T-t+\hbar}dt}%
{\displaystyle\int_{t_{1}}^{t_{2}}\frac{1}{T-t+\hbar}dt}=\frac{-\text{ln}%
\left(  T-t_{3}+\hbar\right)  +\text{ln}\left(  T-t_{2}+\hbar\right)
}{-\text{ln}\left(  T-t_{2}+\hbar\right)  +\text{ln}\left(  T-t_{1}%
+\hbar\right)  }\text{ }%
\]
and
\[%
\begin{array}
[c]{ll}%
D & =M\left(  \left(  t_{2}-t_{1}\right)  \displaystyle\int_{t_{1}}^{t_{2}%
}\frac{\left\Vert e^{\Phi/2}g\left(  \cdot,t\right)  \right\Vert _{0}^{2}%
}{\left\Vert f\left(  \cdot,t\right)  \right\Vert _{0}^{2}}%
dt+\displaystyle\int_{t_{1}}^{t_{2}}\frac{\left\Vert e^{\Phi/2}g\left(
\cdot,t\right)  \right\Vert _{0}}{\left\Vert f\left(  \cdot,t\right)
\right\Vert _{0}}dt\right) \\
& \quad+\displaystyle\int_{t_{2}}^{t_{3}}\frac{1}{T-t+\hbar}%
dt\displaystyle\int_{t_{2}}^{t_{3}}\frac{\left\Vert e^{\Phi/2}g\left(
\cdot,t\right)  \right\Vert _{0}^{2}}{\left\Vert f\left(  \cdot,t\right)
\right\Vert _{0}^{2}}dt+\displaystyle\int_{t_{2}}^{t_{3}}\frac{\left\Vert
e^{\Phi/2}g\left(  \cdot,t\right)  \right\Vert _{0}}{\left\Vert f\left(
\cdot,t\right)  \right\Vert _{0}}dt\text{ .}%
\end{array}
\]
Indeed, we shall distinguish two cases: $t\in\left[  t_{1},t_{2}\right]  $;
$t\in\left[  t_{2},t_{3}\right]  $. For $t_{1}\leq t\leq t_{2}$, we integrate
$\left(  \left(  T-t+\hbar\right)  \mathbf{N}\left(  t\right)  \right)
^{\prime}\leq\left(  T-t+\hbar\right)  \frac{\left\Vert e^{\Phi/2}g\left(
\cdot,t\right)  \right\Vert _{0}^{2}}{\left\Vert f\left(  \cdot,t\right)
\right\Vert _{0}^{2}}$ over $\left(  t,t_{2}\right)  $ to get
\[
\left(  \frac{T-t_{2}+\hbar}{T-t+\hbar}\right)  \mathbf{N}\left(
t_{2}\right)  -\int_{t_{1}}^{t_{2}}\frac{\left\Vert e^{\Phi/2}g\left(
\cdot,s\right)  \right\Vert _{0}^{2}}{\left\Vert f\left(  \cdot,s\right)
\right\Vert _{0}^{2}}ds\leq\mathbf{N}\left(  t\right)  \text{ .}%
\]
Then we solve
\[
\frac{1}{2}\frac{d}{dt}\left\Vert f\right\Vert _{0}^{2}+\left[  \left(
\frac{T-t_{2}+\hbar}{T-t+\hbar}\right)  \mathbf{N}\left(  t_{2}\right)
-\int_{t_{1}}^{t_{2}}\frac{\left\Vert e^{\Phi/2}g\left(  \cdot,s\right)
\right\Vert _{0}^{2}}{\left\Vert f\left(  \cdot,s\right)  \right\Vert _{0}%
^{2}}ds-\frac{\left\Vert e^{\Phi/2}g\right\Vert _{0}}{\left\Vert f\right\Vert
_{0}}\right]  \left\Vert f\right\Vert _{0}^{2}\leq0
\]
and integrate it over $\left(  t_{1},t_{2}\right)  $ to obtain
\[%
\begin{array}
[c]{ll}%
e^{\displaystyle2\mathbf{N}\left(  t_{2}\right)  \int_{t_{1}}^{t_{2}}%
\frac{T-t_{2}+\hbar}{T-t+\hbar}dt} & \leq\displaystyle\frac{\left\Vert
f\left(  \cdot,t_{1}\right)  \right\Vert _{0}^{2}}{\left\Vert f\left(
\cdot,t_{2}\right)  \right\Vert _{0}^{2}}\\
& \quad\times e^{2\left(  t_{2}-t_{1}\right)  \left(  \displaystyle\int
_{t_{1}}^{t_{2}}\frac{\left\Vert e^{\Phi/2}g\left(  \cdot,t\right)
\right\Vert _{0}^{2}}{\left\Vert f\left(  \cdot,t\right)  \right\Vert _{0}%
^{2}}dt\right)  +\displaystyle2\int_{t_{1}}^{t_{2}}\frac{\left\Vert e^{\Phi
/2}g\left(  \cdot,t\right)  \right\Vert _{0}}{\left\Vert f\left(
\cdot,t\right)  \right\Vert _{0}}dt}\text{ .}%
\end{array}
\]
For $t_{2}\leq t\leq t_{3}$, we integrate $\left(  \left(  T-t+\hbar\right)
\mathbf{N}\left(  t\right)  \right)  ^{\prime}\leq\left(  T-t+\hbar\right)
\frac{\left\Vert e^{\Phi/2}g\left(  \cdot,t\right)  \right\Vert _{0}^{2}%
}{\left\Vert f\left(  \cdot,t\right)  \right\Vert _{0}^{2}}$ over $\left(
t_{2},t\right)  $ to get
\[
\mathbf{N}\left(  t\right)  \leq\frac{T-t_{2}+\hbar}{T-t+\hbar}\left(
\mathbf{N}\left(  t_{2}\right)  +\int_{t_{2}}^{t_{3}}\frac{\left\Vert
e^{\Phi/2}g\left(  \cdot,s\right)  \right\Vert _{0}^{2}}{\left\Vert f\left(
\cdot,s\right)  \right\Vert _{0}^{2}}ds\right)  \text{ .}%
\]
Then we solve
\[
0\leq\frac{1}{2}\frac{d}{dt}\left\Vert f\right\Vert _{0}^{2}+\left[
\frac{T-t_{2}+\hbar}{T-t+\hbar}\left(  \mathbf{N}\left(  t_{2}\right)
+\int_{t_{2}}^{t_{3}}\frac{\left\Vert e^{\Phi/2}g\left(  \cdot,s\right)
\right\Vert _{0}^{2}}{\left\Vert f\left(  \cdot,s\right)  \right\Vert _{0}%
^{2}}ds\right)  +\frac{\left\Vert e^{\Phi/2}g\right\Vert _{0}}{\left\Vert
f\right\Vert _{0}}\right]  \left\Vert f\right\Vert _{0}^{2}%
\]
and integrate it over $\left(  t_{2},t_{3}\right)  $ to obtain
\[%
\begin{array}
[c]{ll}%
\left\Vert f\left(  \cdot,t_{2}\right)  \right\Vert _{0}^{2} & \leq\left\Vert
f\left(  \cdot,t_{3}\right)  \right\Vert _{0}^{2}e^{\displaystyle2\mathbf{N}%
\left(  t_{2}\right)  \int_{t_{2}}^{t_{3}}\frac{T-t_{2}+\hbar}{T-t+\hbar}dt}\\
& \quad\times e^{\displaystyle2\int_{t_{2}}^{t_{3}}\frac{T-t_{2}+\hbar
}{T-t+\hbar}dt\int_{t_{2}}^{t_{3}}\frac{\left\Vert e^{\Phi/2}g\left(
\cdot,t\right)  \right\Vert _{0}^{2}}{\left\Vert f\left(  \cdot,t\right)
\right\Vert _{0}^{2}}dt+\displaystyle2\int_{t_{2}}^{t_{3}}\frac{\left\Vert
e^{\Phi/2}g\left(  \cdot,t\right)  \right\Vert _{0}}{\left\Vert f\left(
\cdot,t\right)  \right\Vert _{0}}dt}\text{ .}%
\end{array}
\]
Finally, combining the case $t_{1}\leq t\leq t_{2}$ and the case $t_{2}\leq
t\leq t_{3}$, we have%
\[%
\begin{array}
[c]{ll}%
\left\Vert f\left(  \cdot,t_{2}\right)  \right\Vert _{0}^{2} & \leq\left\Vert
f\left(  \cdot,t_{3}\right)  \right\Vert _{0}^{2}\left(  \displaystyle\frac
{\left\Vert f\left(  \cdot,t_{1}\right)  \right\Vert _{0}^{2}}{\left\Vert
f\left(  \cdot,t_{2}\right)  \right\Vert _{0}^{2}}\right)  ^{M}\\
& \quad\times e^{\displaystyle2M\left(  t_{2}-t_{1}\right)  \int_{t_{1}%
}^{t_{2}}\frac{\left\Vert e^{\Phi/2}g\left(  \cdot,t\right)  \right\Vert
_{0}^{2}}{\left\Vert f\left(  \cdot,t\right)  \right\Vert _{0}^{2}}%
dt}e^{\displaystyle2M\int_{t_{1}}^{t_{2}}\frac{\left\Vert e^{\Phi/2}g\left(
\cdot,t\right)  \right\Vert _{0}}{\left\Vert f\left(  \cdot,t\right)
\right\Vert _{0}}dt}\\
& \quad\times e^{\displaystyle2\int_{t_{2}}^{t_{3}}\frac{T-t_{2}+\hbar
}{T-t+\hbar}dt\int_{t_{2}}^{t_{3}}\frac{\left\Vert e^{\Phi/2}g\left(
\cdot,t\right)  \right\Vert _{0}^{2}}{\left\Vert f\left(  \cdot,t\right)
\right\Vert _{0}^{2}}dt}e^{\displaystyle2\int_{t_{2}}^{t_{3}}\frac{\left\Vert
e^{\Phi/2}g\left(  \cdot,t\right)  \right\Vert _{0}}{\left\Vert f\left(
\cdot,t\right)  \right\Vert _{0}}dt}%
\end{array}
\]
which implies the desired inequality.

\begin{step}
\label{step2.2.4} The rest term.
\end{step}

We estimate $\displaystyle\frac{\left\Vert e^{\Phi/2}g\right\Vert _{0}^{2}%
}{\left\Vert f\right\Vert _{0}^{2}}$. We begin by giving the following result.
(Recall that we have introduced $0<r<R$, $x_{0}\in\Omega$ and $\delta
\in\left(  0,1\right]  $).

\bigskip

\begin{lemm}
\label{lemma2.2} For any $T-\theta\leq t\leq T$, one has
\[
\frac{\left\Vert u\left(  \cdot,0\right)  \right\Vert ^{2}}{\left\Vert
u\left(  \cdot,t\right)  \right\Vert _{L^{2}\left(  \Omega\cap B_{x_{0}%
,\left(  1+\delta\right)  R}\right)  }^{2}}\leq e^{\left(  1+\delta\right)
\delta\frac{R^{2}}{2\theta}}\text{ }%
\]
where
\[
\frac{1}{\theta}=\frac{2}{\left(  \delta R\right)  ^{2}}%
\normalfont{\text{ln}}\left(  2e^{R^{2}\left(  1+\frac{1}{T}\right)  }%
\frac{\left\Vert u\left(  \cdot,0\right)  \right\Vert ^{2}}{\left\Vert
u\left(  \cdot,T\right)  \right\Vert _{L^{2}\left(  \Omega\cap B_{x_{0}%
,R}\right)  }^{2}}\right)  \text{ ,}%
\]
\textit{\ }with $0<\theta\leq\normalfont{\text{min}}\left(  1,T/2\right)  $ .
\end{lemm}

\bigskip

Indeed, denote $u\left(  x,t\right)  =e^{t\Delta}u_{0}\left(  x\right)  $ with
$u_{0}\in L^{2}(\Omega)$ non-null initial data. Recall that for any locally
Lipschitz function $\xi\left(  x,t\right)  $ such that $\partial_{t}\xi
+\frac{1}{2}\left\vert \nabla\xi\right\vert ^{2}\leq0$, the following integral
$\displaystyle\int_{\Omega}\left\vert u\left(  x,t\right)  \right\vert
^{2}e^{\xi\left(  x,t\right)  }dx$ is a decreasing function in $t$ by integral
maximun principle (see \cite{Gr}). Choose $\xi\left(  x,t\right)
=-\frac{\left\vert x-x_{0}\right\vert ^{2}}{2\left(  T-t+\epsilon\right)  }$,
then
\[
\int_{\Omega}\left\vert u\left(  x,T\right)  \right\vert ^{2}e^{-\frac
{\left\vert x-x_{0}\right\vert ^{2}}{2\epsilon}}dx\leq\int_{\Omega}\left\vert
u\left(  x,t\right)  \right\vert ^{2}e^{-\frac{\left\vert x-x_{0}\right\vert
^{2}}{2\left(  T-t+\epsilon\right)  }}dx\text{ .}%
\]
It implies that%
\[%
\begin{array}
[c]{ll}%
\left\Vert u\left(  \cdot,T\right)  \right\Vert _{L^{2}\left(  \Omega\cap
B_{x_{0},R}\right)  }^{2} & \leq\displaystyle e^{\frac{R^{2}}{2\epsilon}}%
\int_{\Omega\cap B_{x_{0},R}}\left\vert u\left(  x,T\right)  \right\vert
^{2}e^{-\frac{\left\vert x-x_{0}\right\vert ^{2}}{2\epsilon}}dx\\
& \leq\displaystyle e^{\frac{R^{2}}{2\epsilon}}\int_{\Omega}\left\vert
u\left(  x,t\right)  \right\vert ^{2}e^{-\frac{\left\vert x-x_{0}\right\vert
^{2}}{2\left(  T-t+\epsilon\right)  }}dx\\
& \leq\displaystyle e^{\frac{R^{2}}{2\epsilon}}\left\Vert u\left(
\cdot,t\right)  \right\Vert _{L^{2}\left(  \Omega\cap B_{x_{0},\left(
1+\delta\right)  R}\right)  }^{2}+\displaystyle e^{\frac{R^{2}}{2\epsilon}%
}e^{-\frac{R^{2}\left(  1+\delta\right)  ^{2}}{2\left(  T-t+\epsilon\right)
}}\left\Vert u\left(  \cdot,0\right)  \right\Vert ^{2}\text{ .}%
\end{array}
\]
Choose $T/2\leq T-\epsilon\delta\leq t\leq T$ with $0<\epsilon\leq T/2$ and
$\delta\in\left(  0,1\right]  $, then we get that
\[
\left\Vert u\left(  \cdot,T\right)  \right\Vert _{L^{2}\left(  \Omega\cap
B_{x_{0},R}\right)  }^{2}\leq e^{\frac{R^{2}}{2\epsilon}}\left\Vert u\left(
\cdot,t\right)  \right\Vert _{L^{2}\left(  \Omega\cap B_{x_{0},\left(
1+\delta\right)  R}\right)  }^{2}+e^{-\frac{\delta R^{2}}{2\epsilon}%
}\left\Vert u\left(  \cdot,0\right)  \right\Vert ^{2}\text{ .}%
\]
Choose
\[
\epsilon=\frac{\delta R^{2}}{2\text{ln}\left(  2e^{R^{2}\left(  1+\frac{1}%
{T}\right)  }\frac{\left\Vert u\left(  \cdot,0\right)  \right\Vert ^{2}%
}{\left\Vert u\left(  \cdot,T\right)  \right\Vert _{L^{2}\left(  \Omega\cap
B_{x_{0},R}\right)  }^{2}}\right)  }\leq\text{min}\left(  1,T/2\right)  \text{
,}%
\]
that is,
\[
e^{-\frac{\delta R^{2}}{2\epsilon}}\left\Vert u\left(  \cdot,0\right)
\right\Vert ^{2}=\frac{1}{2}e^{-R^{2}\left(  1+\frac{1}{T}\right)  }\left\Vert
u\left(  \cdot,T\right)  \right\Vert _{L^{2}\left(  \Omega\cap B_{x_{0}%
,R}\right)  }^{2}%
\]
in order that
\[
\left(  1-\frac{1}{2}e^{-R^{2}\left(  1+\frac{1}{T}\right)  }\right)
\left\Vert u\left(  \cdot,T\right)  \right\Vert _{L^{2}\left(  \Omega\cap
B_{x_{0},R}\right)  }^{2}\leq e^{\frac{R^{2}}{2\epsilon}}\left\Vert u\left(
\cdot,t\right)  \right\Vert _{L^{2}\left(  \Omega\cap B_{x_{0},\left(
1+\delta\right)  R}\right)  }^{2}%
\]
and
\[
e^{-\frac{\delta R^{2}}{2\epsilon}}\left\Vert u\left(  \cdot,0\right)
\right\Vert ^{2}\leq\frac{1}{2}e^{-R^{2}\left(  1+\frac{1}{T}\right)  }\left(
1-\frac{1}{2}e^{-R^{2}\left(  1+\frac{1}{T}\right)  }\right)  ^{-1}%
e^{\frac{R^{2}}{2\epsilon}}\left\Vert u\left(  \cdot,t\right)  \right\Vert
_{L^{2}\left(  \Omega\cap B_{x_{0},\left(  1+\delta\right)  R}\right)  }%
^{2}\text{ .}%
\]
This above inequality implies
\[
\frac{\left\Vert u\left(  \cdot,0\right)  \right\Vert ^{2}}{\left\Vert
u\left(  \cdot,t\right)  \right\Vert _{L^{2}\left(  \Omega\cap B_{x_{0}%
,\left(  1+\delta\right)  R}\right)  }^{2}}\leq e^{\left(  1+\delta\right)
\frac{R^{2}}{2\epsilon}}=e^{\left(  1+\delta\right)  \delta\frac{R^{2}%
}{2\theta}}\text{ ,}%
\]
that is,
\[
\frac{\left\Vert u\left(  \cdot,0\right)  \right\Vert ^{2}}{\left\Vert
u\left(  \cdot,t\right)  \right\Vert _{L^{2}\left(  \Omega\cap B_{x_{0}%
,\left(  1+\delta\right)  R}\right)  }^{2}}\leq e^{\left(  1+\frac{1}{\delta
}\right)  R^{2}\left(  1+\frac{1}{T}\right)  }\left(  \frac{2\left\Vert
u\left(  \cdot,0\right)  \right\Vert ^{2}}{\left\Vert u\left(  \cdot,T\right)
\right\Vert _{L^{2}\left(  \Omega\cap B_{x_{0},R}\right)  }^{2}}\right)
^{1+\frac{1}{\delta}}%
\]
as long as $T/2\leq T-\theta\leq t\leq T$ with
\[
\frac{1}{\theta}=\frac{2}{\left(  \delta R\right)  ^{2}}\text{ln}\left(
2e^{R^{2}\left(  1+\frac{1}{T}\right)  }\frac{\left\Vert u\left(
\cdot,0\right)  \right\Vert ^{2}}{\left\Vert u\left(  \cdot,T\right)
\right\Vert _{L^{2}\left(  \Omega\cap B_{x_{0},R}\right)  }^{2}}\right)
\text{ .}%
\]
Notice that $\theta\leq$min$\left(  1,T/2\right)  $. This completes the proof
of Lemma \ref{lemma2.2}. The interested reader may wish here to compare this
lemma's proof with \cite[p.660]{BP} or \cite[p.216]{EFV}.

\bigskip

Now we can estimate $\displaystyle\frac{\left\Vert e^{\Phi/2}g\right\Vert
_{0}^{2}}{\left\Vert f\right\Vert _{0}^{2}}$, by regularizing effect, as
follows.%
\[%
\begin{array}
[c]{ll}
& \quad\displaystyle\frac{\left\Vert e^{\Phi/2}g\left(  \cdot,t\right)
\right\Vert _{0}^{2}}{\left\Vert f\left(  \cdot,t\right)  \right\Vert _{0}%
^{2}}\\
& =\frac{\displaystyle\int_{\Omega\cap B_{x_{0},\left(  1+2\delta\right)  R}%
}\left\vert -2\nabla\chi\cdot\nabla u\left(  x,t\right)  -\Delta\chi u\left(
x,t\right)  \right\vert ^{2}e^{\Phi\left(  x,t\right)  }dx}{\displaystyle\int
_{\Omega\cap B_{x_{0},\left(  1+2\delta\right)  R}}\left\vert \chi u\left(
x,t\right)  \right\vert ^{2}e^{\Phi\left(  x,t\right)  }dx}\\
& \leq\frac{\displaystyle\int_{\Omega\cap\left\{  \left(  1+3\delta/2\right)
R\leq\left\vert x-x_{0}\right\vert \leq R_{0}\right\}  }\left\vert
-2\nabla\chi\cdot\nabla u\left(  x,t\right)  -\Delta\chi u\left(  x,t\right)
\right\vert ^{2}e^{\Phi\left(  x,t\right)  }dx}{\displaystyle\int_{\Omega\cap
B_{x_{0},\left(  1+\delta\right)  R}}\left\vert u\left(  x,t\right)
\right\vert ^{2}e^{\Phi\left(  x,t\right)  }dx}\\
& \leq\text{exp}\left[  -\underset{\left\vert x-x_{0}\right\vert \leq\left(
1+\delta\right)  R}{\text{min}}\Phi\left(  x,t\right)  +\underset{\left(
1+3\delta/2\right)  R\leq\left\vert x-x_{0}\right\vert \leq R_{0}}{\text{max}%
}\Phi\left(  x,t\right)  \right]  \displaystyle\frac{C\left(  1+\frac{1}%
{t}\right)  \left\Vert u\left(  \cdot,0\right)  \right\Vert ^{2}}{\left\Vert
u\left(  \cdot,t\right)  \right\Vert _{L^{2}\left(  \Omega\cap B_{x_{0}%
,\left(  1+\delta\right)  R}\right)  }^{2}}\\
& \leq\text{exp}\left[  -\underset{\left\vert x-x_{0}\right\vert \leq\left(
1+\delta\right)  R}{\text{min}}\Phi\left(  x,t\right)  +\underset{\left(
1+3\delta/2\right)  R\leq\left\vert x-x_{0}\right\vert \leq R_{0}}{\text{max}%
}\Phi\left(  x,t\right)  \right]  C\left(  1+\frac{1}{t}\right)  e^{\left(
1+\delta\right)  \delta\frac{R^{2}}{2\theta}}%
\end{array}
\]
as long as $T/2\leq T-\theta\leq t\leq T$.

\begin{step}
\label{step2.2.5} First assumption on the weight function.
\end{step}

We choose a weight function $\Phi\left(  x,t\right)  =\displaystyle\frac
{\varphi\left(  x\right)  }{T-t+\hbar}$ such that
\[
\underset{\left(  1+3\delta/2\right)  R\leq\left\vert x-x_{0}\right\vert \leq
R_{0}}{\text{max}}\varphi\left(  x\right)  -\underset{\left\vert
x-x_{0}\right\vert \leq\left(  1+\delta\right)  R}{\text{min}}\varphi\left(
x\right)  <0
\]
in order that
\[%
\begin{array}
[c]{ll}
& \quad-\underset{\left\vert x-x_{0}\right\vert \leq\left(  1+\delta\right)
R}{\text{min}}\Phi\left(  x,t\right)  +\underset{\left(  1+3\delta/2\right)
R\leq\left\vert x-x_{0}\right\vert \leq R_{0}}{\text{max}}\Phi\left(
x,t\right)  +\displaystyle\left(  1+\delta\right)  \delta\frac{R^{2}}{2\theta
}\\
& =\displaystyle\frac{-1}{T-t+\hbar}\left\vert \underset{\left\vert
x-x_{0}\right\vert \leq\left(  1+\delta\right)  R}{\text{min}}\varphi\left(
x\right)  -\underset{\left(  1+3\delta/2\right)  R\leq\left\vert
x-x_{0}\right\vert \leq R_{0}}{\text{max}}\varphi\left(  x\right)  \right\vert
+\displaystyle\left(  1+\delta\right)  \delta\frac{R^{2}}{2\theta}\\
& \leq\displaystyle\left(  1+\delta\right)  \delta\frac{R^{2}}{2\theta
}-\displaystyle\frac{1}{\left(  1+2\ell\right)  \hbar}\left\vert
\underset{\left\vert x-x_{0}\right\vert \leq\left(  1+\delta\right)
R}{\text{min}}\varphi\left(  x\right)  -\underset{\left(  1+3\delta/2\right)
R\leq\left\vert x-x_{0}\right\vert \leq R_{0}}{\text{max}}\varphi\left(
x\right)  \right\vert \text{ when }T-2\ell\hbar\leq t\\
& <0
\end{array}
\]
by taking
\[
\hbar\leq\theta\frac{1}{\left(  1+2\ell\right)  \left(  1+\delta\right)
\delta R^{2}}\left\vert \underset{\left\vert x-x_{0}\right\vert \leq\left(
1+\delta\right)  R}{\text{min}}\varphi\left(  x\right)  -\underset{\left(
1+3\delta/2\right)  R\leq\left\vert x-x_{0}\right\vert \leq R_{0}}{\text{max}%
}\varphi\left(  x\right)  \right\vert :=\theta C_{\left(  \ell,\varphi\right)
}\text{ .}%
\]
Here $\ell>1$. Combining with the previous Step \ref{step2.2.4}, one conclude
that for any $\hbar\leq\theta$min$\left(  C_{\left(  \ell,\varphi\right)
},1/\left(  2\ell\right)  \right)  $ and any $T-2\ell\hbar\leq t$,
\[
\frac{\left\Vert e^{\Phi/2}g\left(  \cdot,t\right)  \right\Vert _{0}^{2}%
}{\left\Vert f\left(  \cdot,t\right)  \right\Vert _{0}^{2}}\leq C\left(
1+\frac{1}{t}\right)  \text{ .}%
\]
Next, we choose $t_{3}=T$, $t_{2}=T-\ell\hbar$, $t_{1}=T-2\ell\hbar$, with
$\hbar\leq\theta$min$\left(  C_{\left(  \ell,\varphi\right)  },1/\left(
2\ell\right)  \right)  $. Therefore, the inequality of Step \ref{step2.2.3}
\[
\left(  \left\Vert f\left(  \cdot,t_{2}\right)  \right\Vert _{0}^{2}\right)
^{1+M}\leq\left(  \left\Vert f\left(  \cdot,t_{1}\right)  \right\Vert _{0}%
^{2}\right)  ^{M}\left\Vert f\left(  \cdot,t_{3}\right)  \right\Vert _{0}%
^{2}e^{2D}%
\]
becomes%
\[
\left(  \left\Vert f\left(  \cdot,T-\ell\hbar\right)  \right\Vert _{0}%
^{2}\right)  ^{1+M_{\ell}}\leq e^{2C_{\ell}\frac{1}{T}}\left(  \left\Vert
f\left(  \cdot,T-2\ell\hbar\right)  \right\Vert _{0}^{2}\right)  ^{M_{\ell}%
}\left\Vert f\left(  \cdot,T\right)  \right\Vert _{0}^{2}%
\]
as long as $\hbar\leq\theta$min$\left(  C_{\left(  \ell,\varphi\right)
},1/\left(  2\ell\right)  \right)  $. Here $C_{\ell}>0$ is a constant
depending on $\ell$ and recall that
\[
M_{\ell}=\frac{\text{ln}\left(  \ell+1\right)  }{\text{ln}\left(  \frac
{2\ell+1}{\ell+1}\right)  }\text{ .}%
\]

\begin{step}
\label{step2.2.6}
\end{step}

Let $\omega_{0}$ be a nonempty open subset of $B_{x_{0},r}$. Now by taking off
the weight function $\Phi\left(  x,t\right)  =\displaystyle\frac
{\varphi\left(  x\right)  }{T-t+\hbar}$ from the integrals, we have that for
any $0<\hbar\leq\theta$min$\left(  C_{\left(  \ell,\varphi\right)  },1/\left(
2\ell\right)  \right)  $,
\[%
\begin{array}
[c]{ll}
& \quad\left\Vert u\left(  \cdot,T-\ell\hbar\right)  \right\Vert
_{L^{2}\left(  \Omega\cap B_{x_{0},\left(  1+\delta\right)  R}\right)
}^{1+M_{\ell}}\\
& \leq\text{exp}\frac{1}{2\hbar}\left[  -\frac{1+M_{\ell}}{1+\ell}%
\underset{x\in\overline{\Omega}\cap\overline{B}_{x_{0},\left(  1+\delta
\right)  R}}{\text{min}}\varphi\left(  x\right)  +\frac{M_{\ell}}{1+2\ell
}\underset{x\in\overline{\Omega}\cap\overline{B}_{x_{0},R_{0}}}{\text{max}%
}\varphi\left(  x\right)  +\underset{x\in\overline{\omega_{0}}}{\text{max}%
}\varphi\left(  x\right)  \right] \\
& \quad\times e^{C_{\ell}\frac{1}{T}}\left\Vert u\left(  \cdot,0\right)
\right\Vert ^{M_{\ell}}\left\Vert u\left(  \cdot,T\right)  \right\Vert
_{L^{2}\left(  \omega_{0}\right)  }\\
& \quad+\text{exp}\frac{1}{2\hbar}\left[  -\frac{1+M_{\ell}}{1+\ell}%
\underset{x\in\overline{\Omega}\cap\overline{B}_{x_{0},\left(  1+\delta
\right)  R}}{\text{min}}\varphi\left(  x\right)  +\frac{M_{\ell}}{1+2\ell
}\underset{x\in\overline{\Omega}\cap\overline{B}_{x_{0},R_{0}}}{\text{max}%
}\varphi\left(  x\right)  +\underset{x\in\left.  \left(  \overline{\Omega}%
\cap\overline{B}_{x_{0},R_{0}}\right)  \right\backslash \omega_{0}}%
{\text{max}}\varphi\left(  x\right)  \right] \\
& \quad\times e^{C_{\ell}\frac{1}{T}}\left\Vert u\left(  \cdot,0\right)
\right\Vert ^{1+M_{\ell}}\text{ .}%
\end{array}
\]
But, by Lemma \ref{lemma2.2}, observe that%
\[
\frac{\left\Vert u\left(  \cdot,0\right)  \right\Vert ^{2}}{\left\Vert
u\left(  \cdot,t\right)  \right\Vert _{L^{2}\left(  \Omega\cap B_{x_{0}%
,\left(  1+\delta\right)  R}\right)  }^{2}}\leq e^{\left(  1+\delta\right)
\delta\frac{R^{2}}{2\theta}}%
\]
which gives, with $C_{\left(  \delta,R\right)  }:=\left(  1+\delta\right)
\delta\frac{R^{2}}{4}$ and $2\ell\hbar\leq\theta$,
\[
\left\Vert u\left(  \cdot,0\right)  \right\Vert \leq e^{\frac{1}{\theta
}C_{\left(  \delta,R\right)  }}\left\Vert u\left(  \cdot,T-\ell\hbar\right)
\right\Vert _{L^{2}\left(  \Omega\cap B_{x_{0},\left(  1+\delta\right)
R}\right)  }\text{ .}%
\]
Since $\left\Vert u\left(  \cdot,T\right)  \right\Vert \leq\left\Vert u\left(
\cdot,0\right)  \right\Vert $, we can see that%
\[
e^{-\frac{1}{\theta}C_{\left(  \delta,R\right)  }}\left\Vert u\left(
\cdot,T\right)  \right\Vert \leq\left\Vert u\left(  \cdot,T-\ell\hbar\right)
\right\Vert _{L^{2}\left(  \Omega\cap B_{x_{0},\left(  1+\delta\right)
R}\right)  }%
\]
and conclude that
\[%
\begin{array}
[c]{ll}
& \quad\left(  e^{-\frac{1}{\theta}C_{\left(  \delta,R\right)  }}\left\Vert
u\left(  \cdot,T\right)  \right\Vert \right)  ^{1+M_{\ell}}\\
& \leq\text{exp}\frac{1}{2\hbar}\left[  -\frac{1+M_{\ell}}{1+\ell}%
\underset{x\in\overline{\Omega}\cap\overline{B}_{x_{0},\left(  1+\delta
\right)  R}}{\text{min}}\varphi\left(  x\right)  +\frac{M_{\ell}}{1+2\ell
}\underset{x\in\overline{\Omega}\cap\overline{B}_{x_{0},R_{0}}}{\text{max}%
}\varphi\left(  x\right)  +\underset{x\in\overline{\omega_{0}}}{\text{max}%
}\varphi\left(  x\right)  \right] \\
& \quad\times e^{C_{\ell}\frac{1}{T}}\left\Vert u\left(  \cdot,0\right)
\right\Vert ^{M_{\ell}}\left\Vert u\left(  \cdot,T\right)  \right\Vert
_{L^{2}\left(  \omega_{0}\right)  }\\
& \quad+\text{exp}\frac{1}{2\hbar}\left[  -\frac{1+M_{\ell}}{1+\ell}%
\underset{x\in\overline{\Omega}\cap\overline{B}_{x_{0},\left(  1+\delta
\right)  R}}{\text{min}}\varphi\left(  x\right)  +\frac{M_{\ell}}{1+2\ell
}\underset{x\in\overline{\Omega}\cap\overline{B}_{x_{0},R_{0}}}{\text{max}%
}\varphi\left(  x\right)  +\underset{x\in\left.  \left(  \overline{\Omega}%
\cap\overline{B}_{x_{0},R_{0}}\right)  \right\backslash \omega_{0}}%
{\text{max}}\varphi\left(  x\right)  \right] \\
& \quad\times e^{C_{\ell}\frac{1}{T}}\left\Vert u\left(  \cdot,0\right)
\right\Vert ^{1+M_{\ell}}\text{ .}%
\end{array}
\]

\begin{step}
\label{step2.2.7} Second assumption on the weight function.
\end{step}

We construct $\varphi\left(  x\right)  $ and choose $\ell>1$ sufficiently
large in order that
\[
-\frac{1+M_{\ell}}{1+\ell}\underset{x\in\overline{\Omega}\cap\overline
{B}_{x_{0},\left(  1+\delta\right)  R}}{\text{min}}\varphi\left(  x\right)
+\frac{M_{\ell}}{1+2\ell}\underset{x\in\overline{\Omega}\cap\overline
{B}_{x_{0},R_{0}}}{\text{max}}\varphi\left(  x\right)  +\underset{x\in\left.
\left(  \overline{\Omega}\cap\overline{B}_{x_{0},R_{0}}\right)
\right\backslash \omega_{0}}{\text{max}}\varphi\left(  x\right)  <0\text{ .}%
\]
Consequently, there are $C_{1}>0$ and $C_{2}>0$ such that for any $\hbar>0$
with $\hbar\leq\theta$min$\left(  C_{\left(  \ell,\varphi\right)  },1/\left(
2\ell\right)  \right)  :=\theta C_{3}$,
\[
\left(  e^{-\frac{1}{\theta}C_{\left(  \delta,R\right)  }}\left\Vert u\left(
\cdot,T\right)  \right\Vert \right)  ^{1+M_{\ell}}\leq e^{C_{1}\frac{1}{\hbar
}}\left\Vert u\left(  \cdot,0\right)  \right\Vert ^{M_{\ell}}\left\Vert
u\left(  \cdot,T\right)  \right\Vert _{L^{2}\left(  \omega_{0}\right)
}+e^{-C_{2}\frac{1}{\hbar}}\left\Vert u\left(  \cdot,0\right)  \right\Vert
^{1+M_{\ell}}\text{ .}%
\]
On the other hand, for any $\hbar\geq\theta C_{3}$, $1\leq e^{\frac{C_{2}%
}{C_{3}}\frac{1}{\theta}}e^{-C_{2}\frac{1}{\hbar}}$. Therefore for any
$\hbar>0$,%
\[
\left(  e^{-\frac{1}{\theta}C_{\left(  \delta,R\right)  }}\left\Vert u\left(
\cdot,T\right)  \right\Vert \right)  ^{1+M_{\ell}}\leq e^{C_{1}\frac{1}{\hbar
}}\left\Vert u\left(  \cdot,0\right)  \right\Vert ^{M_{\ell}}\left\Vert
u\left(  \cdot,T\right)  \right\Vert _{L^{2}\left(  \omega_{0}\right)
}+e^{\frac{C_{2}}{C_{3}}\frac{1}{\theta}}e^{-C_{2}\frac{1}{\hbar}}\left\Vert
u\left(  \cdot,0\right)  \right\Vert ^{1+M_{\ell}}\text{ .}%
\]
Finally, we choose $\hbar>0$ such that
\[
e^{\frac{C_{2}}{C_{3}}\frac{1}{\theta}}e^{-C_{2}\frac{1}{\hbar}}\left\Vert
u\left(  \cdot,0\right)  \right\Vert ^{1+M_{\ell}}=\frac{1}{2}\left(
e^{-\frac{1}{\theta}C_{\left(  \delta,R\right)  }}\left\Vert u\left(
\cdot,T\right)  \right\Vert \right)  ^{1+M_{\ell}}\text{ ,}%
\]
that is,
\[
e^{C_{2}\frac{1}{\hbar}}:=2e^{\frac{C_{2}}{C_{3}}\frac{1}{\theta}}\left(
\frac{\left\Vert u\left(  \cdot,0\right)  \right\Vert }{e^{-\frac{1}{\theta
}C_{\left(  \delta,R\right)  }}\left\Vert u\left(  \cdot,T\right)  \right\Vert
}\right)  ^{1+M_{\ell}}%
\]
in order that%
\[%
\begin{array}
[c]{ll}%
\left(  e^{-C\frac{1}{\theta}}\left\Vert u\left(  \cdot,T\right)  \right\Vert
\right)  ^{1+M_{\ell}} & \leq2\left(  2e^{\frac{C_{2}}{C_{3}}\frac{1}{\theta}%
}\left(  \frac{\left\Vert u\left(  \cdot,0\right)  \right\Vert }{e^{-\frac
{1}{\theta}C_{\left(  \delta,R\right)  }}\left\Vert u\left(  \cdot,T\right)
\right\Vert }\right)  ^{1+M_{\ell}}\right)  ^{\frac{C_{1}}{C_{2}}}\\
& \quad\times\left\Vert u\left(  \cdot,0\right)  \right\Vert ^{M_{\ell}%
}\left\Vert u\left(  \cdot,T\right)  \right\Vert _{L^{2}\left(  \omega
_{0}\right)  }\text{ ,}%
\end{array}
\]
that is,
\[
e^{-C\frac{1}{\theta}}\left\Vert u\left(  \cdot,T\right)  \right\Vert
\leq2^{1+\frac{C_{1}}{C_{2}}}e^{\frac{C_{2}}{C_{3}}\frac{1}{\theta}}\left(
\frac{\left\Vert u\left(  \cdot,0\right)  \right\Vert }{e^{-\frac{1}{\theta
}C_{\left(  \delta,R\right)  }}\left\Vert u\left(  \cdot,T\right)  \right\Vert
}\right)  ^{M_{\ell}+\left(  1+M_{\ell}\right)  \frac{C_{1}}{C_{2}}}\left\Vert
u\left(  \cdot,T\right)  \right\Vert _{L^{2}\left(  \omega_{0}\right)  }\text{
.}%
\]
As a consequence, we obtain that for some $c>0$,
\[
\left\Vert u\left(  \cdot,T\right)  \right\Vert ^{1+c}\leq ce^{c\frac
{1}{\theta}}\left\Vert u\left(  \cdot,0\right)  \right\Vert ^{c}\left\Vert
u\left(  \cdot,T\right)  \right\Vert _{L^{2}\left(  \omega_{0}\right)  }\text{
.}%
\]
But recall the definition of $\theta$ in Lemma \ref{lemma2.2} saying that
\[
\frac{1}{\theta}=\frac{2}{\left(  \delta R\right)  ^{2}}\text{ln}\left(
2e^{R^{2}\left(  1+\frac{1}{T}\right)  }\frac{\left\Vert u\left(
\cdot,0\right)  \right\Vert ^{2}}{\left\Vert u\left(  \cdot,T\right)
\right\Vert _{L^{2}\left(  \Omega\cap B_{x_{0},R}\right)  }^{2}}\right)
\text{ .}%
\]
Therefore,%
\[
\left\Vert u\left(  \cdot,T\right)  \right\Vert ^{1+c}\leq\left[
2e^{R^{2}\left(  1+\frac{1}{T}\right)  }\frac{\left\Vert u\left(
\cdot,0\right)  \right\Vert ^{2}}{\left\Vert u\left(  \cdot,T\right)
\right\Vert _{L^{2}\left(  \Omega\cap B_{x_{0},R}\right)  }^{2}}\right]
^{c}c\left\Vert u\left(  \cdot,0\right)  \right\Vert ^{c}\left\Vert u\left(
\cdot,T\right)  \right\Vert _{L^{2}\left(  \omega_{0}\right)  }\text{ .}%
\]
which gives for some $K>0$, the following inequality%
\[
\left\Vert u\left(  \cdot,T\right)  \right\Vert _{L^{2}\left(  \Omega\cap
B_{x_{0},R}\right)  }^{1+K}\leq e^{K\left(  1+\frac{1}{T}\right)  }\left\Vert
u\left(  \cdot,0\right)  \right\Vert ^{K}\left\Vert u\left(  \cdot,T\right)
\right\Vert _{L^{2}\left(  \omega_{0}\right)  }\text{ }%
\]
and yields to the desired conclusion of Lemma \ref{lemma2.1}.

\bigskip

\section{The weight function}

\bigskip

In the previous section, the observation estimate at one time was derived by
using appropriate assumptions on the weight function $\Phi$ and by solving a
system of differential inequalities. Now, our goal is to explore different
explicit choices of weight function $\Phi$.

\bigskip

The weight function $\Phi_{\hbar}$ used in a series of results for the
doubling property or frequency function for heat equations was based on the
backward heat kernel (we also refer to \cite{BP} for parabolic equations where
the Euclidian distance is replaced by the geodesic distance). Precisely,%
\[
e^{\Phi_{\hbar}\left(  x,t\right)  }=G_{\hbar}\left(  x,t\right)  =\frac
{1}{\left(  T-t+\hbar\right)  ^{n/2}}e^{\frac{-\left\vert x-x_{0}\right\vert
^{2}}{4\left(  T-t+\hbar\right)  }}%
\]
or simply
\[
\Phi_{\hbar}\left(  x,t\right)  =\frac{-\left\vert x-x_{0}\right\vert ^{2}%
}{4\left(  T-t+\hbar\right)  }-\frac{n}{2}\text{ln}\left(  T-t+\hbar\right)
\text{ .}%
\]
It leads to the following differential inequalities (see \cite[Lemma 2 at
p.487]{PWZ}):

Define for $z\in H^{1}\left(  0,T;L^{2}\left(  \Omega\cap B_{x_{0},R_{0}%
}\right)  \right)  \cap L^{2}\left(  0,T;H^{2}\cap H_{0}^{1}\left(  \Omega\cap
B_{x_{0},R_{0}}\right)  \right)  $ and $t\in\left(  0,T\right]  $,
\[
\mathbf{N}_{\hbar}\left(  t\right)  =\frac{\displaystyle\int_{\Omega\cap
B_{x_{0},R_{0}}}\left\vert \nabla z\left(  x,t\right)  \right\vert
^{2}G_{\hbar}\left(  x,t\right)  dx}{\displaystyle\int_{\Omega\cap
B_{x_{0},R_{0}}}\left\vert z\left(  x,t\right)  \right\vert ^{2}G_{\hbar
}\left(  x,t\right)  dx}\text{, whenever }\int_{\Omega\cap B_{x_{0},R_{0}}%
}\left\vert z\left(  x,t\right)  \right\vert ^{2}dx\neq0\text{ .}%
\]
The following two properties hold.

\begin{description}
\item[$i)$]
\[%
\begin{array}
[c]{ll}
& \quad\displaystyle\frac{1}{2}\frac{d}{dt}\displaystyle\int_{\Omega\cap
B_{x_{0},R_{0}}}\left\vert z\left(  x,t\right)  \right\vert ^{2}G_{\hbar
}\left(  x,t\right)  dx+\displaystyle\int_{\Omega\cap B_{R_{0}}}\left\vert
\nabla z\left(  x,t\right)  \right\vert ^{2}G_{\hbar}\left(  x,t\right)  dx\\
& =\displaystyle\int_{\Omega\cap B_{x_{0},R_{0}}}z\left(  x,t\right)  \left(
\partial_{t}-\Delta\right)  z\left(  x,t\right)  G_{\hbar}\left(  x,t\right)
dx\text{ .}%
\end{array}
\]

\item[$ii)$] When $\Omega\cap B_{x_{0},R_{0}}$ is star-shaped with respect to
$x_{0}$,
\[
\frac{d}{dt}\mathbf{N}_{\hbar}\left(  t\right)  \leq\frac{1}{T-t+\hbar
}\mathbf{N}_{\hbar}\left(  t\right)  +\frac{\displaystyle\int_{\Omega\cap
B_{x_{0},R_{0}}}\left\vert \left(  \partial_{t}-\Delta\right)  z\left(
x,t\right)  \right\vert ^{2}G_{\hbar}\left(  x,t\right)  dx}{\displaystyle\int
_{\Omega\cap B_{x_{0},R_{0}}}\left\vert z\left(  x,t\right)  \right\vert
^{2}G_{\hbar}\left(  x,t\right)  dx}\text{ .}%
\]

\end{description}

\bigskip

The differential inequalities obtained with the Carleman commutator are given
in Step \ref{step2.1.2} and Step \ref{step2.2.2} of the previous Section 2:

Define for $f\in H^{1}\left(  0,T;L^{2}\left(  \Omega\cap B_{x_{0},R_{0}%
}\right)  \right)  \cap L^{2}\left(  0,T;H^{2}\cap H_{0}^{1}\left(  \Omega\cap
B_{x_{0},R_{0}}\right)  \right)  $ and $t\in\left(  0,T\right]  $,%
\[
\left\{
\begin{array}
[c]{ll}%
\mathcal{A}f=-\nabla\Phi\cdot\nabla f-\frac{1}{2}\Delta\Phi f\text{ ,} & \\
\mathcal{S}f=\Delta f+\eta f\text{ where }\eta=\frac{1}{2}\left(  \partial
_{t}\Phi+\frac{1}{2}\left\vert \nabla\Phi\right\vert ^{2}\right)  \text{ ,} &
\end{array}
\right.
\]
and%
\[
\mathbf{N}=\frac{\left\langle -\mathcal{S}f,f\right\rangle _{0}}{\left\Vert
f\right\Vert _{0}^{2}}\text{ .}%
\]
The following two properties hold.

\begin{description}
\item[$i)$]
\[
\frac{1}{2}\frac{d}{dt}\left\Vert f\right\Vert _{0}^{2}+\mathbf{N}\left\Vert
f\right\Vert _{0}^{2}=\left\langle \partial_{t}f-\mathcal{S}f-\mathcal{A}%
f,f\right\rangle _{0}\text{ .}%
\]

\item[$ii)$]
\[%
\begin{array}
[c]{ll}%
\displaystyle\frac{d}{dt}\mathbf{N} & \leq\displaystyle\frac{1}{\left\Vert
f\right\Vert _{0}^{2}}\left\langle -\left(  \mathcal{S}^{\prime}+\left[
\mathcal{S},\mathcal{A}\right]  \right)  f,f\right\rangle _{0}%
-\displaystyle\frac{1}{\left\Vert f\right\Vert _{0}^{2}}\displaystyle\int
_{\partial\left(  \Omega\cap B_{x_{0},R_{0}}\right)  }\partial_{\nu
}f\mathcal{A}fd\sigma\\
& \quad+\displaystyle\frac{1}{\left\Vert f\right\Vert _{0}^{2}}\left\Vert
\partial_{t}f-\mathcal{S}f-\mathcal{A}f\right\Vert _{0}^{2}\text{ .}%
\end{array}
\]

\end{description}

\bigskip

We will assume that $\displaystyle\int_{\partial\left(  \Omega\cap
B_{x_{0},R_{0}}\right)  }\partial_{\nu}f\mathcal{A}fd\sigma\geq0$ by the
star-shaped property of $\Omega\cap B_{x_{0},R_{0}}$. Now we focus our
attention on the term $\left\langle -\left(  \mathcal{S}^{\prime}+\left[
\mathcal{S},\mathcal{A}\right]  \right)  f,f\right\rangle _{0}$. We decompose
our presentation into three parts.

\begin{parti}
\label{part3.1} Key formula.
\end{parti}

We claim that:%
\[%
\begin{array}
[c]{ll}%
\left\langle -\left(  \mathcal{S}^{\prime}+\left[  \mathcal{S},\mathcal{A}%
\right]  \right)  f,f\right\rangle _{0} & =-2\displaystyle\int_{\Omega\cap
B_{x_{0},R_{0}}}\nabla f\cdot\nabla^{2}\Phi\nabla fdx\\
& \quad+\displaystyle\frac{1}{2}\int_{\Omega\cap B_{x_{0},R_{0}}}\Delta
^{2}\Phi\left\vert f\right\vert ^{2}dx-\displaystyle\int_{\Omega\cap
B_{x_{0},R_{0}}}\left(  \partial_{t}\eta+\nabla\Phi\cdot\nabla\eta\right)
\left\vert f\right\vert ^{2}dx\text{ }%
\end{array}
\]
which is, with the computation of $\partial_{t}\eta+\nabla\Phi\cdot\nabla\eta
$,%
\[%
\begin{array}
[c]{ll}
& \quad\left\langle -\left(  \mathcal{S}^{\prime}+\left[  \mathcal{S}%
,\mathcal{A}\right]  \right)  f,f\right\rangle _{0}\\
& =-2\displaystyle\int_{\Omega\cap B_{x_{0},R_{0}}}\nabla f\cdot\nabla^{2}%
\Phi\nabla fdx\\
& \quad+\displaystyle\frac{1}{2}\int_{\Omega\cap B_{x_{0},R_{0}}}\left(
\Delta^{2}\Phi-\partial_{t}^{2}\Phi-2\nabla\Phi\cdot\nabla\partial_{t}%
\Phi-\nabla\Phi\cdot\nabla^{2}\Phi\nabla\Phi\right)  \left\vert f\right\vert
^{2}dx\text{ .}%
\end{array}
\]

Proof of the claim .- First, $\mathcal{S}^{\prime}f=\partial_{t}\eta f$. Next,
we compute $\left[  \mathcal{S},\mathcal{A}\right]  f:=\mathcal{SA}%
f-\mathcal{AS}f$. Precisely, with standard summation notations,
\[%
\begin{array}
[c]{ll}%
\mathcal{SA}f & =\Delta\left(  -\nabla\Phi\cdot\nabla f-\frac{1}{2}\Delta\Phi
f\right)  +\eta\left(  -\nabla\Phi\cdot\nabla f-\frac{1}{2}\Delta\Phi f\right)
\\
& =-\Delta\nabla\Phi\cdot\nabla f-2\partial_{i}\nabla\Phi\cdot\partial
_{i}\nabla f-\nabla\Phi\cdot\Delta\nabla f-\frac{1}{2}\Delta^{2}\Phi
f-\nabla\Delta\Phi\cdot\nabla f-\frac{1}{2}\Delta\Phi\Delta f\\
& \quad-\eta\nabla\Phi\cdot\nabla f-\frac{1}{2}\eta\Delta\Phi f\text{ ,}\\
\mathcal{AS}f & =-\nabla\Phi\cdot\nabla\left(  \Delta f+\eta f\right)
-\frac{1}{2}\Delta\Phi\left(  \Delta f+\eta f\right) \\
& =-\nabla\Phi\cdot\nabla\Delta f-\nabla\Phi\cdot\nabla\eta f-\eta\nabla
\Phi\cdot\nabla f-\frac{1}{2}\Delta\Phi\Delta f-\frac{1}{2}\Delta\Phi\eta
f\text{ .}%
\end{array}
\]
This implies that%
\[
\left[  \mathcal{S},\mathcal{A}\right]  f=-2\partial_{i}\nabla\Phi
\cdot\partial_{i}\nabla f+\nabla\Phi\cdot\nabla\eta f-\frac{1}{2}\Delta
^{2}\Phi f-2\Delta\nabla\Phi\cdot\nabla f\text{ .}%
\]
Therefore, we obtain that
\[
-\left(  \mathcal{S}^{\prime}+\left[  \mathcal{S},\mathcal{A}\right]  \right)
f=2\partial_{i}\nabla\Phi\cdot\partial_{i}\nabla f-\left(  \partial_{t}%
\eta+\nabla\Phi\cdot\nabla\eta\right)  f+\frac{1}{2}\Delta^{2}\Phi
f+2\Delta\nabla\Phi\cdot\nabla f\text{ .}%
\]
Furthermore, by one integration by parts we have
\[
\left\langle \partial_{i}\nabla\Phi\cdot\partial_{i}\nabla f,f\right\rangle
_{0}=\frac{1}{2}\int_{\Omega\cap B_{x_{0},R_{0}}}\Delta^{2}\Phi\left\vert
f\right\vert ^{2}dx-\int_{\Omega\cap B_{x_{0},R_{0}}}\nabla f\cdot\nabla
^{2}\Phi\nabla fdx
\]
and%
\[
\left\langle \Delta\nabla\Phi\cdot\nabla f,f\right\rangle _{0}=-\frac{1}%
{2}\int_{\Omega\cap B_{x_{0},R_{0}}}\Delta^{2}\Phi\left\vert f\right\vert
^{2}dx\text{ .}%
\]
Combining the above equalities yields the desired formula. Then the claim follows.

\bigskip

Example linked with the heat kernel .- If
\[
\Phi\left(  x,t\right)  =\frac{-\left\vert x-x_{0}\right\vert ^{2}}{4\left(
T-t+\hbar\right)  }-\frac{n}{2}\text{ln}\left(  T-t+\hbar\right)  \text{ ,}%
\]
then we have
\[
\left\langle -\left(  \mathcal{S}^{\prime}+\left[  \mathcal{S},\mathcal{A}%
\right]  \right)  f,f\right\rangle _{0}=\frac{1}{\Upsilon}\left\langle
-\mathcal{S}f,f\right\rangle _{0}\text{ ,}%
\]
and
\[
\int_{\partial\left(  \Omega\cap B_{x_{0},R_{0}}\right)  }\partial_{\nu
}f\mathcal{A}fd\sigma=\frac{1}{2\Upsilon}\int_{\partial\left(  \Omega\cap
B_{x_{0},R_{0}}\right)  }\left\vert \partial_{\nu}f\right\vert ^{2}\left(
x-x_{0}\right)  \cdot\overrightarrow{\nu}d\sigma\geq0\text{ }%
\]
by the star-shaped property of $\Omega\cap B_{x_{0},R_{0}}$. Here and from
now, $\Upsilon\left(  t\right)  :=T-t+\hbar$ and $\overrightarrow{\nu}$ is the
outward unit normal vector to $\partial\left(  \Omega\cap B_{x_{0},R_{0}%
}\right)  $.

\begin{parti}
\label{part3.2} A particular form of the weight function.
\end{parti}

Assume that $\Phi\left(  x,t\right)  =\displaystyle\frac{\varphi\left(
x\right)  }{T-t+\hbar}$. Then, we can see that%
\[%
\begin{array}
[c]{ll}
& \quad\left\langle -\left(  \mathcal{S}^{\prime}+\left[  \mathcal{S}%
,\mathcal{A}\right]  \right)  f,f\right\rangle _{0}-\displaystyle\frac
{1}{\Upsilon}\left\langle -\mathcal{S}f,f\right\rangle _{0}\\
& =-\displaystyle\frac{1}{\Upsilon}\int_{\Omega\cap B_{x_{0},R_{0}}}\nabla
f\cdot\left(  2\nabla^{2}\varphi+I_{d}\right)  \nabla fdx+\displaystyle\frac
{1}{2\Upsilon}\int_{\Omega\cap B_{x_{0},R_{0}}}\Delta^{2}\varphi\left\vert
f\right\vert ^{2}dx\\
& \quad-\displaystyle\frac{1}{2\Upsilon^{3}}\int_{\Omega\cap B_{x_{0},R_{0}}%
}\left(  \varphi+\left\vert \nabla\varphi\right\vert ^{2}+\frac{1}{2}%
\nabla\varphi\cdot\left(  2\nabla^{2}\varphi+I_{d}\right)  \nabla
\varphi\right)  \left\vert f\right\vert ^{2}dx\text{ .}%
\end{array}
\]
Indeed,
\[%
\begin{array}
[c]{ll}%
\left\langle -\left(  \mathcal{S}^{\prime}+\left[  \mathcal{S},\mathcal{A}%
\right]  \right)  f,f\right\rangle _{0} & =-\displaystyle\frac{2}{\Upsilon
}\int_{\Omega\cap B_{x_{0},R_{0}}}\nabla f\cdot\nabla^{2}\varphi\nabla
fdx+\displaystyle\frac{1}{2\Upsilon}\int_{\Omega\cap B_{x_{0},R_{0}}}%
\Delta^{2}\varphi\left\vert f\right\vert ^{2}dx\\
& \quad-\displaystyle\frac{1}{\Upsilon^{3}}\int_{\Omega\cap B_{x_{0},R_{0}}%
}\left(  \varphi+\left\vert \nabla\varphi\right\vert ^{2}+\frac{1}{2}%
\nabla\varphi\cdot\nabla^{2}\varphi\nabla\varphi\right)  \left\vert
f\right\vert ^{2}dx\text{ ,}%
\end{array}
\]
and%
\[
\frac{1}{\Upsilon}\left\langle -\mathcal{S}f,f\right\rangle _{0}=\frac
{1}{\Upsilon}\int_{\Omega\cap B_{x_{0},R_{0}}}\left\vert \nabla f\right\vert
^{2}dx-\frac{1}{\Upsilon^{3}}\int_{\Omega\cap B_{x_{0},R_{0}}}\left(  \frac
{1}{2}\varphi+\frac{1}{4}\left\vert \nabla\varphi\right\vert ^{2}\right)
\left\vert f\right\vert ^{2}dx\text{ .}%
\]

\bigskip

Example of a weight function for localization with balls .- If
\[
\Phi\left(  x,t\right)  =\frac{-\left\vert x-x_{0}\right\vert ^{2}}{4\left(
T-t+\hbar\right)  }\text{ that is, }\varphi\left(  x\right)  =-\frac{1}%
{4}\left\vert x-x_{0}\right\vert ^{2}\text{ ,}%
\]
then we have
\[
\left\langle -\left(  \mathcal{S}^{\prime}+\left[  \mathcal{S},\mathcal{A}%
\right]  \right)  f,f\right\rangle _{0}=\frac{1}{\Upsilon}\left\langle
-\mathcal{S}f,f\right\rangle _{0}\text{ ,}%
\]
and
\[
\int_{\partial\left(  \Omega\cap B_{x_{0},R_{0}}\right)  }\partial_{\nu
}f\mathcal{A}fd\sigma=\frac{1}{2\Upsilon}\int_{\partial\left(  \Omega\cap
B_{x_{0},R_{0}}\right)  }\left\vert \partial_{\nu}f\right\vert ^{2}\left(
x-x_{0}\right)  \cdot\overrightarrow{\nu}d\sigma\geq0\text{ }%
\]
by the star-shaped property of $\Omega\cap B_{x_{0},R_{0}}$. One conclude
that, with such weight function $\Phi$, the assumptions of Step
\ref{step2.2.3} of the previous Section 2 are satisfied and therefore
\[
\left\{
\begin{array}
[c]{ll}%
\displaystyle\left\vert \frac{1}{2}\frac{d}{dt}\left\Vert f\left(
\cdot,t\right)  \right\Vert _{0}^{2}+\mathbf{N}\left(  t\right)  \left\Vert
f\left(  \cdot,t\right)  \right\Vert _{0}^{2}\right\vert \leq\left\Vert
e^{\Phi/2}g\left(  \cdot,t\right)  \right\Vert _{0}\left\Vert f\left(
\cdot,t\right)  \right\Vert _{0}\text{ ,} & \\
\displaystyle\frac{d}{dt}\mathbf{N}\left(  t\right)  \leq\displaystyle\frac
{1}{T-t+\hbar}\mathbf{N}\left(  t\right)  +\displaystyle\frac{\left\Vert
e^{\Phi/2}g\left(  \cdot,t\right)  \right\Vert _{0}^{2}}{\left\Vert f\left(
\cdot,t\right)  \right\Vert _{0}^{2}}\text{ .} &
\end{array}
\right.
\]
Now we check the assumptions on $\varphi\left(  x\right)  =-\frac{1}%
{4}\left\vert x-x_{0}\right\vert ^{2}$ at Step \ref{step2.2.5} and Step
\ref{step2.2.7} of the previous Section 2. We observe that
\[%
\begin{array}
[c]{ll}
& \quad\underset{\left(  1+3\delta/2\right)  R\leq\left\vert x-x_{0}%
\right\vert \leq R_{0}}{\text{max}}\varphi\left(  x\right)  -\underset
{\left\vert x-x_{0}\right\vert \leq\left(  1+\delta\right)  R}{\text{min}%
}\varphi\left(  x\right) \\
& =-\frac{1}{4}\left(  1+3\delta/2\right)  ^{2}R^{2}+\frac{1}{4}\left(
1+\delta\right)  ^{2}R^{2}<0
\end{array}
\]
and
\[%
\begin{array}
[c]{ll}
& -\frac{1+M_{\ell}}{1+\ell}\underset{x\in\overline{\Omega}\cap\overline
{B}_{x_{0},\left(  1+\delta\right)  R}}{\text{min}}\varphi\left(  x\right)
+\frac{M_{\ell}}{1+2\ell}\underset{x\in\overline{\Omega}\cap\overline
{B}_{x_{0},R_{0}}}{\text{max}}\varphi\left(  x\right)  +\underset{x\in\left.
\left(  \overline{\Omega}\cap\overline{B}_{x_{0},R_{0}}\right)
\right\backslash \omega_{0}}{\text{max}}\varphi\left(  x\right) \\
& \leq\left(  1+\frac{\text{ln}\left(  \ell+1\right)  }{\text{ln}\left(
\frac{2\ell+1}{\ell+1}\right)  }\right)  \frac{1}{1+\ell}\frac{1}{4}\left(
1+\delta\right)  ^{2}R^{2}-\frac{1}{4}r^{2}<0
\end{array}
\]
by choosing $\omega_{0}=B_{x_{0},r}\Subset\Omega$ with $0<r<R$ and by taking
$\ell>1$ sufficiently large.

\begin{parti}
\label{part3.3} The weight function for localization with annulus.
\end{parti}

Assume that $\varphi\left(  x\right)  =-a\left\vert x-x_{0}\right\vert
^{2}+b\left\vert x-x_{0}\right\vert ^{s}-c$ for some $a,b,c>0$ and $1\leq
s<2$. We would like to check the assumptions of the previous Section 2 and
find the adequate parameters $a,b,c,s$. First, we observe that the formula in
the previous Part \ref{part3.2}
\[%
\begin{array}
[c]{ll}
& \quad\left\langle -\left(  \mathcal{S}^{\prime}+\left[  \mathcal{S}%
,\mathcal{A}\right]  \right)  f,f\right\rangle _{0}-\displaystyle\frac
{1}{\Upsilon}\left\langle -\mathcal{S}f,f\right\rangle _{0}\\
& =-\displaystyle\frac{1}{\Upsilon}\int_{\Omega\cap B_{x_{0},R_{0}}}\nabla
f\cdot\left(  2\nabla^{2}\varphi+I_{d}\right)  \nabla fdx+\displaystyle\frac
{1}{2\Upsilon}\int_{\Omega\cap B_{x_{0},R_{0}}}\Delta^{2}\varphi\left\vert
f\right\vert ^{2}dx\\
& \quad-\displaystyle\frac{1}{2\Upsilon^{3}}\int_{\Omega\cap B_{x_{0},R_{0}}%
}\left(  \varphi+\left\vert \nabla\varphi\right\vert ^{2}+\frac{1}{2}%
\nabla\varphi\cdot\left(  2\nabla^{2}\varphi+I_{d}\right)  \nabla
\varphi\right)  \left\vert f\right\vert ^{2}dx
\end{array}
\]
gives%
\[%
\begin{array}
[c]{ll}
& \quad\left\langle -\left(  \mathcal{S}^{\prime}+\left[  \mathcal{S}%
,\mathcal{A}\right]  \right)  f,f\right\rangle _{0}-\displaystyle\frac
{1}{\Upsilon}\left\langle -\mathcal{S}f,f\right\rangle _{0}\\
& =-\displaystyle\frac{1}{\Upsilon}\left(  -4a+1\right)  \int_{\Omega\cap
B_{x_{0},R_{0}}}\left\vert \nabla f\right\vert ^{2}dx\\
& \quad-\displaystyle\frac{2bs}{\Upsilon}\left[  \displaystyle\int_{\Omega\cap
B_{x_{0},R_{0}}}\left\vert x-x_{0}\right\vert ^{s-2}\left\vert \nabla
f\right\vert ^{2}dx\right] \\
& \quad\quad\left.  \displaystyle-\left(  2-s\right)  \displaystyle\int
_{\Omega\cap B_{x_{0},R_{0}}}\left\vert x-x_{0}\right\vert ^{s-4}\left\vert
\left(  x-x_{0}\right)  \cdot\nabla f\right\vert ^{2}dx\right] \\
& \quad-\displaystyle\frac{1}{2\Upsilon}bs\left(  2-s\right)  \left(
n+s-2\right)  \left(  n+s-4\right)  \displaystyle\int_{\Omega\cap
B_{x_{0},R_{0}}}\left\vert x-x_{0}\right\vert ^{s-4}\left\vert f\right\vert
^{2}dx\\
& \quad+\displaystyle\frac{1}{2\Upsilon^{3}}c\displaystyle\int_{\Omega\cap
B_{x_{0},R_{0}}}\left\vert f\right\vert ^{2}dx\\
& \quad+\displaystyle\frac{1}{2\Upsilon^{3}}a\left(  1-2a\right)  \left(
1-4a\right)  \displaystyle\int_{\Omega\cap B_{x_{0},R_{0}}}\left\vert
x-x_{0}\right\vert ^{2}\left\vert f\right\vert ^{2}dx\\
& \quad+\displaystyle\frac{1}{2\Upsilon^{3}}b\left(  -1+6as-4a^{2}%
s-4a^{2}s^{2}\right)  \displaystyle\int_{\Omega\cap B_{x_{0},R_{0}}}\left\vert
x-x_{0}\right\vert ^{s}\left\vert f\right\vert ^{2}dx\\
& \quad-\displaystyle\frac{1}{2\Upsilon^{3}}\left(  bs\right)  ^{2}\left(
\frac{3}{2}+2a-4as\right)  \displaystyle\int_{\Omega\cap B_{x_{0},R_{0}}%
}\left\vert x-x_{0}\right\vert ^{2s-2}\left\vert f\right\vert ^{2}dx\\
& \quad-\displaystyle\frac{1}{2\Upsilon^{3}}\left(  bs\right)  ^{3}\left(
s-1\right)  \displaystyle\int_{\Omega\cap B_{x_{0},R_{0}}}\left\vert
x-x_{0}\right\vert ^{3s-4}\left\vert f\right\vert ^{2}dx\text{ .}%
\end{array}
\]
We start to choose $a=\frac{1}{4}$. Next we treat the third line of the above
formula by using Cauchy-Schwarz inequality, we find that
\[%
\begin{array}
[c]{ll}
& \quad\left\langle -\left(  \mathcal{S}^{\prime}+\left[  \mathcal{S}%
,\mathcal{A}\right]  \right)  f,f\right\rangle _{0}-\displaystyle\frac
{1}{\Upsilon}\left\langle -\mathcal{S}f,f\right\rangle _{0}\\
& \leq-\displaystyle\frac{1}{2\Upsilon}bs\left(  2-s\right)  \left(
n+s-2\right)  \left(  n+s-4\right)  \displaystyle\int_{\Omega\cap
B_{x_{0},R_{0}}}\left\vert x-x_{0}\right\vert ^{s-4}\left\vert f\right\vert
^{2}dx\\
& \quad+\displaystyle\frac{1}{2\Upsilon^{3}}c\displaystyle\int_{\Omega\cap
B_{x_{0},R_{0}}}\left\vert f\right\vert ^{2}dx\\
& \quad+\displaystyle\frac{1}{2\Upsilon^{3}}b\left(  -1+\frac{5}{4}s-\frac
{1}{4}s^{2}\right)  \displaystyle\int_{\Omega\cap B_{x_{0},R_{0}}}\left\vert
x-x_{0}\right\vert ^{s}\left\vert f\right\vert ^{2}dx\\
& \quad-\displaystyle\frac{1}{2\Upsilon^{3}}\left(  bs\right)  ^{2}\left(
2-s\right)  \displaystyle\int_{\Omega\cap B_{x_{0},R_{0}}}\left\vert
x-x_{0}\right\vert ^{2s-2}\left\vert f\right\vert ^{2}dx\\
& \quad-\displaystyle\frac{1}{2\Upsilon^{3}}\left(  bs\right)  ^{3}\left(
s-1\right)  \displaystyle\int_{\Omega\cap B_{x_{0},R_{0}}}\left\vert
x-x_{0}\right\vert ^{3s-4}\left\vert f\right\vert ^{2}dx\text{ .}%
\end{array}
\]
In order that $\left\langle -\left(  \mathcal{S}^{\prime}+\left[
\mathcal{S},\mathcal{A}\right]  \right)  f,f\right\rangle _{0}%
-\displaystyle\frac{1}{\Upsilon}\left\langle -\mathcal{S}f,f\right\rangle
_{0}\leq0$, we can take $n\geq3$ and $s=1$ with $c\leq\left(  bs\right)
^{2}\left(  2-s\right)  =b^{2}$. Another choice is $n\geq3$ and $s=\frac{4}%
{3}$ which gives%
\[%
\begin{array}
[c]{ll}
& \quad\left\langle -\left(  \mathcal{S}^{\prime}+\left[  \mathcal{S}%
,\mathcal{A}\right]  \right)  f,f\right\rangle _{0}-\displaystyle\frac
{1}{\Upsilon}\left\langle -\mathcal{S}f,f\right\rangle _{0}\\
& \leq\quad\displaystyle\frac{1}{\Upsilon^{3}}\frac{c}{2}\displaystyle\int
_{\Omega\cap B_{x_{0},R_{0}}}\left\vert f\right\vert ^{2}dx\\
& \quad+\displaystyle\frac{1}{\Upsilon^{3}}\frac{b}{9}\displaystyle\int
_{\Omega\cap B_{x_{0},R_{0}}}\left\vert x-x_{0}\right\vert ^{4/3}\left\vert
f\right\vert ^{2}dx-\displaystyle\frac{1}{\Upsilon^{3}}\frac{1}{3}\left(
\frac{4}{3}b\right)  ^{2}\displaystyle\int_{\Omega\cap B_{x_{0},R_{0}}%
}\left\vert x-x_{0}\right\vert ^{2/3}\left\vert f\right\vert ^{2}dx\\
& \quad-\displaystyle\frac{1}{\Upsilon^{3}}\frac{1}{6}\left(  \frac{4}%
{3}b\right)  ^{3}\displaystyle\int_{\Omega\cap B_{x_{0},R_{0}}}\left\vert
f\right\vert ^{2}dx\text{ ,}%
\end{array}
\]
and finally, we can choose $\frac{c}{2}\leq\frac{1}{6}\left(  \frac{4}%
{3}b\right)  ^{3}$ and $\frac{b}{9}R_{0}^{2/3}\leq\frac{1}{3}\left(  \frac
{4}{3}b\right)  ^{2}$.

\bigskip

Now, set $\varphi\left(  x\right)  =-\frac{1}{4}\left\vert x-x_{0}\right\vert
^{2}+\frac{1}{4}\left\vert x-x_{0}\right\vert ^{4/3}-\left(  \frac{1}%
{3}\right)  ^{4}$ and $R=1$, $R_{0}=\left(  \frac{4}{3}\right)  ^{3/2}%
\approx1.53$ and $\delta=\frac{R_{0}-1}{2}\approx0.26$. Here $1<(1+3\delta
/2)R\approx1.40<R_{0}=(1+2\delta)R$. We can see that for $n\geq3$,
$\left\langle -\left(  \mathcal{S}^{\prime}+\left[  \mathcal{S},\mathcal{A}%
\right]  \right)  f,f\right\rangle _{0}-\displaystyle\frac{1}{\Upsilon
}\left\langle -\mathcal{S}f,f\right\rangle _{0}\leq0$.

\bigskip

We write $\varphi\left(  x\right)  =W\left(  \left\vert x-x_{0}\right\vert
\right)  $ with $W\left(  \rho\right)  =-\frac{1}{4}\rho^{2}+\frac{1}{4}%
\rho^{4/3}-\left(  \frac{1}{3}\right)  ^{4}$. We have $W\left(  0\right)
=W\left(  1\right)  =-\left(  \frac{1}{3}\right)  ^{4}$, $W^{\prime}\left(
\left(  2/3\right)  ^{3/2}\right)  =0$ and $\rho\mapsto W\left(  \rho\right)
$ is strictly decreasing for $\rho\geq1$.

\bigskip

Finally, we check the assumptions on $\varphi$ at Step \ref{step2.2.5} and
Step \ref{step2.2.7} of the previous Section 2: We observe that
\[%
\begin{array}
[c]{ll}
& \quad\underset{\left(  1+3\delta/2\right)  R\leq\left\vert x-x_{0}%
\right\vert \leq R_{0}}{\text{max}}\varphi\left(  x\right)  -\underset
{\left\vert x-x_{0}\right\vert \leq\left(  1+\delta\right)  R}{\text{min}%
}\varphi\left(  x\right) \\
& \leq W\left(  \left(  1+3\delta/2\right)  R\right)  -W\left(  \left(
1+\delta\right)  R\right)  <0
\end{array}
\]
(because $\rho\mapsto W\left(  \rho\right)  $ is strictly decreasing for
$\rho\geq1=R$) and
\[%
\begin{array}
[c]{ll}
& -\frac{1+M_{\ell}}{1+\ell}\underset{x\in\overline{\Omega}\cap\overline
{B}_{x_{0},\left(  1+\delta\right)  R}}{\text{min}}\varphi\left(  x\right)
+\frac{M_{\ell}}{1+2\ell}\underset{x\in\overline{\Omega}\cap\overline
{B}_{x_{0},R_{0}}}{\text{max}}\varphi\left(  x\right)  +\underset{x\in\left.
\left(  \overline{\Omega}\cap\overline{B}_{x_{0},R_{0}}\right)
\right\backslash \omega_{0}}{\text{max}}\varphi\left(  x\right) \\
& \leq-\left(  1+\frac{\text{ln}\left(  \ell+1\right)  }{\text{ln}\left(
\frac{2\ell+1}{\ell+1}\right)  }\right)  \frac{1}{1+\ell}W\left(  \left(
1+\delta\right)  R\right)  +\frac{\text{ln}\left(  \ell+1\right)  }%
{\text{ln}\left(  \frac{2\ell+1}{\ell+1}\right)  }\frac{1}{1+2\ell}W\left(
\left(  2/3\right)  ^{3/2}\right)  +W\left(  r_{0}\right)  <0
\end{array}
\]
by choosing $\omega_{0}=\left\{  x;r_{0}<\left\vert x-x_{0}\right\vert
<r\right\}  \Subset\Omega$ with $0<r_{0}<r<1$, $W\left(  r_{0}\right)
=W\left(  r\right)  \in\left(  -\left(  \frac{1}{3}\right)  ^{4},0\right)  $
and by taking $\ell>1$ sufficiently large.

\bigskip

\bigskip

\section{Proof of Theorem 1.1}

\bigskip

\bigskip

The observability estimate in Theorem \ref{theorem1.1} can be deduced from the
observation inequality at one time (see \cite{PW2} or the following Lemma
\ref{lemma4.1}). It was noticed in \cite{AEWZ} that the spectral inequality in
Theorem \ref{theorem1.1} is a consequence of the observation inequality at one
time (see Lemma A in Appendix (see page \pageref{appendix})).

\bigskip

\begin{lemm}
\label{lemma4.1} Let $\omega$ be a nonempty open subset of $\Omega$. Let
$p\in\left[  1,2\right]  $, $\gamma>0$, $\beta\in\left(  0,1\right)  $,
$C_{\beta}=\frac{1+\beta}{\beta\left[  \left(  1+\beta\right)  ^{1/\left(
2\gamma\right)  }-1\right]  ^{\gamma}}$ and $c$, $K>0$. Suppose that for any
$u_{0}\in L^{2}\left(  \Omega\right)  $ and any $T>0$,
\[
\left\Vert e^{T\Delta}u_{0}\right\Vert _{L^{2}\left(  \Omega\right)  }%
\leq\left(  ce^{\frac{1}{T^{\gamma}}K}\left\Vert e^{T\Delta}u_{0}\right\Vert
_{L^{p}\left(  \omega\right)  }\right)  ^{\beta}\left\Vert u_{0}\right\Vert
_{L^{2}\left(  \Omega\right)  }^{1-\beta}\text{ .}%
\]
Then for any $u_{0}\in L^{2}\left(  \Omega\right)  $ and any\thinspace$T>0$,
one has%
\[
\left\Vert e^{T\Delta}u_{0}\right\Vert _{L^{2}\left(  \Omega\right)  }%
\leq\dfrac{c}{K^{1/\gamma}}e^{\frac{1}{T^{\gamma}}\left(  1+1/\gamma\right)
KC_{\beta}}\int_{0}^{T}\left\Vert e^{t\Delta}u_{0}\right\Vert _{L^{p}\left(
\omega\right)  }dt\text{ .}%
\]

\end{lemm}

\bigskip

The above lemma is somehow standard, but we still give the proof here to make
a self-contained discussion.

\bigskip

Proof of Lemma \ref{lemma4.1} .- First, by Young inequality, the following
interpolation estimate
\[
\left\Vert u\left(  \cdot,T\right)  \right\Vert _{L^{2}\left(  \Omega\right)
}\leq\left(  ce^{\frac{K}{T^{\gamma}}}\left\Vert u\left(  \cdot,T\right)
\right\Vert _{L^{p}\left(  \omega\right)  }\right)  ^{\beta}\left\Vert
u\left(  \cdot,0\right)  \right\Vert _{L^{2}\left(  \Omega\right)  }^{1-\beta}%
\]
implies that for any $\varepsilon>0$, we have%
\[
\left\Vert u\left(  \cdot,T\right)  \right\Vert _{L^{2}\left(  \Omega\right)
}\leq\dfrac{1}{\varepsilon^{\frac{1-\beta}{\beta}}}ce^{\frac{K}{T^{\gamma}}%
}\left\Vert u\left(  \cdot,T\right)  \right\Vert _{L^{p}\left(  \omega\right)
}+\varepsilon\left\Vert u\left(  \cdot,0\right)  \right\Vert _{L^{2}\left(
\Omega\right)  }\text{ .}%
\]
Next introduce a decreasing sequence $\left(  T_{m}\right)  _{m\geq0}$ of
positive real numbers defined by
\[
T_{m}=\dfrac{T}{z^{m}}\text{ with }z>1\text{ .}%
\]
Take $0<T_{m+2}<T_{m+1}\leq t<T_{m}<\cdot\cdot\cdot<T$ and apply the
observation estimate at one time $t$ with initial time $T_{m+2}$. We find
that
\[
\left\Vert u\left(  \cdot,t\right)  \right\Vert _{L^{2}\left(  \Omega\right)
}\leq\dfrac{1}{\varepsilon^{\frac{1-\beta}{\beta}}}ce^{\frac{K}{\left(
t-T_{m+2}\right)  ^{\gamma}}}\left\Vert u\left(  \cdot,t\right)  \right\Vert
_{L^{p}\left(  \omega\right)  }+\varepsilon\left\Vert u\left(  \cdot
,T_{m+2}\right)  \right\Vert _{L^{2}\left(  \Omega\right)  }\text{ .}%
\]
Since $\left\Vert u\left(  \cdot,T_{m}\right)  \right\Vert \leq\left\Vert
u\left(  \cdot,t\right)  \right\Vert $, we deduce that
\[
\left\Vert u\left(  \cdot,T_{m}\right)  \right\Vert _{L^{2}\left(
\Omega\right)  }\leq\dfrac{1}{\varepsilon^{\frac{1-\beta}{\beta}}}ce^{\frac
{K}{\left(  t-T_{m+2}\right)  ^{\gamma}}}\left\Vert u\left(  \cdot,t\right)
\right\Vert _{L^{p}\left(  \omega\right)  }+\varepsilon\left\Vert u\left(
\cdot,T_{m+2}\right)  \right\Vert _{L^{2}\left(  \Omega\right)  }\text{ .}%
\]
Now, integrate the above inequality over $\left(  T_{m+1},T_{m}\right)  $, it
yields that
\[
\left\Vert u\left(  \cdot,T_{m}\right)  \right\Vert _{L^{2}\left(
\Omega\right)  }\leq\dfrac{1}{\varepsilon^{\frac{1-\beta}{\beta}}%
}{\displaystyle}\frac{c}{T_{m}-T_{m+1}}e^{\frac{K}{\left(  T_{m+1}%
-T_{m+2}\right)  ^{\gamma}}}{\displaystyle\int_{T_{m+1}}^{T_{m}}}\left\Vert
u\left(  \cdot,t\right)  \right\Vert _{L^{p}\left(  \omega\right)
}dt+\varepsilon\left\Vert u\left(  \cdot,T_{m+2}\right)  \right\Vert
_{L^{2}\left(  \Omega\right)  }%
\]
which implies, since ${\displaystyle}\frac{c}{T_{m}-T_{m+1}}e^{\frac
{K}{\left(  T_{m+1}-T_{m+2}\right)  ^{\gamma}}}={\displaystyle}\frac{c}%
{z}\frac{z^{m+2}}{\left(  z-1\right)  T}e^{\left(  K^{1/\gamma}\frac{z^{m+2}%
}{\left(  z-1\right)  T}\right)  ^{\gamma}}$,
\[
\left\Vert u\left(  \cdot,T_{m}\right)  \right\Vert _{L^{2}\left(
\Omega\right)  }\leq\dfrac{1}{\varepsilon^{\frac{1-\beta}{\beta}}%
}{\displaystyle}\frac{c}{zK^{1/\gamma}}e^{\frac{\left(  1+1/\gamma\right)
K}{\left(  z-1\right)  ^{\gamma}T^{\gamma}}z^{\gamma\left(  m+2\right)  }%
}{\displaystyle\int_{T_{m+1}}^{T_{m}}}\left\Vert u\left(  \cdot,t\right)
\right\Vert _{L^{p}\left(  \omega\right)  }dt+\varepsilon\left\Vert u\left(
\cdot,T_{m+2}\right)  \right\Vert _{L^{2}\left(  \Omega\right)  }\text{ ,}%
\]
that is,
\[%
\begin{array}
[c]{ll}
& \quad\varepsilon^{\frac{1-\beta}{\beta}}e^{-\frac{\left(  1+1/\gamma\right)
K}{\left(  z-1\right)  ^{\gamma}T^{\gamma}}z^{\gamma\left(  m+2\right)  }%
}\left\Vert u\left(  \cdot,T_{m}\right)  \right\Vert _{L^{2}\left(
\Omega\right)  }-\varepsilon^{\frac{1}{\beta}}e^{-\frac{\left(  1+1/\gamma
\right)  K}{\left(  z-1\right)  ^{\gamma}T^{\gamma}}z^{\gamma\left(
m+2\right)  }}\left\Vert u\left(  \cdot,T_{m+2}\right)  \right\Vert
_{L^{2}\left(  \Omega\right)  }\\
& \leq\displaystyle{\frac{c}{zK^{1/\gamma}}\int_{T_{m+1}}^{T_{m}}}\left\Vert
u\left(  \cdot,t\right)  \right\Vert _{L^{p}\left(  \omega\right)  }dt\text{
.}%
\end{array}
\]
Replacing $m$ by $2m$, we can see that%
\[%
\begin{array}
[c]{ll}
& \quad\varepsilon^{\frac{1-\beta}{\beta}}e^{-\frac{\left(  1+1/\gamma\right)
K}{\left(  z-1\right)  ^{\gamma}T^{\gamma}}z^{\gamma\left(  2m+2\right)  }%
}\left\Vert u\left(  \cdot,T_{2m}\right)  \right\Vert _{L^{2}\left(
\Omega\right)  }-\varepsilon^{\frac{1}{\beta}}e^{-\frac{\left(  1+1/\gamma
\right)  K}{\left(  z-1\right)  ^{\gamma}T^{\gamma}}z^{\gamma\left(
2m+2\right)  }}\left\Vert u\left(  \cdot,T_{2m+2}\right)  \right\Vert
_{L^{2}\left(  \Omega\right)  }\\
& \leq\displaystyle{\frac{c}{zK^{1/\gamma}}\int_{T_{2m+1}}^{T_{2m}}}\left\Vert
u\left(  \cdot,t\right)  \right\Vert _{L^{p}\left(  \omega\right)  }dt\text{
.}%
\end{array}
\]
We write $A_{m}=e^{-\frac{\left(  1+1/\gamma\right)  K}{\left(  z-1\right)
^{\gamma}T^{\gamma}}z^{\gamma\left(  2m+2\right)  }}$ and choose
$\varepsilon=A_{m}$, in order to get
\[
A_{m}^{\frac{1}{\beta}}\left\Vert u\left(  \cdot,T_{2m}\right)  \right\Vert
_{L^{2}\left(  \Omega\right)  }-A_{m}^{1+\frac{1}{\beta}}\left\Vert u\left(
\cdot,T_{2m+2}\right)  \right\Vert _{L^{2}\left(  \Omega\right)  }%
\leq{\displaystyle\frac{c}{zK^{1/\gamma}}\int_{T_{2m+1}}^{T_{2m}}}\left\Vert
u\left(  \cdot,t\right)  \right\Vert _{L^{p}\left(  \omega\right)  }dt\text{
.}%
\]
Our task is to have
\[
A_{m}^{1+\frac{1}{\beta}}=A_{m+1}^{\frac{1}{\beta}}\text{ , that is,
}e^{-\frac{\left(  1+1/\gamma\right)  K}{\left(  z-1\right)  ^{\gamma
}T^{\gamma}}z^{2\gamma m}z^{2\gamma}\left(  1+\frac{1}{\beta}\right)
}=e^{-\frac{\left(  1+1/\gamma\right)  K}{\left(  z-1\right)  ^{\gamma
}T^{\gamma}}z^{2\gamma m}z^{4\gamma}\frac{1}{\beta}}\text{ ,}%
\]
in order to get, with $X_{m}=A_{m}^{\frac{1}{\beta}}\left\Vert u\left(
\cdot,T_{2m}\right)  \right\Vert _{L^{2}\left(  \Omega\right)  }$,
\[
X_{m}-X_{m+1}\leq{\displaystyle\frac{c}{zK^{1/\gamma}}\int_{T_{2m+1}}^{T_{2m}%
}}\left\Vert u\left(  \cdot,t\right)  \right\Vert _{L^{p}\left(
\omega\right)  }dt\text{ .}%
\]
To this end, we take $z^{2\gamma}\frac{1}{\beta}=1+\frac{1}{\beta}$. It
remains to sum the telescoping series from $m=0$ to $+\infty$ to complete the
proof of Lemma \ref{lemma4.1} and to find that
\[
\left\Vert u\left(  \cdot,T\right)  \right\Vert _{L^{2}\left(  \Omega\right)
}\leq{\displaystyle}\frac{c}{zK^{1/\gamma}}e^{\frac{\left(  1+1/\gamma\right)
K}{T^{\gamma}}C_{\beta}}{\displaystyle\int\nolimits_{0}^{T}}\left\Vert
u\left(  \cdot,t\right)  \right\Vert _{L^{p}\left(  \omega\right)  }dt
\]
with $C_{\beta}={\displaystyle}\frac{\beta+1}{\beta\left[  \left(
\beta+1\right)  ^{\frac{1}{2\gamma}}-1\right]  ^{\gamma}}$.

\bigskip

With the help of Lemma \ref{lemma4.1} and the analysis done in Section 2 for a
convex domain $\Omega\subset\mathbb{R}^{n}$ or a star-shaped domain with
respect to $x_{0}\in\Omega$, we are ready to show Theorem \ref{theorem1.1}. It
suffices to prove the observation at one point of Lemma \ref{lemma4.1} with
$\gamma=1$ and $p=2$.

\bigskip

Let $\hbar>0$. Set
\[
\Phi\left(  x,t\right)  =-\frac{\left\vert x-x_{0}\right\vert ^{2}}{4\left(
T-t+\hbar\right)  }\text{ .}%
\]
The differential inequalities are (see Part \ref{part3.2} of Section 3 and its
example):
\[
\frac{1}{2}\frac{d}{dt}\int_{\Omega}\left\vert u\left(  x,t\right)
\right\vert ^{2}e^{\Phi\left(  x,t\right)  }dx+\mathbf{N}\left(  t\right)
\int_{\Omega}\left\vert u\left(  x,t\right)  \right\vert ^{2}e^{\Phi\left(
x,t\right)  }dx=0\text{ ;}%
\]
Since $\Omega$ is convex or star-shaped w.r.t. $x_{0}$,
\[
\frac{d}{dt}\mathbf{N}\left(  t\right)  \leq\frac{1}{T-t+\hbar}\mathbf{N}%
\left(  t\right)  \text{ .}%
\]
By solving such differential inequalities, we have: For any $0<t_{1}%
<t_{2}<t_{3}\leq T$,\textit{ }%
\[
\left(  \displaystyle\int_{\Omega}\left\vert u\left(  x,t_{2}\right)
\right\vert ^{2}e^{\Phi\left(  x,t_{2}\right)  }dx\right)  ^{1+M}%
\leq\displaystyle\int_{\Omega}\left\vert u\left(  x,t_{3}\right)  \right\vert
^{2}e^{\Phi\left(  x,t_{3}\right)  }dx\left(  \displaystyle\int_{\Omega
}\left\vert u\left(  x,t_{1}\right)  \right\vert ^{2}e^{\Phi\left(
x,t_{1}\right)  }dx\right)  ^{M}%
\]
where%
\[
M=\frac{-\text{ln}\left(  T-t_{3}+\hbar\right)  +\text{ln}\left(
T-t_{2}+\hbar\right)  }{-\text{ln}\left(  T-t_{2}+\hbar\right)  +\text{ln}%
\left(  T-t_{1}+\hbar\right)  }\text{ .}%
\]
Choose $t_{3}=T$, $t_{2}=T-\ell\hbar$, $t_{1}=T-2\ell\hbar$ with $0<2\ell
\hbar<T$ and $\ell>1$, and denote
\[
M_{\ell}=\frac{\text{ln}\left(  \ell+1\right)  }{\text{ln}\left(  \frac
{2\ell+1}{\ell+1}\right)  }\text{ ,}%
\]
then
\[
\left(  \int_{\Omega}\left\vert u\left(  x,T-\ell\hbar\right)  \right\vert
^{2}e^{-\frac{\left\vert x-x_{0}\right\vert ^{2}}{4\left(  \ell+1\right)
\hbar}}dx\right)  ^{1+M_{\ell}}\leq\left(  \int_{\Omega}\left\vert u\left(
x,0\right)  \right\vert ^{2}dx\right)  ^{M_{\ell}}\int_{\Omega}\left\vert
u\left(  x,T\right)  \right\vert ^{2}e^{-\frac{\left\vert x-x_{0}\right\vert
^{2}}{4\hbar}}dx
\]
which implies%
\[
\left(  \int_{\Omega}\left\vert u\left(  x,T\right)  \right\vert
^{2}dx\right)  ^{1+M_{\ell}}\leq e^{\frac{R^{2}\left(  1+M_{\ell}\right)
}{4\left(  \ell+1\right)  \hbar}}\left(  \int_{\Omega}\left\vert u\left(
x,0\right)  \right\vert ^{2}dx\right)  ^{M_{\ell}}\int_{\Omega}\left\vert
u\left(  x,T\right)  \right\vert ^{2}e^{-\frac{\left\vert x-x_{0}\right\vert
^{2}}{4\hbar}}dx\text{ .}%
\]
Here and throughout the proof of Theorem \ref{theorem1.1}, $R:=\underset
{x\in\overline{\Omega}}{\text{max}}\left\vert x-x_{0}\right\vert $. Next, we
split $\displaystyle\int_{\Omega}\left\vert u\left(  x,T\right)  \right\vert
^{2}e^{-\frac{\left\vert x-x_{0}\right\vert ^{2}}{4\hbar}}dx$ into two parts:
With $B_{x_{0},r}:=\left\{  x;\left\vert x-x_{0}\right\vert <r\right\}
\Subset\Omega$ where $r<R$, we can see that
\[
\int_{\Omega}\left\vert u\left(  x,T\right)  \right\vert ^{2}e^{-\frac
{\left\vert x-x_{0}\right\vert ^{2}}{4\hbar}}dx\leq\int_{B_{x_{0},r}%
}\left\vert u\left(  x,T\right)  \right\vert ^{2}dx+e^{-\frac{r^{2}}{4\hbar}%
}\int_{\Omega}\left\vert u\left(  x,0\right)  \right\vert ^{2}dx\text{ .}%
\]
Therefore, taking the above estimates into consideration yields that%
\[%
\begin{array}
[c]{ll}
& \quad\left(  \displaystyle\int_{\Omega}\left\vert u\left(  x,T\right)
\right\vert ^{2}dx\right)  ^{1+M_{\ell}}\\
& \leq\left(  \displaystyle\int_{\Omega}\left\vert u\left(  x,0\right)
\right\vert ^{2}dx\right)  ^{M_{\ell}}\\
& \quad\times\left(  e^{\frac{R^{2}\left(  1+M_{\ell}\right)  }{4\left(
\ell+1\right)  \hbar}}\displaystyle\int_{B_{x_{0},r}}\left\vert u\left(
x,T\right)  \right\vert ^{2}dx+e^{\frac{R^{2}\left(  1+M_{\ell}\right)
}{4\left(  \ell+1\right)  \hbar}}e^{-\frac{r^{2}}{4\hbar}}\displaystyle\int
_{\Omega}\left\vert u\left(  x,0\right)  \right\vert ^{2}dx\right)  \text{ .}%
\end{array}
\]
But for $\ell>1$, $\frac{\left(  1+M_{\ell}\right)  }{4\left(  \ell+1\right)
}\leq\frac{1}{2\left(  \ell+1\right)  }\frac{\text{ln}\left(  \ell+1\right)
}{\text{ln}\left(  \frac{2\ell+1}{\ell+1}\right)  }\leq\frac{1}{2\text{ln}%
\left(  3/2\right)  }\frac{\text{ln}\left(  2\ell\right)  }{\ell}\leq
\frac{2^{\varepsilon}}{2\varepsilon\text{ln}\left(  3/2\right)  }\frac{1}%
{\ell^{1-\varepsilon}}$ $\forall\varepsilon\in\left(  0,1\right)  $. Our
choice of $\ell$:
\[
\ell:=\left(  \frac{R^{2}}{r^{2}}\right)  ^{1/\left(  1-\varepsilon\right)
}\left(  \frac{2^{2+\varepsilon}}{\varepsilon\text{ln}\left(  3/2\right)
}\right)  ^{1/\left(  1-\varepsilon\right)  }\text{ }\forall\varepsilon
\in\left(  0,1\right)
\]
gives
\[
\frac{R^{2}\left(  1+M_{\ell}\right)  }{4\left(  \ell+1\right)  \hbar}%
\leq\frac{r^{2}}{8\hbar}\text{ .}%
\]
One the one hand, it implies that for any $2\ell\hbar<T$%
\[%
\begin{array}
[c]{ll}%
\left(  \displaystyle\int_{\Omega}\left\vert u\left(  x,T\right)  \right\vert
^{2}dx\right)  ^{1+M_{\ell}} & \leq\left(  \displaystyle\int_{\Omega
}\left\vert u\left(  x,0\right)  \right\vert ^{2}dx\right)  ^{M_{\ell}}\\
& \quad\times\left(  e^{\frac{r^{2}}{8\hbar}}\displaystyle\int_{B_{x_{0},r}%
}\left\vert u\left(  x,T\right)  \right\vert ^{2}dx+e^{-\frac{r^{2}}{8\hbar}%
}\displaystyle\int_{\Omega}\left\vert u\left(  x,0\right)  \right\vert
^{2}dx\right)  \text{ .}%
\end{array}
\]
On the other hand, $\left\Vert u\left(  \cdot,T\right)  \right\Vert
\leq\left\Vert u\left(  \cdot,0\right)  \right\Vert $ and for any $2\ell
\hbar\geq T$, $1\leq e^{\frac{r^{2}\ell}{4T}}e^{-\frac{r^{2}}{8\hbar}}$.
Therefore we conclude that for any $\hbar>0$,%
\[%
\begin{array}
[c]{ll}%
\left(  \displaystyle\int_{\Omega}\left\vert u\left(  x,T\right)  \right\vert
^{2}dx\right)  ^{1+M_{\ell}} & \leq\left(  \displaystyle\int_{\Omega
}\left\vert u\left(  x,0\right)  \right\vert ^{2}dx\right)  ^{M_{\ell}}\\
& \quad\times\left(  e^{\frac{r^{2}}{8\hbar}}\displaystyle\int_{B_{x_{0},r}%
}\left\vert u\left(  x,T\right)  \right\vert ^{2}dx+e^{\frac{r^{2}\ell}{4T}%
}e^{-\frac{r^{2}}{8\hbar}}\displaystyle\int_{\Omega}\left\vert u\left(
x,0\right)  \right\vert ^{2}dx\right)  \text{ .}%
\end{array}
\]
Finally, we choose $\hbar>0$ such that
\[
e^{\frac{r^{2}}{8\hbar}}:=2e^{\frac{r^{2}\ell}{4T}}\left(  \frac
{\displaystyle\int_{\Omega}\left\vert u\left(  x,0\right)  \right\vert ^{2}%
dx}{\displaystyle\int_{\Omega}\left\vert u\left(  x,T\right)  \right\vert
^{2}dx}\right)  ^{1+M_{\ell}}%
\]
in order that
\[
\displaystyle\int_{\Omega}\left\vert u\left(  x,T\right)  \right\vert
^{2}dx\leq\left(  4e^{\frac{r^{2}\ell}{4T}}\displaystyle\int_{B_{x_{0},r}%
}\left\vert u\left(  x,T\right)  \right\vert ^{2}dx\right)  ^{\frac
{1}{2\left(  1+M_{\ell}\right)  }}\left(  \displaystyle\int_{\Omega}\left\vert
u\left(  x,0\right)  \right\vert ^{2}dx\right)  ^{\frac{1+2M_{\ell}}{2\left(
1+M_{\ell}\right)  }}\text{ ,}%
\]
that is,
\[
\left\Vert u\left(  x,T\right)  \right\Vert _{L^{2}\left(  \Omega\right)
}\leq\left(  2e^{\frac{r^{2}\ell}{8T}}\left\Vert u\left(  x,T\right)
\right\Vert _{L^{2}\left(  B_{x_{0},r}\right)  }\right)  ^{\frac{1}{2\left(
1+M_{\ell}\right)  }}\left(  \left\Vert u\left(  x,0\right)  \right\Vert
_{L^{2}\left(  \Omega\right)  }\right)  ^{\frac{1+2M_{\ell}}{2\left(
1+M_{\ell}\right)  }}\text{ .}%
\]
Now, we can apply Lemma \ref{lemma4.1} with $\gamma=1$, $p=2$ and Lemma A in
Appendix (see page \pageref{appendix}) with $c=2$, $K=\frac{r^{2}\ell}{8}$,
$\beta=\frac{1}{2\left(  1+M_{\ell}\right)  }$. Consequently, we obtain that
\[
\left\Vert e^{T\Delta}u_{0}\right\Vert _{L^{2}\left(  \Omega\right)  }%
\leq\frac{16}{r^{2}\ell}e^{\frac{r^{2}\ell}{4T}C_{\beta}}\int_{0}%
^{T}\left\Vert e^{t\Delta}u_{0}\right\Vert _{L^{2}\left(  B_{x_{0},r}\right)
}dt\text{ }%
\]
and
\[
\sum\limits_{\lambda_{i}\leq\lambda}\left\vert a_{i}\right\vert ^{2}%
\leq4e^{4\sqrt{\lambda\left(  1+2M_{\ell}\right)  \frac{r^{2}\ell}{8}}}%
{\int\nolimits_{B_{x_{0},r}}}\left\vert \sum\limits_{\lambda_{i}\leq\lambda
}a_{i}e_{i}\left(  x\right)  \right\vert ^{2}dx\text{ .}%
\]
We can see that $C_{\beta}\leq$constant$\left(  M_{\ell}\right)  ^{2}$. By the
definition of $M_{\ell}$ and of $\ell$, we have $M_{\ell}\leq\frac
{\text{ln}\left(  \ell+1\right)  }{\text{ln}\left(  3/2\right)  }\leq
$constant$\frac{1}{r^{\varepsilon}}$. Therefore, we conclude that
\[
\left\Vert e^{T\Delta}u_{0}\right\Vert _{L^{2}\left(  \Omega\right)  }\leq
K_{\varepsilon}e^{\frac{K_{\varepsilon}}{r^{\varepsilon}}\frac{1}{T}}\int
_{0}^{T}\left\Vert e^{t\Delta}u_{0}\right\Vert _{L^{2}\left(  B_{x_{0}%
,r}\right)  }dt\text{ }%
\]
and%
\[
\sum\limits_{\lambda_{i}\leq\lambda}\left\vert a_{i}\right\vert ^{2}%
\leq4e^{\frac{K_{\varepsilon}}{r^{\varepsilon}}\sqrt{\lambda}}{\int
\nolimits_{B_{x_{0},r}}}\left\vert \sum\limits_{\lambda_{i}\leq\lambda}%
a_{i}e_{i}\left(  x\right)  \right\vert ^{2}dx\text{ .}%
\]
This completes the proof of Theorem \ref{theorem1.1}.

\bigskip

\bigskip

\section{Proof of Theorem 1.2}

\bigskip

\bigskip

Let $n\geq3$ and consider a $C^{2}$ bounded domain $\Omega\subset
\mathbb{R}^{n}$ such that $0\in\Omega$, and let $\omega\subset\Omega$ be a
nonempty open set. To simplify the presentation, we assume that $0\notin
\overline{\omega}$, that can always be done, taking if necessary a smaller
set. Let $R_{0}=\left(  \frac{4}{3}\right)  ^{3/2}\approx1.53$. We also assume
that the unit ball $\overline{B}_{0,R_{0}}$ is included in $\Omega$ and
$\overline{B}_{0,R_{0}}\cap\overline{\omega}$ is empty. This can always be
done by a scaling argument.

\bigskip

We are interested in the following heat equation with an inverse square
potential%
\[
\left\{
\begin{array}
[c]{ll}%
{\partial}_{t}u-\Delta u-\frac{\mu}{\left\vert x\right\vert ^{2}}u=0\text{ ,}
& \quad\text{in}~\Omega\times\left(  0,T\right)  \text{ ,}\\
u=0\text{ ,} & \quad\text{on}~\partial\Omega\times\left(  0,T\right)  \text{
,}\\
u\left(  \cdot,0\right)  =u_{0}\text{ , } & \quad\text{in}~\Omega\text{ ,}%
\end{array}
\right.
\]
where $u_{0}\in L^{2}\left(  \Omega\right)  $, $T>0$ and $\mu<\mu^{\ast
}\left(  n\right)  :=\frac{\left(  n-2\right)  ^{2}}{4}$. It is well-known
that this is a well-posed problem \cite{VZ}. In particular, $u\in C\left(
\left[  0,T\right]  ;L^{2}\left(  \Omega\right)  \right)  \cap L^{2}\left(
0,T;H_{0}^{1}\left(  \Omega\right)  \right)  $ and for any $t\in\left(
0,T\right]  $, we have
\[
\int_{\Omega}\left\vert u\left(  x,t\right)  \right\vert ^{2}dx\leq
\int_{\Omega}\left\vert u_{0}\left(  x\right)  \right\vert ^{2}dx\text{ ,}%
\]
and the regularizing effect%
\[
\int_{\Omega}\left\vert \nabla u\left(  x,t\right)  \right\vert ^{2}%
dx\leq\frac{C}{t}\int_{\Omega}\left\vert u_{0}\left(  x\right)  \right\vert
^{2}dx\text{ .}%
\]

\bigskip

Applying Lemma A in Appendix (see page \pageref{appendix}), we obtain the
following result.

\bigskip

\begin{lemm}
\label{lemma5.1} Let $\beta\in\left(  0,1\right)  $ and $c$, $K>0$. Suppose
that for any $u_{0}\in L^{2}\left(  \Omega\right)  $ and any $T>0$,
\[
\left\Vert u\left(  \cdot,T\right)  \right\Vert _{L^{2}\left(  \Omega\right)
}\leq\left(  ce^{\frac{K}{T}}\left\Vert u\left(  \cdot,T\right)  \right\Vert
_{L^{1}\left(  \omega\right)  }\right)  ^{\beta}\left(  \left\Vert
u_{0}\right\Vert _{L^{2}\left(  \Omega\right)  }\right)  ^{1-\beta}\text{ .}%
\]
Then for any $\left(  a_{j}\right)  _{j\geq1}\in\mathbb{R}$ and any
$\lambda>0$, one has%
\[
\sqrt{\sum\limits_{\lambda_{j}\leq\lambda}\left\vert a_{j}\right\vert ^{2}%
}\leq ce^{2\sqrt{\frac{1-\beta}{\beta}K\lambda}}\left\Vert \sum
\limits_{\lambda_{j}\leq\lambda}a_{j}e_{j}\right\Vert _{L^{1}\left(
\omega\right)  }\text{ .}%
\]
Here $\left(  \lambda_{j},e_{j}\right)  $ denotes the eigenbasis of the
Schr\"{o}dinger operator $-\Delta-\frac{\mu}{\left\vert x\right\vert ^{2}}$
with Dirichlet boundary condition
\[
\left\{
\begin{array}
[c]{rll}%
-\Delta e_{j}-\frac{\mu}{\left\vert x\right\vert ^{2}}e_{j} & =\lambda
_{j}e_{j}\text{ ,} & \text{ in }\Omega\text{ ,}\\
e_{j} & =0\text{ ,} & \text{ on }\partial\Omega\text{ .}%
\end{array}
\right.
\]

\end{lemm}

\bigskip

Now we are ready to prove Theorem \ref{theorem1.2}. It suffices to check the
assumption of the above Lemma \ref{lemma5.1} when
\[
\mu\leq\left\vert
\begin{array}
[c]{ll}%
\frac{7}{2\cdot3^{3}}:=\mu^{\ast}\left(  3\right)  -\frac{13}{4\cdot3^{3}%
}\text{ ,} & \text{if }n=3\text{ ,}\\
\frac{1}{4}\left(  n-1\right)  \left(  n-3\right)  :=\mu^{\ast}\left(
n\right)  -\frac{1}{4}\text{ ,} & \text{if }n\geq4\text{ .}%
\end{array}
\right.
\]

\bigskip

Recall that $R_{0}=\left(  \frac{4}{3}\right)  ^{3/2}\approx1.53$ and let
$\delta=\frac{R_{0}-1}{2}\approx0.26$. We have assumed that the unit ball
$B_{0,R_{0}}$ is included in $\Omega$. Let $\chi\in C_{0}^{\infty}\left(
B_{0,R_{0}}\right)  $ be a cut-off function satisfying $0\leq\chi\leq1$ and
$\chi=1$ on $\left\{  x;\left\vert x\right\vert \leq1+3\delta/2\right\}  $.
Introduce $z=\chi u$. It solves
\[
\partial_{t}z-\Delta z-\frac{\mu}{\left\vert x\right\vert ^{2}}z=g:=-2\nabla
\chi\cdot\nabla u-\Delta\chi u\text{ ,}%
\]
and furthermore, $z_{\left\vert \partial B_{0,R_{0}}\right.  }=\partial_{\nu
}z_{\left\vert \partial B_{0,R_{0}}\right.  }=0$. Let $\Phi$ be a sufficiently
smooth function of $\left(  x,t\right)  \in\mathbb{R}^{n}\times\mathbb{R}_{t}$
and set
\[
f\left(  x,t\right)  =z\left(  x,t\right)  e^{\Phi\left(  x,t\right)
/2}\text{ .}%
\]
We look for the equation solved by $f$ by computing $e^{\Phi\left(
x,t\right)  /2}\left(  \partial_{t}-\Delta\right)  \left(  e^{-\Phi\left(
x,t\right)  /2}f\left(  x,t\right)  \right)  $. It gives
\[
\partial_{t}f-\Delta f-\frac{1}{2}f\left(  \partial_{t}\Phi+\frac{1}%
{2}\left\vert \nabla\Phi\right\vert ^{2}-\frac{2\mu}{\left\vert x\right\vert
^{2}}\right)  +\nabla\Phi\cdot\nabla f+\frac{1}{2}\Delta\Phi f=e^{\Phi
/2}g\quad\text{in}~B_{0,R_{0}}\times\left(  0,T\right)  \text{ ,}%
\]
and furthermore, $f_{\left\vert \partial B_{x_{0},R_{0}}\right.  }%
=\partial_{\nu}f_{\left\vert \partial B_{x_{0},R_{0}}\right.  }=0$. Introduce
\[
\left\{
\begin{array}
[c]{ll}%
\mathcal{A}f=-\nabla\Phi\cdot\nabla f-\frac{1}{2}\Delta\Phi f\text{ ,} & \\
\mathcal{S}f=\Delta f+\left(  \frac{1}{2}\partial_{t}\Phi+\frac{1}%
{4}\left\vert \nabla\Phi\right\vert ^{2}+\frac{\mu}{\left\vert x\right\vert
^{2}}\right)  f\text{ .} &
\end{array}
\right.
\]
Then, it holds%
\[
\left\{
\begin{array}
[c]{ll}%
\left\langle \mathcal{A}f,v\right\rangle _{0}=-\left\langle \mathcal{A}%
v,f\right\rangle _{0}\text{ ,} & \\
\left\langle \mathcal{S}f,v\right\rangle _{0}=\left\langle \mathcal{S}%
v,f\right\rangle _{0}\text{ for any }v\in H_{0}^{1}\left(  B_{0,R_{0}}\right)
\text{ ,} &
\end{array}
\right.
\]
where $\left\langle \cdot,\cdot\right\rangle _{0}$ is the usual scalar product
in $L^{2}\left(  B_{0,R_{0}}\right)  $ and $\left\Vert \cdot\right\Vert _{0}$
will denote the corresponding norm. Furthermore, we have
\[
\partial_{t}f-\mathcal{S}f-\mathcal{A}f=e^{\Phi/2}g\text{ .}%
\]
Multiplying by $f$ the above equation, integrating over $B_{0,R_{0}}$, it
follows that%
\[
\frac{1}{2}\frac{d}{dt}\left\Vert f\right\Vert _{0}^{2}+\left\langle
-\mathcal{S}f,f\right\rangle _{0}=\left\langle f,e^{\Phi/2}g\right\rangle
_{0}\text{ .}%
\]
Introduce the frequency function $t\mapsto\mathbf{N}\left(  t\right)  $
defined by%
\[
\mathbf{N}=\frac{\left\langle -\mathcal{S}f,f\right\rangle _{0}}{\left\Vert
f\right\Vert _{0}^{2}}\text{ .}%
\]
Then, we have
\[
\frac{1}{2}\frac{d}{dt}\left\Vert f\right\Vert _{0}^{2}+\mathbf{N}\left\Vert
f\right\Vert _{0}^{2}=\left\langle e^{\Phi/2}g,f\right\rangle _{0}\text{ ,}%
\]
and the derivative of $\mathbf{N}$ satisfies (see Step \ref{step2.2.2} in
Section 2):%
\[
\frac{d}{dt}\mathbf{N}\leq\frac{1}{\left\Vert f\right\Vert _{0}^{2}%
}\left\langle -\left(  \mathcal{S}^{\prime}+\left[  \mathcal{S},\mathcal{A}%
\right]  \right)  f,f\right\rangle _{0}+\frac{1}{\left\Vert f\right\Vert
_{0}^{2}}\left\Vert e^{\Phi/2}g\right\Vert _{0}^{2}\text{ .}%
\]
Notice that the boundary terms have vanished since $f_{\left\vert \partial
B_{x_{0},R_{0}}\right.  }=\partial_{\nu}f_{\left\vert \partial B_{x_{0},R_{0}%
}\right.  }=0$.

\bigskip

The estimate of $\displaystyle\frac{\left\Vert e^{\Phi/2}g\right\Vert _{0}%
^{2}}{\left\Vert f\right\Vert _{0}^{2}}$ can be obtained in a similar way than
in Step \ref{step2.2.4} and Step \ref{step2.2.5} of Section 2. Indeed, first,
we check that
\[
\frac{d}{dt}\int_{\Omega}\left\vert u\left(  x,t\right)  \right\vert
^{2}e^{\xi\left(  x,t\right)  }dx\leq0\text{ }\quad\text{with }\xi\left(
x,t\right)  =-\frac{\left\vert x\right\vert ^{2}}{2\left(  T-t+\epsilon
\right)  }\text{ ,}%
\]
as follows:%
\[%
\begin{array}
[c]{ll}
& \quad\displaystyle\frac{1}{2}\frac{d}{dt}\displaystyle\int_{\Omega
}\left\vert u\right\vert ^{2}e^{\xi}dx\\
& =\displaystyle\int_{\Omega}u\partial_{t}ue^{\xi}dx+\displaystyle\frac{1}%
{2}\int_{\Omega}\left\vert u\right\vert ^{2}\partial_{t}\xi e^{\xi}dx\\
& =\displaystyle\int_{\Omega}u\left(  \Delta u+\frac{\mu}{\left\vert
x\right\vert ^{2}}u\right)  e^{\xi}dx+\displaystyle\frac{1}{2}\int_{\Omega
}\left\vert u\right\vert ^{2}\partial_{t}\xi e^{\xi}dx\\
& =-\displaystyle\int_{\Omega}\left\vert \nabla u\right\vert ^{2}e^{\xi
}dx-\displaystyle\int_{\Omega}u\nabla u\cdot\nabla\xi e^{\xi}%
dx+\displaystyle\int_{\Omega}\frac{\mu}{\left\vert x\right\vert ^{2}%
}\left\vert ue^{\xi/2}\right\vert ^{2}dx+\displaystyle\frac{1}{2}\int_{\Omega
}\left\vert u\right\vert ^{2}\partial_{t}\xi e^{\xi}dx\\
& \leq-\displaystyle\int_{\Omega}\left\vert \nabla u\right\vert ^{2}e^{\xi
}dx-\displaystyle\int_{\Omega}u\nabla u\cdot\nabla\xi e^{\xi}%
dx+\displaystyle\int_{\Omega}\left\vert \nabla\left(  ue^{\xi/2}\right)
\right\vert ^{2}dx+\displaystyle\frac{1}{2}\int_{\Omega}\left\vert
u\right\vert ^{2}\partial_{t}\xi e^{\xi}dx\\
& =-\displaystyle\int_{\Omega}\left\vert \nabla u\right\vert ^{2}e^{\xi
}dx-\displaystyle\int_{\Omega}u\nabla u\cdot\nabla\xi e^{\xi}%
dx+\displaystyle\int_{\Omega}\left\vert \nabla ue^{\xi/2}+u\frac{1}{2}%
\nabla\xi e^{\xi/2}\right\vert ^{2}dx+\displaystyle\frac{1}{2}\int_{\Omega
}\left\vert u\right\vert ^{2}\partial_{t}\xi e^{\xi}dx\\
& =\displaystyle\frac{1}{4}\int_{\Omega}\left\vert u\right\vert ^{2}\left\vert
\nabla\xi\right\vert ^{2}e^{\xi}dx+\displaystyle\frac{1}{2}\int_{\Omega
}\left\vert u\right\vert ^{2}\partial_{t}\xi e^{\xi}dx=0\text{ ,}%
\end{array}
\]
where in the fifth line we used Hardy inequality. Therefore, we have
\[
\int_{\Omega}\left\vert u\left(  x,T\right)  \right\vert ^{2}e^{-\frac
{\left\vert x-x_{0}\right\vert ^{2}}{2\epsilon}}dx\leq\int_{\Omega}\left\vert
u\left(  x,t\right)  \right\vert ^{2}e^{-\frac{\left\vert x-x_{0}\right\vert
^{2}}{2\left(  T-t+\epsilon\right)  }}dx\text{ ,}%
\]
which implies that Lemma \ref{lemma2.2} is still true for any $u$ solution of
the heat equation with an inverse square potential.

Next we can estimate $\displaystyle\frac{\left\Vert e^{\Phi/2}g\right\Vert
_{0}^{2}}{\left\Vert f\right\Vert _{0}^{2}}$ as follows.%
\[%
\begin{array}
[c]{ll}
& \quad\displaystyle\frac{\left\Vert e^{\Phi/2}g\left(  \cdot,t\right)
\right\Vert _{0}^{2}}{\left\Vert f\left(  \cdot,t\right)  \right\Vert _{0}%
^{2}}\\
& \leq\frac{\displaystyle\int_{\Omega\cap\left\{  \left(  1+3\delta/2\right)
R\leq\left\vert x\right\vert \leq R_{0}\right\}  }\left\vert -2\nabla\chi
\cdot\nabla u\left(  x,t\right)  -\Delta\chi u\left(  x,t\right)  \right\vert
^{2}e^{\Phi\left(  x,t\right)  }dx}{\displaystyle\int_{\Omega\cap B_{0,\left(
1+\delta\right)  R}}\left\vert u\left(  x,t\right)  \right\vert ^{2}%
e^{\Phi\left(  x,t\right)  }dx}\\
& \leq\text{exp}\left[  -\underset{\left\vert x\right\vert \leq\left(
1+\delta\right)  R}{\text{min}}\Phi\left(  x,t\right)  +\underset{\left(
1+3\delta/2\right)  R\leq\left\vert x\right\vert \leq R_{0}}{\text{max}}%
\Phi\left(  x,t\right)  \right]  \displaystyle\frac{C\left(  1+\frac{1}%
{t}\right)  \left\Vert u\left(  \cdot,0\right)  \right\Vert ^{2}}{\left\Vert
u\left(  \cdot,t\right)  \right\Vert _{L^{2}\left(  \Omega\cap B_{0,\left(
1+\delta\right)  R}\right)  }^{2}}\\
& \leq\text{exp}\left[  -\underset{\left\vert x\right\vert \leq\left(
1+\delta\right)  R}{\text{min}}\Phi\left(  x,t\right)  +\underset{\left(
1+3\delta/2\right)  R\leq\left\vert x\right\vert \leq R_{0}}{\text{max}}%
\Phi\left(  x,t\right)  \right]  C\left(  1+\frac{1}{t}\right)  e^{\left(
1+\delta\right)  \delta\frac{R^{2}}{2\theta}}%
\end{array}
\]
as long as $T/2\leq T-\theta\leq t\leq T$, where in the third line we used the
regularizing effect of a gradient term for the solution $u$ of the heat
equation with an inverse square potential. Therefore, the conclusion of Step
\ref{step2.2.5} still holds: Under the assumptions of Step \ref{step2.2.5}, we
have
\[
\frac{\left\Vert e^{\Phi/2}g\left(  \cdot,t\right)  \right\Vert _{0}^{2}%
}{\left\Vert f\left(  \cdot,t\right)  \right\Vert _{0}^{2}}\leq C\left(
1+\frac{1}{t}\right)  \text{ .}%
\]

\bigskip

The difficulty with the heat equation with an inverse square potential comes
with the estimate of $\left\langle -\left(  \mathcal{S}^{\prime}+\left[
\mathcal{S},\mathcal{A}\right]  \right)  f,f\right\rangle _{0}%
-\displaystyle\frac{1}{\Upsilon}\left\langle -\mathcal{S}f,f\right\rangle
_{0}$. Notice also that the treatment far from the point $0\in\Omega$ where
the inverse square potential have its singularities can be done in the same
way than for the heat equation with a potential in $L^{\infty}\left(
\Omega\times\left(  0,T\right)  \right)  $ (see \cite{PWZ}). Our main task is
to treat the assumptions of Step \ref{step2.2.5} and Step \ref{step2.2.7} of
Section 2, carefully with a suitable choice of $\Phi$ (see also Part
\ref{part3.3} of Section 3).

\bigskip

We claim that:%
\[%
\begin{array}
[c]{ll}
& \quad\left\langle -\left(  \mathcal{S}^{\prime}+\left[  \mathcal{S}%
,\mathcal{A}\right]  \right)  f,f\right\rangle _{0}\\
& =-2\displaystyle\int_{\Omega\cap B_{0,R_{0}}}\nabla f\cdot\nabla^{2}%
\Phi\nabla fdx+2\mu\displaystyle\int_{\Omega\cap B_{0,R_{0}}}\nabla\Phi
\cdot\frac{x}{\left\vert x\right\vert ^{4}}\left\vert f\right\vert ^{2}dx\\
& \quad+\displaystyle\frac{1}{2}\int_{\Omega\cap B_{0,R_{0}}}\left(
\Delta^{2}\Phi-\partial_{t}^{2}\Phi-2\nabla\Phi\cdot\nabla\partial_{t}%
\Phi-\nabla\Phi\cdot\nabla^{2}\Phi\nabla\Phi\right)  \left\vert f\right\vert
^{2}dx\text{ .}%
\end{array}
\]

Proof of the claim .- First, $\mathcal{S}^{\prime}f=\left(  \frac{1}%
{2}\partial_{t}^{2}\Phi+\frac{1}{2}\nabla\Phi\cdot\nabla\partial_{t}%
\Phi\right)  f$. Next, we compute $\left[  \mathcal{S},\mathcal{A}\right]
f=\mathcal{SA}f-\mathcal{AS}f$ and get
\[
\left[  \mathcal{S},\mathcal{A}\right]  f=-2\partial_{i}\nabla\Phi
\cdot\partial_{i}\nabla f+\nabla\Phi\cdot\nabla\left(  \frac{1}{2}\partial
_{t}\Phi+\frac{1}{4}\left\vert \nabla\Phi\right\vert ^{2}+\frac{\mu
}{\left\vert x\right\vert ^{2}}\right)  f-\frac{1}{2}\Delta^{2}\Phi
f-2\Delta\nabla\Phi\cdot\nabla f\text{ }%
\]
(which corresponds to the formula in the claim of Part \ref{part3.1} in
Section 3 with $\eta=\frac{1}{2}\partial_{t}\Phi+\frac{1}{4}\left\vert
\nabla\Phi\right\vert ^{2}+\frac{\mu}{\left\vert x\right\vert ^{2}}$).
Therefore,
\[%
\begin{array}
[c]{ll}%
-\left(  \mathcal{S}^{\prime}+\left[  \mathcal{S},\mathcal{A}\right]  \right)
f & =2\partial_{i}\nabla\Phi\cdot\partial_{i}\nabla f+2\Delta\nabla\Phi
\cdot\nabla f+\frac{1}{2}\Delta^{2}\Phi f\\
& \quad-\left(  \frac{1}{2}\partial_{t}^{2}\Phi+\frac{1}{2}\nabla\Phi
\cdot\nabla\partial_{t}\Phi+\nabla\Phi\cdot\nabla\left(  \frac{1}{2}%
\partial_{t}\Phi+\frac{1}{4}\left\vert \nabla\Phi\right\vert ^{2}+\frac{\mu
}{\left\vert x\right\vert ^{2}}\right)  \right)  f\\
& =2\partial_{i}\nabla\Phi\cdot\partial_{i}\nabla f+2\Delta\nabla\Phi
\cdot\nabla f+\frac{1}{2}\Delta^{2}\Phi f\\
& \quad-\left(  \frac{1}{2}\partial_{t}^{2}\Phi+\nabla\Phi\cdot\nabla
\partial_{t}\Phi+\frac{1}{2}\nabla\Phi\cdot\nabla^{2}\Phi\nabla\Phi+\nabla
\Phi\cdot\nabla\left(  \frac{\mu}{\left\vert x\right\vert ^{2}}\right)
\right)  f\text{ .}%
\end{array}
\]
Furthermore, by one integration by parts we have
\[
2\left\langle \partial_{i}\nabla\Phi\cdot\partial_{i}\nabla f,f\right\rangle
_{0}=\int_{\Omega\cap B_{0,R_{0}}}\Delta^{2}\Phi\left\vert f\right\vert
^{2}dx-2\int_{\Omega\cap B_{0,R_{0}}}\nabla f\cdot\nabla^{2}\Phi\nabla fdx
\]
and%
\[
2\left\langle \Delta\nabla\Phi\cdot\nabla f,f\right\rangle _{0}=-\int
_{\Omega\cap B_{0,R_{0}}}\Delta^{2}\Phi\left\vert f\right\vert ^{2}dx\text{ .}%
\]
Combining the above equalities yields the desired claim.

\bigskip

Assume that $\Phi\left(  x,t\right)  =\displaystyle\frac{\varphi\left(
x\right)  }{T-t+\hbar}$ and recall that $\Upsilon\left(  t\right)
:=T-t+\hbar$. Then, we can see that%
\[%
\begin{array}
[c]{ll}
& \quad\left\langle -\left(  \mathcal{S}^{\prime}+\left[  \mathcal{S}%
,\mathcal{A}\right]  \right)  f,f\right\rangle _{0}-\displaystyle\frac
{1}{\Upsilon}\left\langle -\mathcal{S}f,f\right\rangle _{0}\\
& =-\displaystyle\frac{1}{\Upsilon}\int_{\Omega\cap B_{0,R_{0}}}\nabla
f\cdot\left(  2\nabla^{2}\varphi+I_{d}\right)  \nabla fdx\\
& \quad+\displaystyle\frac{1}{\Upsilon}\int_{\Omega\cap B_{0,R_{0}}}\left(
\frac{1}{2}\Delta^{2}\varphi+2\mu\nabla\varphi\cdot\frac{x}{\left\vert
x\right\vert ^{4}}+\frac{\mu}{\left\vert x\right\vert ^{2}}\right)  \left\vert
f\right\vert ^{2}dx\\
& \quad-\displaystyle\frac{1}{2\Upsilon^{3}}\int_{\Omega\cap B_{0,R_{0}}%
}\left(  \varphi+\left\vert \nabla\varphi\right\vert ^{2}+\frac{1}{2}%
\nabla\varphi\cdot\left(  2\nabla^{2}\varphi+I_{d}\right)  \nabla
\varphi\right)  \left\vert f\right\vert ^{2}dx\text{ .}%
\end{array}
\]
Indeed,
\[%
\begin{array}
[c]{ll}
& \quad\left\langle -\left(  \mathcal{S}^{\prime}+\left[  \mathcal{S}%
,\mathcal{A}\right]  \right)  f,f\right\rangle _{0}\\
& =-\displaystyle\frac{2}{\Upsilon}\int_{\Omega\cap B_{x_{0},R_{0}}}\nabla
f\cdot\nabla^{2}\varphi\nabla fdx+\displaystyle\frac{1}{2\Upsilon}\int
_{\Omega\cap B_{x_{0},R_{0}}}\Delta^{2}\varphi\left\vert f\right\vert ^{2}dx\\
& \quad+\displaystyle\frac{1}{\Upsilon}\int_{\Omega\cap B_{0,R_{0}}}2\mu
\nabla\varphi\cdot\frac{x}{\left\vert x\right\vert ^{4}}\left\vert
f\right\vert ^{2}dx\\
& \quad-\displaystyle\frac{1}{\Upsilon^{3}}\int_{\Omega\cap B_{x_{0},R_{0}}%
}\left(  \varphi+\left\vert \nabla\varphi\right\vert ^{2}+\frac{1}{2}%
\nabla\varphi\cdot\nabla^{2}\varphi\nabla\varphi\right)  \left\vert
f\right\vert ^{2}dx\text{ ,}%
\end{array}
\]
and%
\[%
\begin{array}
[c]{ll}%
\displaystyle\frac{1}{\Upsilon}\left\langle -\mathcal{S}f,f\right\rangle _{0}
& =\displaystyle\frac{1}{\Upsilon}\int_{\Omega\cap B_{0,R_{0}}}\left(
\left\vert \nabla f\right\vert ^{2}-\frac{\mu}{\left\vert x\right\vert ^{2}%
}\left\vert f\right\vert ^{2}\right)  dx\\
& -\displaystyle\frac{1}{\Upsilon^{3}}\int_{\Omega\cap B_{0,R_{0}}}\left(
\frac{1}{2}\varphi+\frac{1}{4}\left\vert \nabla\varphi\right\vert ^{2}\right)
\left\vert f\right\vert ^{2}dx\text{ .}%
\end{array}
\]

\bigskip

Assume that $\varphi\left(  x\right)  =-a\left\vert x\right\vert
^{2}+b\left\vert x\right\vert ^{s}-c$ for some $a,b,c>0$ and $1\leq s<2$. We
would like to check the assumptions of Section 2 and find the adequat
parameters $a,b,c,s$.

\bigskip

First, we observe that the identity
\[%
\begin{array}
[c]{ll}
& \quad\left\langle -\left(  \mathcal{S}^{\prime}+\left[  \mathcal{S}%
,\mathcal{A}\right]  \right)  f,f\right\rangle _{0}-\displaystyle\frac
{1}{\Upsilon}\left\langle -\mathcal{S}f,f\right\rangle _{0}\\
& =-\displaystyle\frac{1}{\Upsilon}\int_{\Omega\cap B_{0,R_{0}}}\nabla
f\cdot\left(  2\nabla^{2}\varphi+I_{d}\right)  \nabla fdx\\
& \quad+\displaystyle\frac{1}{\Upsilon}\int_{\Omega\cap B_{0,R_{0}}}\left(
\frac{1}{2}\Delta^{2}\varphi+2\mu\nabla\varphi\cdot\frac{x}{\left\vert
x\right\vert ^{4}}+\frac{\mu}{\left\vert x\right\vert ^{2}}\right)  \left\vert
f\right\vert ^{2}dx\\
& \quad-\displaystyle\frac{1}{2\Upsilon^{3}}\int_{\Omega\cap B_{0,R_{0}}%
}\left(  \varphi+\left\vert \nabla\varphi\right\vert ^{2}+\frac{1}{2}%
\nabla\varphi\cdot\left(  2\nabla^{2}\varphi+I_{d}\right)  \nabla
\varphi\right)  \left\vert f\right\vert ^{2}dx\text{ }%
\end{array}
\]
gives%
\[%
\begin{array}
[c]{ll}
& \quad\left\langle -\left(  \mathcal{S}^{\prime}+\left[  \mathcal{S}%
,\mathcal{A}\right]  \right)  f,f\right\rangle _{0}-\displaystyle\frac
{1}{\Upsilon}\left\langle -\mathcal{S}f,f\right\rangle _{0}\\
& =-\displaystyle\frac{1}{\Upsilon}\left(  -4a+1\right)  \int_{\Omega\cap
B_{0,R_{0}}}\left(  \left\vert \nabla f\right\vert ^{2}-\frac{\mu}{\left\vert
x\right\vert ^{2}}\left\vert f\right\vert ^{2}\right)  dx\\
& \quad-\displaystyle\frac{2bs}{\Upsilon}\left[  \displaystyle\int_{\Omega\cap
B_{0,R_{0}}}\left\vert x\right\vert ^{s-2}\left\vert \nabla f\right\vert
^{2}dx-\left(  2-s\right)  \displaystyle\int_{\Omega\cap B_{0,R_{0}}%
}\left\vert x\right\vert ^{s-4}\left\vert x\cdot\nabla f\right\vert
^{2}dx\right] \\
& \quad+\displaystyle\frac{1}{2\Upsilon}bs\left[  4\mu-\left(  2-s\right)
\left(  n+s-2\right)  \left(  n+s-4\right)  \right]  \displaystyle\int
_{\Omega\cap B_{0,R_{0}}}\left\vert x\right\vert ^{s-4}\left\vert f\right\vert
^{2}dx\\
& \quad+\displaystyle\frac{1}{2\Upsilon^{3}}c\displaystyle\int_{\Omega\cap
B_{0,R_{0}}}\left\vert f\right\vert ^{2}dx\\
& \quad+\displaystyle\frac{1}{2\Upsilon^{3}}a\left(  1-2a\right)  \left(
1-4a\right)  \displaystyle\int_{\Omega\cap B_{0,R_{0}}}\left\vert x\right\vert
^{2}\left\vert f\right\vert ^{2}dx\\
& \quad+\displaystyle\frac{1}{2\Upsilon^{3}}b\left(  -1+6as-4a^{2}%
s-4a^{2}s^{2}\right)  \displaystyle\int_{\Omega\cap B_{0,R_{0}}}\left\vert
x\right\vert ^{s}\left\vert f\right\vert ^{2}dx\\
& \quad-\displaystyle\frac{1}{2\Upsilon^{3}}\left(  bs\right)  ^{2}\left(
\frac{3}{2}+2a-4as\right)  \displaystyle\int_{\Omega\cap B_{0,R_{0}}%
}\left\vert x\right\vert ^{2s-2}\left\vert f\right\vert ^{2}dx\\
& \quad-\displaystyle\frac{1}{2\Upsilon^{3}}\left(  bs\right)  ^{3}\left(
s-1\right)  \displaystyle\int_{\Omega\cap B_{0,R_{0}}}\left\vert x\right\vert
^{3s-4}\left\vert f\right\vert ^{2}dx\text{ .}%
\end{array}
\]
When $n=3$, $a=\frac{1}{4}$, $b=\frac{1}{4}$, $c=\left(  \frac{1}{3}\right)
^{4}$ and $s=\frac{4}{3}$, since $R_{0}=\left(  \frac{4}{3}\right)  ^{3/2}$,
we have%
\[%
\begin{array}
[c]{ll}
& \quad\left\langle -\left(  \mathcal{S}^{\prime}+\left[  \mathcal{S}%
,\mathcal{A}\right]  \right)  f,f\right\rangle _{0}-\displaystyle\frac
{1}{T-t+\hbar}\left\langle -\mathcal{S}f,f\right\rangle _{0}\\
& \leq\displaystyle\frac{2bs}{\Upsilon}\left(  \mu-\frac{7}{2\cdot3^{3}%
}\right)  \displaystyle\int_{\Omega\cap B_{0,R_{0}}}\left\vert x\right\vert
^{-8/3}\left\vert f\right\vert ^{2}dx\\
& \leq0
\end{array}
\]
with our assumption on $\mu$. Since $\varphi\left(  x\right)  =W\left(
\left\vert x\right\vert \right)  $ with $W\left(  \rho\right)  =-\frac{1}%
{4}\rho^{2}+\frac{1}{4}\rho^{4/3}-\left(  \frac{1}{3}\right)  ^{4}$, we have
$W\left(  0\right)  =W\left(  1\right)  =-\left(  \frac{1}{3}\right)  ^{4}$,
$W^{\prime}\left(  \left(  2/3\right)  ^{3/2}\right)  =0$ and $\rho\mapsto
W\left(  \rho\right)  $ is strictly decreasing for $\rho\geq1$ and the
assumptions on $\varphi$ at Step \ref{step2.2.5} and Step \ref{step2.2.7} of
Section 2 hold by choosing $\omega_{0}=\left\{  x;r_{0}<\left\vert
x\right\vert <r\right\}  $ with $0<r_{0}<r<1$, $W\left(  r_{0}\right)
=W\left(  r\right)  \in\left(  -\left(  \frac{1}{3}\right)  ^{4},0\right)  $
and by taking $\ell>1$ sufficiently large.

When $n\geq4$, $a=\frac{1}{4}$, $b=\frac{1}{4}$, $c=b^{2}$ and $s=1$, we have
\[%
\begin{array}
[c]{ll}
& \quad\left\langle -\left(  \mathcal{S}^{\prime}+\left[  \mathcal{S}%
,\mathcal{A}\right]  \right)  f,f\right\rangle _{0}-\displaystyle\frac
{1}{T-t+\hbar}\left\langle -\mathcal{S}f,f\right\rangle _{0}\\
& \leq\displaystyle\frac{2bs}{\Upsilon}\left(  \mu-\frac{1}{4}\left(
n-1\right)  \left(  n-3\right)  \right)  \displaystyle\int_{\Omega\cap
B_{0,R_{0}}}\left\vert x\right\vert ^{-3}\left\vert f\right\vert ^{2}dx\\
& \leq0
\end{array}
\]
with our assumption on $\mu$. Since $\rho\mapsto-\frac{1}{4}\rho^{2}+\frac
{1}{4}\rho-\frac{1}{16}$ is a non positive function and is strictly decreasing
for $\rho\geq1$, the assumptions on $\varphi$ at Step \ref{step2.2.5} and Step
\ref{step2.2.7} of Section 2 hold by choosing $\omega_{0}=\left\{
x;r_{0}<\left\vert x\right\vert <r\right\}  $ with $0<r_{0}<r<1$, and by
taking $\ell>1$ sufficiently large.

\bigskip

One conclude that for any $u_{0}\in L^{2}\left(  \Omega\right)  $ and any
$T>0$,
\[
\left\Vert u\left(  \cdot,T\right)  \right\Vert _{L^{2}\left(  B_{0,R}\right)
}\leq\left(  ce^{\frac{K}{T}}\left\Vert u\left(  \cdot,T\right)  \right\Vert
_{L^{2}\left(  \omega_{0}\right)  }\right)  ^{\beta}\left\Vert u_{0}%
\right\Vert _{L^{2}\left(  \Omega\right)  }^{1-\beta}\text{ .}%
\]
Since $0\notin\omega_{0}\Subset\Omega$, we can replace $\omega_{0}$ by any
nonempty open subset $\widetilde{\omega}$ of $\Omega$ by propagation of
smallness. The treatment far from the point $0\in\Omega$ where the inverse
square potential have its singularities can be done in the same way than for
the heat equation with a potential in $L^{\infty}\left(  \Omega\times\left(
0,T\right)  \right)  $ (see \cite{PWZ}) and we also have
\[
\left\Vert u\left(  \cdot,T\right)  \right\Vert _{L^{2}\left(  \Omega
\left\backslash B_{0,R}\right.  \right)  }\leq\left(  ce^{\frac{K}{T}%
}\left\Vert u\left(  \cdot,T\right)  \right\Vert _{L^{2}\left(  \widetilde
{\omega}\right)  }\right)  ^{\beta}\left\Vert u_{0}\right\Vert _{L^{2}\left(
\Omega\right)  }^{1-\beta}\text{ .}%
\]
Finally, we will replace $\left\Vert u\left(  \cdot,T\right)  \right\Vert
_{L^{2}\left(  \widetilde{\omega}\right)  }$ by $\left\Vert u\left(
\cdot,T\right)  \right\Vert _{L^{1}\left(  \omega\right)  }$ thanks to Nash
inequality: Here $\widetilde{\omega}\Subset\omega$. Let $\phi\in C_{0}%
^{\infty}\left(  \omega\right)  $ be such that $0\leq\phi\leq1$ and $\phi=1$
on $\widetilde{\omega}$. Then we have
\[%
\begin{array}
[c]{ll}%
\left\Vert u\left(  \cdot,T\right)  \right\Vert _{L^{2}\left(  \widetilde
{\omega}\right)  } & \leq\left\Vert \phi u\left(  \cdot,T\right)  \right\Vert
_{L^{2}\left(  \omega\right)  }\\
& \leq\left(  e\left\Vert \phi u\left(  \cdot,T\right)  \right\Vert
_{L^{1}\left(  \omega\right)  }\right)  ^{\frac{2}{2+n}}\left(  \frac
{1}{2\sqrt{\pi}}\left\Vert \phi u\left(  \cdot,T\right)  \right\Vert
_{H_{0}^{1}\left(  \omega\right)  }\right)  ^{\frac{n}{2+n}}\\
& \leq\left(  e\left\Vert u\left(  \cdot,T\right)  \right\Vert _{L^{1}\left(
\omega\right)  }\right)  ^{\frac{2}{2+n}}\left(  \frac{1}{2\sqrt{\pi}%
}\left\Vert \phi u\left(  \cdot,T\right)  \right\Vert _{H_{0}^{1}\left(
\Omega\right)  }\right)  ^{\frac{n}{2+n}}\\
& \leq C_{n}\left\Vert u\left(  \cdot,T\right)  \right\Vert _{L^{1}\left(
\omega\right)  }^{\frac{2}{2+n}}\left\Vert u\left(  \cdot,T\right)
\right\Vert _{H_{0}^{1}\left(  \Omega\right)  }^{\frac{n}{2+n}}\\
& \leq C_{n}\left\Vert u\left(  \cdot,T\right)  \right\Vert _{L^{1}\left(
\omega\right)  }^{\frac{2}{2+n}}\left(  \frac{C}{\sqrt{T}}\left\Vert
u_{0}\right\Vert _{L^{2}\left(  \Omega\right)  }\right)  ^{\frac{n}{2+n}%
}\text{ .}%
\end{array}
\]
This completes the proof.

\bigskip

\bigskip

\section*{Appendix}

\label{appendix}

\bigskip

Let $H$ be a real Hilbert space endowed with an inner product $\left\langle
\cdot,\cdot\right\rangle $, and $A$ be a linear self-adjoint operator from
$D(A)$ into $H$, where $D(A)$ being the domain of $A$ is a subspace of $H$. We
assume that $A$ is an isomorphism from $D(A)$ (equipped with the graph norm)
onto $H$, that $A^{-1}$ is a linear compact operator in $H$ and that
$\left\langle Av,v\right\rangle >0$ $\forall v\in D(A)$, $v\neq0$. Introduce
the set $\left\{  \lambda_{j}\right\}  _{j=1}^{\infty}$ for the family of all
eigenvalues of $A$ so that%
\[
0<\lambda_{1}\leq\lambda_{2}\leq\cdot\cdot\leq\lambda_{m}\leq\lambda_{m+1}%
\leq\cdot\cdot\cdot\text{ and }\underset{j\rightarrow\infty}{\text{lim}%
}\lambda_{j}=\infty\text{ ,}%
\]
and let $\left\{  e_{j}\right\}  _{j=1}^{\infty}$ be the family of the
corresponding orthogonal normalized eigenfunctions: $Ae_{j}=\lambda_{j}e_{j}$,
$e_{j}\in D(A)$ and $\left\langle e_{j},e_{i}\right\rangle =\delta_{i,j}$.

\bigskip

By Lumer-Phillips theorem, $-A$ generates on $H$ a strongly continuous
semigroup $S:t\mapsto S(t)=e^{-tA}$. For any $t\geq0$ and any $u_{0}\in H$, we
have that $S\left(  t\right)  u_{0}=\sum_{j\geq1}\left\langle u_{0}%
,e_{j}\right\rangle e^{-\lambda_{j}t}e_{j}:=u\left(  \cdot,t\right)  $ and
$u\in C\left(  \left[  0,+\infty\right)  ;H\right)  \cap C^{1}\left(  \left(
0,+\infty\right)  ;D(A)\right)  $ is the unique solution of $\partial
_{t}u+Au=0$ with $u\left(  \cdot,0\right)  =u_{0}$.

\bigskip

Below, $H:=L^{2}\left(  \Omega\right)  $ where $\Omega$ is a bounded open set
of $\mathbb{R}^{n}$.

\bigskip

\begin{appen}
Let $\omega$ be a nonempty open subset of $\Omega$. Let $p\in\left[
1,2\right]  $, $\beta\in\left(  0,1\right)  $ and $c$, $K$, $\gamma>0$.
Suppose that for any $u_{0}\in L^{2}\left(  \Omega\right)  $ and any $T>0$,
\[
\left\Vert u\left(  \cdot,T\right)  \right\Vert _{L^{2}\left(  \Omega\right)
}\leq\left(  ce^{\frac{1}{T^{\gamma}}K}\left\Vert u\left(  \cdot,T\right)
\right\Vert _{L^{p}\left(  \omega\right)  }\right)  ^{\beta}\left\Vert
u_{0}\right\Vert _{L^{2}\left(  \Omega\right)  }^{1-\beta}\text{ .}%
\]
Then for any $\left(  a_{j}\right)  _{j\geq1}\in\mathbb{R}$ and any
$\lambda>0$, one has%
\[
\sqrt{\sum\limits_{\lambda_{j}\leq\lambda}\left\vert a_{j}\right\vert ^{2}%
}\leq ce^{\lambda^{\frac{\gamma}{1+\gamma}}2\left(  \frac{1-\beta}{\beta
}\right)  ^{\frac{\gamma}{1+\gamma}}K^{\frac{1}{1+\gamma}}}\left\Vert
\sum\limits_{\lambda_{j}\leq\lambda}a_{j}e_{j}\right\Vert _{L^{p}\left(
\omega\right)  }\text{ .}%
\]

\end{appen}

\bigskip

Indeed, we choose $u_{0}=\sum\limits_{\lambda_{j}\leq\lambda}a_{j}%
e^{\lambda_{j}T}e_{j}$ and apply
\[
\left\Vert u\left(  \cdot,T\right)  \right\Vert _{L^{2}\left(  \Omega\right)
}\leq\left(  ce^{K\frac{1}{T^{\gamma}}}\left\Vert u\left(  \cdot,T\right)
\right\Vert _{L^{p}\left(  \omega\right)  }\right)  ^{\beta}\left\Vert
u_{0}\right\Vert _{L^{2}\left(  \Omega\right)  }^{1-\beta}\text{ ,}%
\]
to get
\[%
\begin{array}
[c]{ll}%
\sqrt{\displaystyle\sum\limits_{\lambda_{j}\leq\lambda}\left\vert
a_{j}\right\vert ^{2}} & \leq\left(  ce^{K\frac{1}{T^{\gamma}}}\left\Vert
\sum\limits_{\lambda_{j}\leq\lambda}a_{j}e_{j}\right\Vert _{L^{p}\left(
\omega\right)  }\right)  ^{\beta}\left(  \sqrt{\displaystyle\sum
\limits_{\lambda_{j}\leq\lambda}\left\vert a_{j}\right\vert ^{2}%
e^{2\lambda_{j}T}}\right)  ^{1-\beta}\\
& \leq\left(  ce^{K\frac{1}{T^{\gamma}}}\left\Vert \sum\limits_{\lambda
_{j}\leq\lambda}a_{j}e_{j}\right\Vert _{L^{p}\left(  \omega\right)  }\right)
^{\beta}e^{\lambda T\left(  1-\beta\right)  }\left(  \sqrt{\displaystyle\sum
\limits_{\lambda_{j}\leq\lambda}\left\vert a_{j}\right\vert ^{2}}\right)
^{1-\beta}\text{ .}%
\end{array}
\text{ }%
\]
Therefore,
\[
\sqrt{\sum\limits_{\lambda_{j}\leq\lambda}\left\vert a_{j}\right\vert ^{2}%
}\leq ce^{K\frac{1}{T^{\gamma}}+\lambda T\left(  \frac{1-\beta}{\beta}\right)
}\left\Vert \sum\limits_{\lambda_{j}\leq\lambda}a_{j}e_{j}\right\Vert
_{L^{p}\left(  \omega\right)  }\text{ .}%
\]
We conclude by choosing
\[
T=\left[  \left(  \frac{\beta}{1-\beta}\right)  \frac{K}{\lambda}\right]
^{\frac{1}{1+\gamma}}\text{ .}%
\]

\bigskip

Remark .- Conversely, suppose that there are constants $p\in\left[
1,2\right]  $ and $D_{1}$, $D_{2}$, $\gamma>0$ such that any $\left(
a_{j}\right)  _{j\geq0}\in\ell^{2}$ and any $\lambda>\lambda_{1}$,
\[
\sqrt{\sum\limits_{\lambda_{j}\leq\lambda}\left\vert a_{j}\right\vert ^{2}%
}\leq D_{1}e^{\lambda^{\frac{\gamma}{1+\gamma}}D_{2}}\left\Vert \sum
\limits_{\lambda_{j}\leq\lambda}a_{j}e_{j}\right\Vert _{L^{p}\left(
\omega\right)  }\text{ .}%
\]
Then for any $\beta\in\left(  0,1\right)  $ and any $T>0$,%
\[
\left\Vert u\left(  \cdot,T\right)  \right\Vert _{L^{2}\left(  \Omega\right)
}\leq D_{3}e^{\frac{1}{T^{\gamma}}D_{4}}\left\Vert u\left(  \cdot,T\right)
\right\Vert _{L^{p}\left(  \omega\right)  }^{\beta}\left\Vert u\left(
\cdot,0\right)  \right\Vert _{L^{2}\left(  \Omega\right)  }^{1-\beta}%
\]
with
\[
D_{3}=2\left(  1+\text{max}\left(  1,\left\vert \omega\right\vert ^{\frac
{1}{p}-\frac{1}{2}}\right)  D_{1}\right)  \text{ and }D_{4}=\left(
D_{2}\right)  ^{1+\gamma}\frac{1}{\left(  1-\beta\right)  ^{\gamma}}\text{ .}%
\]

\bigskip

Indeed, let $\alpha:=\frac{\gamma}{1+\gamma}$ and $u_{0}:=u\left(
\cdot,0\right)  =\sum\limits_{j\geq1}a_{j}e_{j}$. First, we have
\[
\left\Vert u\left(  \cdot,T\right)  \right\Vert _{L^{2}\left(  \Omega\right)
}^{2}=\sum_{j\geq1}\left\vert a_{j}e^{-\lambda_{j}T}\right\vert ^{2}%
=\sum_{\lambda_{j}\leq\lambda}\left\vert a_{j}e^{-\lambda_{j}T}\right\vert
^{2}+\sum_{\lambda_{j}>\lambda}\left\vert a_{j}e^{-\lambda_{j}T}\right\vert
^{2}%
\]
and%
\[
\left\Vert u\left(  \cdot,T\right)  \right\Vert _{L^{2}\left(  \Omega\right)
}\leq\sqrt{\sum_{\lambda_{j}\leq\lambda}\left\vert a_{j}e^{-\lambda_{j}%
T}\right\vert ^{2}}+\sqrt{\sum_{\lambda_{j}>\lambda}\left\vert a_{j}%
e^{-\lambda_{j}T}\right\vert ^{2}}\text{ .}%
\]
Next, we apply the estimate on the sum of eigenfunctions with $a_{j}$ replaced
by $a_{j}e^{-\lambda_{j}T}$ in order to get%
\[%
\begin{array}
[c]{ll}%
\left\Vert u\left(  \cdot,T\right)  \right\Vert _{L^{2}\left(  \Omega\right)
} & \leq D_{1}e^{D_{2}\lambda^{\alpha}}\left\Vert \displaystyle\sum
_{\lambda_{j}\leq\lambda}a_{j}e^{-\lambda_{j}T}e_{j}\right\Vert _{L^{p}\left(
\omega\right)  }+\sqrt{\displaystyle\sum_{\lambda_{j}>\lambda}\left\vert
a_{j}e^{-\lambda_{j}T}\right\vert ^{2}}\\
& \leq D_{1}e^{D_{2}\lambda^{\alpha}}\left\Vert \displaystyle\sum_{\lambda
_{j}\leq\lambda}a_{j}e^{-\lambda_{j}T}e_{j}+\displaystyle\sum_{\lambda
_{j}>\lambda}a_{j}e^{-\lambda_{j}T}e_{j}\right\Vert _{L^{p}\left(
\omega\right)  }\\
& \quad+D_{1}e^{D_{2}\lambda^{\alpha}}\left\Vert \displaystyle\sum
_{\lambda_{j}>\lambda}a_{j}e^{-\lambda_{j}T}e_{j}\right\Vert _{L^{p}\left(
\omega\right)  }+\sqrt{\displaystyle\sum_{\lambda_{j}>\lambda}\left\vert
a_{j}e^{-\lambda_{j}T}\right\vert ^{2}}\\
& \leq D_{1}e^{D_{2}\lambda^{\alpha}}\left\Vert u\left(  \cdot,T\right)
\right\Vert _{L^{p}\left(  \omega\right)  }+\left(  1+\left\vert
\omega\right\vert ^{\frac{1}{p}-\frac{1}{2}}D_{1}\right)  e^{D_{2}%
\lambda^{\alpha}}e^{-\lambda T}\sqrt{\displaystyle\sum_{\lambda_{j}>\lambda
}\left\vert a_{j}\right\vert ^{2}}\\
& \leq D_{1}e^{D_{2}\lambda^{\alpha}}\left\Vert u\left(  \cdot,T\right)
\right\Vert _{L^{p}\left(  \omega\right)  }+\left(  1+\left\vert
\omega\right\vert ^{\frac{1}{p}-\frac{1}{2}}D_{1}\right)  e^{D_{2}%
\lambda^{\alpha}}e^{-\lambda T}\left\Vert u\left(  \cdot,0\right)  \right\Vert
_{L^{2}\left(  \Omega\right)  }\text{ .}%
\end{array}
\]
Now, by Young inequality, for any $\epsilon>0$,
\[
D_{2}\lambda^{\alpha}=\frac{D_{2}}{\left(  \epsilon T\right)  ^{\alpha}%
}\left(  \epsilon\lambda T\right)  ^{\alpha}\leq\epsilon\lambda T+\left(
\frac{D_{2}}{\left(  \epsilon T\right)  ^{\alpha}}\right)  ^{\frac{1}%
{1-\alpha}}\text{ .}%
\]
Therefore,
\[%
\begin{array}
[c]{ll}%
\left\Vert u\left(  \cdot,T\right)  \right\Vert _{L^{2}\left(  \Omega\right)
} & \leq\left(  1+\text{max}\left(  1,\left\vert \omega\right\vert ^{\frac
{1}{p}-\frac{1}{2}}\right)  D_{1}\right)  e^{\left(  \frac{D_{2}}{\left(
\epsilon T\right)  ^{\alpha}}\right)  ^{\frac{1}{1-\alpha}}}\\
& \quad\times\left(  e^{\epsilon\lambda T}\left\Vert u\left(  \cdot,T\right)
\right\Vert _{L^{p}\left(  \omega\right)  }+e^{\left(  \epsilon-1\right)
\lambda T}\left\Vert u\left(  \cdot,0\right)  \right\Vert _{L^{2}\left(
\Omega\right)  }\right)  \text{ .}%
\end{array}
\]
Choosing $0<\epsilon<1$ and optimizing with respect to $\lambda T$ by taking
\[
\lambda T=\text{ln}\left(  \frac{\left\Vert u\left(  \cdot,0\right)
\right\Vert _{L^{2}\left(  \Omega\right)  }}{\left\Vert u\left(
\cdot,T\right)  \right\Vert _{L^{p}\left(  \omega\right)  }}\right)  \text{ ,}%
\]
yield
\[%
\begin{array}
[c]{ll}%
\left\Vert u\left(  \cdot,T\right)  \right\Vert _{L^{2}\left(  \Omega\right)
} & \leq\left(  1+\text{max}\left(  1,\left\vert \omega\right\vert ^{\frac
{1}{p}-\frac{1}{2}}\right)  D_{1}\right)  e^{\left(  \frac{D_{2}}{\left(
\epsilon T\right)  ^{\alpha}}\right)  ^{\frac{1}{1-\alpha}}}\\
& \quad\times\left(  2\left[  \frac{\left\Vert u\left(  \cdot,0\right)
\right\Vert _{L^{2}\left(  \Omega\right)  }}{\left\Vert u\left(
\cdot,T\right)  \right\Vert _{L^{p}\left(  \omega\right)  }}\right]
^{\epsilon}\left\Vert u\left(  \cdot,T\right)  \right\Vert _{L^{p}\left(
\omega\right)  }\right)  \text{ .}%
\end{array}
\]
Setting $\beta=1-\epsilon$, we finally have the desired observation estimate
at one time.

\bigskip

\bigskip

\bigskip

\bigskip

\bigskip

\bigskip

\bigskip

\bigskip

\bigskip
\end{document}